\documentclass[11 pt]{amsart}

\usepackage{hyperref}
\usepackage{amssymb, amsmath, amsthm, amsfonts}
\usepackage{verbatim}
\usepackage{bbm}
\usepackage{enumerate}
\usepackage{dsfont}
\usepackage{upgreek}
\usepackage[mathscr]{eucal}
\usepackage{tikz}
\usepackage[normalem]{ulem}
\usepackage{kotex}

\usepackage{enumitem}

\usepackage{textcomp}

\usepackage[backgroundcolor=black!5,linecolor=black]{todonotes}
\usepackage{comment}

\usepackage[explicit,indentafter]{titlesec}
\titleformat{\section}{\centering\normalfont\scshape}{\thesection.}{.5em}{#1}
\titleformat{\subsection}[runin]{\normalfont\itshape}{\textnormal{\thesubsection.}}{.5em}{#1.}
\titleformat{\subsubsection}[runin]{\normalfont\itshape}{\thesubsubsection.}{.5em}{#1.}
\titlespacing{\section}{0em}{1em}{0.5em}
\titlespacing{\subsection}{0em}{.5em}{0.5em}

\usepackage{stackengine}
\newcommand{\ubar}[1]{\stackunder[1.2pt]{$#1$}{\rule{.9ex}{.075ex}}}

\newcommand{\ox}{{\overline x}}
\newcommand{\oy}{{\overline y}}

\newcommand{\ux}{{\underline x}}
\newcommand{\uy}{{\underline y}}

\newcommand{\om}{\omega}
\usepackage{color}
\definecolor{gray}{gray}{0.5}

\newcommand{\cmt}[1]{}

\newcommand{\vertiii}[1]{{\left\vert\kern-0.25ex\left\vert\kern-0.25ex\left\vert #1 
    \right\vert\kern-0.25ex\right\vert\kern-0.25ex\right\vert}}
    
    \newcommand{\detail}[1]{}

	\newcommand{\drawQbg}{
		\fill (Q1) node [left] {$Q_1$} circle [radius=.02em];
		\fill (Q2) node [above left] {$Q_{2}$} circle [radius=\ptsize];
		\fill (Q3) node [right] {$Q_{3}$} circle [radius=\ptsize];
		\fill (Q4g) node [below right] {$Q_{4}$} circle [radius=\ptsize];
		\fill [color=\QbgFillColor,opacity=\QbgFillOpacity] (Q1) -- (Q2) -- (Q3) -- (Q4g) -- cycle;
		\draw [opacity=\QbgLineOpacity] (Q1) -- (Q2) -- (Q3);
		\draw [opacity=\QbgLineOpacity] (Q4g) -- (Q1);
		\if\QbgDrawCritSeg1
			\draw [style=\QbgCritSegStyle,opacity=\QbgCritSegOpacity] (Q3) -- (Q4g);
		\fi
	}

\def\lc{\lesssim}

\NewDocumentEnvironment{amatrix}{>{\SplitArgument{1}{|}}m}
 {\left(\makeamatrix#1}
 {\end{array}\right)}
\NewDocumentCommand{\makeamatrix}{mm}{%
  \IfNoValueTF{#2}
    {\begin{array}{@{}*{#1}{c}@{}}}
    {\begin{array}{@{}*{#1}{c}|*{#2}{c}@{}}}%
}

\makeatletter
\def\widebreve{\mathpalette\wide@breve}
\def\wide@breve#1#2{\sbox\z@{$#1#2$}%
     \mathop{\vbox{\m@th\ialign{##\crcr
\kern0.08em\brevefill#1{0.8\wd\z@}\crcr\noalign{\nointerlineskip}%
                    $\hss#1#2\hss$\crcr}}}\limits}
\def\brevefill#1#2{$\m@th\sbox\tw@{$#1($}%
  \hss\resizebox{#2}{\wd\tw@}{\rotatebox[origin=c]{90}{\upshape(}}\hss$}
\makeatletter

\def\eps{\varepsilon}
\def\bbone{{\mathbbm 1}}

\newcommand{\wt}{\widetilde}

\newcommand{\Be}{\begin{equation}}
\newcommand{\Ee}{\end{equation}}

\newcommand{\Bm}{\begin{multline}}
\newcommand{\Em}{\end{multline}}

\def\intslash{\rlap{\kern  .32em $\mspace {.5mu}\backslash$ }\int}
\def\qsl{{\rlap{\kern  .32em $\mspace {.5mu}\backslash$ }\int_{Q_x}}}

\def\lc{\lesssim}

\def\emph#1{{\it #1 }}

\def\rank{{\mathrm{rank}}}

\def\dist{{\mathrm{dist}}}

\def\supp{{\mathrm{supp}}}

\def\inn#1#2{\langle#1,#2\rangle}

\def\eps{\varepsilon}

\def\la{\lambda}             \def\La{\Lambda}

\def\om{\omega}

\def\fG{{\mathfrak {G}}}

\def\fI{{\mathfrak {I}}}

\def\fK{{\mathfrak {K}}}

\def\fM{{\mathfrak {M}}}

\def\fe{{\mathfrak {e}}}

\def\fg{{\mathfrak {g}}}

\def\fn{{\mathfrak {n}}}

\def\fq{{\mathfrak {q}}}
\def\fr{{\mathfrak {r}}}

\def\fu{{\mathfrak {u}}}

\def\fz{{\mathfrak {z}}}

\def\bbR{{\mathbb {R}}}
\def\bbS{{\mathbb {S}}}

\def\cA{{\mathcal {A}}}
\def\cB{{\mathcal {B}}}
\def\cC{{\mathcal {C}}}

\def\cF{{\mathcal {F}}}
\def\cG{{\mathcal {G}}}

\def\cK{{\mathcal {K}}}
\def\cL{{\mathcal {L}}}

\def\cN{{\mathcal {N}}}

\def\cP{{\mathcal {P}}}
\def\cQ{{\mathcal {Q}}}
\def\cR{{\mathcal {R}}}
\def\cS{{\mathcal {S}}}
\def\cT{{\mathcal {T}}}
\def\cU{{\mathcal {U}}}

\def\cW{{\mathcal {W}}}

\def\emph#1{{\it #1}}
\def\textbf#1{{\bf #1}}

\def\beq{\begin{equation}}
\def\endeq{\end{equation}}

\def\bs{\begin{split}}
\def\es{\end{split}}

\theoremstyle{plain}

\newtheorem{thm}{Theorem}[section]
\newtheorem{prop}[thm]{Proposition}

\newtheorem{lem}[thm]{Lemma}

\newtheorem*{thm*}{Theorem}
\newtheorem*{conj*}{Conjecture}
\newtheorem*{openproblem*}{Open Problem}

\theoremstyle{remark}
\newtheorem{rem}[thm]{Remark}

\numberwithin{equation}{section}

\definecolor{jrcol}{rgb}{0,0,1.} % Blue
\definecolor{ascol}{rgb}{1.,0,0} % Red

\def\R{\mathbb{R}}
\def\C{\mathbb{C}}
\def\Z{\mathbb{Z}}
\def\N{\mathbb{N}}

\begin{document}
\title
[Perturbed  Nevo--Thangavelu maximal functions]
{Perturbed Nevo–Thangavelu maximal functions on two-step nilpotent groups}
%{Notes on tilted averages}
%{The Nevo--Thangavelu maximal function \\on  two-step nilpotent Lie groups}
\author[Jaehyeon Ryu]{Jaehyeon Ryu} 

\address{Jaehyeon Ryu: Department of Mathematics, Ewha Womans University, Seoul 03760, Korea}
\email{jhryu67@ewha.ac.kr}

%\address{Andreas Seeger: Department of Mathematics, University of Wisconsin, 480 Lincoln Drive, Madison, WI, 53706, USA.}
%\email{seeger@math.wisc.edu}

%\begin{thanks}
    %{\it Acknowledgement.} 
%    {J.R.  was supported in part by KIAS Individual Grant MG087001.   
%A.S.  was supported in part by NSF Grant DMS-2348797.
%} 
%\end{thanks}

%\date{\today}


\begin{abstract}
Consider a step-two nilpotent Lie group $G$ identified with $\R^d\times \R^m$ in exponential coordinates. We study maximal operators associated with parabolic dilates of smooth convex hypersurfaces with nonvanishing Gaussian curvature in a tilted $d$-dimensional subspace of $G$. We focus on the low-rank regime, where the rank is determined by the group law and the tilt. When $m=1$, we determine the optimal range of $p$ for which the maximal operator is bounded on $L^p(G)$. For $m\ge2$, we characterize, in terms of a {\L}ojasiewicz exponent, when the known necessary condition on $p$ is also sufficient.
\end{abstract}

\maketitle

\section{Introduction}\label{sec:intro} 
We investigate $L^p$ bounds for spherical maximal operators associated
with linearly tilted subspaces of a two-step nilpotent group. Let $G$ be a simply connected finite-dimensional two-step nilpotent Lie group. The Lie algebra $\fg$ admits the direct sum decomposition $\fg = \fg_1\oplus\mathfrak g_2$ where $\dim\mathfrak g_1=d$, $\dim\mathfrak g_2=m$, and $[\mathfrak g,\mathfrak g]\subseteq\mathfrak g_2 \subseteq Z(\mathfrak g)$. Here $Z(\mathfrak g)$ denotes the center of $\mathfrak g$. We refer to $\fg_1$ and $\fg_2$ as the horizontal and vertical subspaces, respectively.

Using the exponential coordinates on $G$ and fixed bases in the two summands, we write points of $G$ as $x=(\ux,\ox)\in\R^d\times\R^m$. The Baker-Campbell-Hausdorff formula gives the group law
\[
    (\ux,\ox)\cdot(\uy,\oy)
        =\bigl(
            \ux+\uy,\,
            \ox+\oy+\ux^\intercal\vec J\uy
          \bigr),
\]
where $\ux^\intercal\vec J\uy = (\ux^\intercal J_1\uy,\dots, \ux^\intercal J_m\uy)$, and $J_1,\dots,J_m$ are real skew-symmetric $d\times d$ matrices representing one half of the Lie bracket $\frac12[X,Y]$, $X,Y\in \fg$. The associated automorphic dilation is given by
\[
    \delta_t(\ux,\ox)=(t\ux,t^2\ox),
    \qquad t>0.
\]

For a linear map $\Lambda:\R^d\to\R^m$, consider a $d$-dimensional plane $H = \{(\om,\La\om):\om \in \R^d\}$.
The choice $\Lambda=0$ gives the horizontal plane, while nonzero $\Lambda$ introduces a linear tilt. Let $\Omega\subset \R^d$ denote a bounded open convex domain containing $0$. Throughout this paper, we assume that the boundary $\Sigma = \partial\Omega$ is smooth and has nonvanishing Gaussian curvature everywhere. Let $\mu$ be the normalized surface measure on $\Sigma$. We define a measure $\mu^\La$ by
\[
\inn{f}{\mu^\La} = \int_\Sigma f(\om,\La\om)d\mu(\om),\quad f\in \cS(\R^{d+m})
\]
where $\cS(\R^{d+m})$ denotes the Schwartz space. By $\Sigma^\La$ we denote the embedded hypersurface $\{(\om,\La\om):\om\in \Sigma\}$ in $H$.
For each $t>0$, we define the dilate $\mu_t^\La$ by $\inn{f}{\mu_t^\La} = \int f(t\uy,t^2\oy) d\mu^\La(\uy,\oy)$. Let $A_t f$ be the averages of $f\in \cS(\R^{d+m})$ over the dilates of $\Sigma^\La$ defined by the non-commutative convolutions $f\ast \mu_t^\La$, that is,
\[
A_t f(x) = f\ast \mu_t^\La(x) = \int_\Sigma f(\ux -t\om, \ox - t^2\La\om -t\ux^\intercal \vec J \om) d\mu(\om).
\]
The associated maximal operator $\mathfrak M$ is defined by
\[
    \fM f(\ux,\ox)
        =\sup_{t>0}
            \bigl|A_t f(\ux,\ox)\bigr|
\]
initially for $f\in\mathcal S(\R^{d+m})$. We will study the problem of determining conditions on the group $G$ and the tilt $\Lambda$ for the $L^p(G)$ boundedness of $\fM$.

This problem was first considered by Nevo and Thangavelu \cite{NevoThangavelu1997} in the untilted case $\La=0$. The boundedness of $\fM$ was subsequently studied in \cite{MuellerSeeger2004, NarayananThangavelu2004}. See also \cite{JK26} for a restricted weak type bound at the endpoint. 
More recently, in \cite{RyuSeeger}, the authors proved that, for $d\ge3$ and $\La = 0$, $\fM$ is bounded on $L^p$ if and only if $p>\frac{d}{d-1}$ on an arbitrary two-step nilpotent group $G$. For $d=2$, results on $L^p$ boundedness are available only for Heisenberg radial functions (see \cite{BeltranGuoHickmanSeeger, LeeLee}), and the problem for general functions remains open.

The tilted model was introduced in M\"uller--Seeger \cite{MuellerSeeger2004} as a perturbed analog of the Nevo--Thangavelu maximal operator. In that work, the authors obtained $L^p$ bounds for $\fM$ for a partial range of $p$ when $G$ is a M\'etivier group and $\La\ne0$. It was subsequently proved in \cite{RoosSeegerSrivastava-imrn} that, when $G$ is a M\'etivier group with $d\ge3$, the condition $p>\frac{d}{d-1}$ is sufficient for $L^p$ boundedness under a smallness assumption on $\La$. When $G$ is a Heisenberg group with $d\ge3$, the authors of \cite{AndersonCladekPramanikSeeger} established that $\fM$ is bounded on $L^p$ if and only if $p>\frac{d}{d-1}$, without any smallness assumption on $\La$. These results show that, under a nondegeneracy condition on $G$, the boundedness property of $\fM$ is stable under small perturbations by $\La$.

In this paper, we will study perturbed averages on a general two-step group $G$ when $\La\ne0$. In recent work \cite{RS26}, the authors showed, in a general setting for $\Sigma$, that stability can fail in the absence of the nondegeneracy condition on $G$. We state below the main result of \cite{RS26} pertaining to the case where $\Sigma$ is of codimension $1$. We first introduce the notation. Let $\{\la_k\}_{k=1}^m\subset\R^d$ denote the columns of
$\La^\intercal$.
For each $1\le k\le m$, define the $(d+1)\times d$ matrix
\begin{align}\label{def:Si}
    S_k=\begin{pmatrix}
        J_k\\
        \la_k^\intercal
    \end{pmatrix}.
\end{align}
For $\theta=(\theta_1,\dots,\theta_m)\in\R^m$, set
\[
    S^\theta=\sum_{k=1}^m\theta_k S_k, \quad
    r(\theta)=\rank(S^\theta).
\]
For $\La\ne0$, define
\begin{align}\label{i:defr}
    r=\min\bigl\{ r(\theta): \theta\in\R^m,\ \La^\intercal\theta\ne0 \bigr\}.
\end{align}

\begin{thm}[\cite{RS26}]\label{thm:fail-RS26}
Assume $d\ge 2$, $\La\neq 0$, and $1\le r\le d-1$.
    \begin{itemize}[topsep = 2pt, itemsep = 2pt]
        \item [(i)] If $2\le r\le d-1$ and $\fM$ extends to a bounded operator on $L^p(G)$ then $p>\frac{r+1}r$.
        \item [(ii)] If $r=1$ and $\fM$ extends to a bounded operator on $L^p(G)$ then $p = \infty$.
    \end{itemize}
\end{thm}

\subsection*{Main results} Our goal is to determine whether the necessary condition on $p$ in Theorem \ref{thm:fail-RS26} is also sufficient when $2\le r\le d-2$. We refer to this range as the low-rank case. For $d-1\le r\le d$, Stein's example \cite{SteinPNAS1976} gives the necessary condition $p>\frac{d}{d-1}$, as also noted in \cite{LiuYan, RyuSeeger}. Thus, for $r\ge2$, the low-rank case is precisely the regime in which Theorem \ref{thm:fail-RS26} yields a stronger necessary condition and hence a failure of stability for $\fM$. We present our main results separately for $m=1$ and $m\ge2$.

We begin with the case $m=1$. In this setting, $\La$ equals the row vector $\la_1^\intercal$ and $r$ denotes the rank of $S_1$. We establish that the necessary condition in Theorem \ref{thm:fail-RS26} is also sufficient.

\begin{thm}\label{thm:main-m=1}
Assume $d\ge 4$, $m=1$, $\La\neq 0$, and $2\le r\le d-2$. Then $\fM$ extends to a bounded operator on $L^p(G)$ if and only if $\frac{r+1}{r}<p\le\infty$.
\end{thm}

We now turn to the case $m\ge 2$. When $m\ge 2$, the boundedness property of $\fM$ becomes substantially more complicated. Somewhat surprisingly, the algebraic structure of $G$ and $\La$ can give rise to a new range of $p$ in which $\fM$ is unbounded on $L^p(G)$, beyond the necessary range in Theorem \ref{thm:fail-RS26}. By the definition of $r$, the product of the $r$ largest singular values of $S^\theta$ is positive whenever $\Lambda^\top\theta\ne0$. Whether the necessary condition in Theorem \ref{thm:fail-RS26} is sufficient depends on how rapidly this product can decay as $\theta$ approaches $\ker(\Lambda^\top)$. We therefore seek a lower bound for the product in terms of a power of $\operatorname{dist}(\theta,\ker(\Lambda^\top))$. This leads us to the {\L}ojasiewicz exponent, whose definition we now recall.

We follow the exposition and notation of Basu and Mohammad-Nezhad \cite{BM24}. Throughout the paper, let
\[
    \cB=\{\theta\in\R^m:|\theta|\le C_\cB\},
\]
where $C_\cB>0$ is a sufficiently large constant. The set $\cB$ is clearly closed, bounded, and semi-algebraic. Suppose that $f,g:\cB\to\R$ are continuous semi-algebraic functions satisfying $f^{-1}(0)\subseteq g^{-1}(0)$.
Then there exist a constant $C>0$ and a rational number $\rho>0$ such that
\[
    |g(x)|^\rho\le C\,|f(x)|, \qquad x\in\cB.
\]
This inequality is called the {\L}ojasiewicz inequality. Moreover, the infimum of the admissible exponents $\rho$ is called the {\L}ojasiewicz exponent and is denoted by $\cL(f,g\,|\,\cB)$. The inequality originated in the work of {\L}ojasiewicz \cite{Loj59}. For a proof of the semi-algebraic version stated above, see \cite{BCR98}.

For $\theta\in\R^m$, write
\[
    s_1(\theta)\ge s_2(\theta)\ge\cdots\ge s_d(\theta)\ge0
\]
for the singular values of $S^\theta$, and consider the product of the $r$ largest singular values,
\[
    \cP_r(\theta)=\prod_{k=1}^r s_k(\theta).
\]
Let $\Pi_\La$ denote the orthogonal projection of $\R^m$ onto the orthogonal complement of $\ker(\La^\intercal)$, and set
\[
    d_\perp(\theta)=|\Pi_\La\theta| =\operatorname{dist}(\theta,\ker(\La^\intercal)).
\]

As mentioned before, we are interested in lower bounds for $\cP_r(\theta)$ in terms of $d_\perp(\theta)$. To formulate such a bound, we use the fact that the square of the product is comparable to the absolute value of a polynomial in $\theta$. We consider the characteristic polynomial
\[
    p_\theta(\lambda) =\det\bigl((S^\theta)^\intercal S^\theta-\lambda I_d\bigr),
\]
and let $\cC_r(\theta)$ denote the coefficient of $\lambda^{d-r}$ in $p_\theta(\lambda)$. Then $\cC_r$ is a homogeneous polynomial of degree $2r$ in $\theta$, and $|\cC_r(\theta)|$ is comparable to $\cP_r(\theta)^2$. See Lemma~\ref{lem:claim3} for details. In particular, $\cC_r^{-1}(0)\subset d_\perp^{-1}(0)$. Since $\cC_r$ and $d_\perp$ are continuous semi-algebraic functions, the exponent $\cL(\cC_r,d_\perp\,|\,\cB)$ is well defined.

We first state the result regarding a sufficient condition. The following demonstrates that the necessary condition in Theorem~\ref{thm:fail-RS26} is also sufficient for $L^p(G)$ boundedness of $\fM$ when $\cL(\cC_r,d_\perp\,|\,\cB)=2r$.

\begin{thm}\label{thm:main-hsu}
Let $G$ be a two-step nilpotent Lie group.
Assume $d\ge4$, $m\ge2$, $\La\ne0$, and $2\le r\le d-2$.
If
\[
    \cL(\cC_r,d_\perp\,|\,\cB)=2r,
\]
then $\fM$ is bounded on $L^p(G)$ if and only if $\frac{r+1}{r}<p\le\infty$.
\end{thm}

Fix $\theta$ with $d_\perp(\theta)>0$. Since $C_r(a\theta)=a^{2r}C_r(\theta)$ and $d_\perp(a\theta)=a\,d_\perp(\theta)$ for $a>0$, letting $a\to 0$ in the {\L}ojasiewicz inequality shows that every admissible exponent is at least $2r$. Theorem \ref{thm:main-hsu} treats the case of equality.

We next present a new necessary condition. We point out that although $\cC_r$ is homogeneous of degree $2r$, the {\L}ojasiewicz exponent $\cL(\cC_r,d_\perp\,|\,\cB)$ need not equal $2r$. Indeed, in Section \ref{sec:necessary} we construct examples of $G$ and $\La$ for which the exponent is considerably larger than $2r$. Our next theorem shows that whenever $\cL(\cC_r,d_\perp\,|\,\cB)>2r$, an additional restriction on $p$ is necessary.

\begin{thm}\label{thm:main-hne}
Assume $d\ge4$, $m\ge2$, $\La\ne0$, and $2\le r\le d-2$.
Suppose that $\cL(\cC_r,d_\perp\,|\,\cB)>2r$. Then there exists $p_\circ=p_\circ(G,\La,\Sigma)>\frac{r+1}{r}$ such that $\fM$ is unbounded on $L^p(G)$ whenever $\frac{r+1}{r}<p<p_\circ$.
\end{thm}

\subsection*{Remarks} \hfill

(a) When $m=2$, the exponent condition $\cL(\cC_r,d_\perp\,|\,\cB)=2r$ in Theorem \ref{thm:main-hsu} holds (see Remark \ref{rem:loja-2r}). Thus, under the other assumptions of the theorem, $\mathfrak M$ is bounded on $L^p(G)$ if and only if $p>(r+1)/r$.

(b) For $m=1$, our argument also applies when $d-1\le r\le d$. For $d\ge4$, it yields $L^p(G)$ boundedness of $\fM$ in the range $\frac{d-1}{d-2}<p\le\infty$. See Remark \ref{rem:lpest-hrank} for the corresponding local estimate. We conjecture that boundedness also holds in the remaining range $\frac{d}{d-1}<p\le\frac{d-1}{d-2}$. Establishing this conjecture appears to require a more refined analysis of the oscillatory integral operator discussed in Section \ref{sec:l2estspf}, without freezing one of the variables as in our argument. This would make the analysis considerably more involved.

(c) If $m=1$ and $r\le d-3$, we can also obtain a restricted weak type estimate at the endpoint $p=\frac{r+1}{r}$ for the local maximal operator $\fM_I f(x)=\sup_{t\in I}|A_t f(x)|$, where $I\subset(0,\infty)$ is a compact interval (see Remark \ref{rem:endpoint}). For $m\ge2$, however, our argument does not yield the corresponding endpoint estimate.

(d) The main result of \cite{RS26} addresses the more general setting in which $\Sigma$ is a $C^1$ submanifold of arbitrary dimension. For general submanifolds, the condition $\frac{r+1}{r}<p\le\infty$ is not sufficient for boundedness without additional geometric assumptions. The importance of the geometry is manifested by the failure of $L^p$ boundedness for every finite $p$ when the averaging surface is the boundary of a cube. Testing $\fM$ on $f_\eps(y) = \chi_{[-\eps,\eps]\times[-1/2,1/2]^{d+m-1}}(y)$ shows that $\fM$ is unbounded on $L^p(G)$ unless $p=\infty$.

(e) Theorem \ref{thm:main-hne} asserts the existence of $p_\circ$. A closer examination of the proof may lead to the optimal value of $p_\circ$, but we do not pursue it here. As the proof in Section \ref{sec:necessary} indicates, the value of $p_\circ$ depends not only on the {\L}ojasiewicz exponent but also on the behavior of the singular values $\{s_k(\theta)\}$ in a subtle way. It would be interesting to obtain an explicit characterization of the optimal value of $p_\circ$.

\subsection*{Outline of the paper}

In Sections \ref{sec:reductions} and \ref{sec:l2estspf}, we establish Theorems \ref{thm:main-m=1} and \ref{thm:main-hsu}. We work in the setting of general $m$ to prove both theorems simultaneously. Section \ref{sec:reductions} reduces the proofs to estimates for an oscillatory integral operator, following the argument in \cite{RyuSeeger}. We first show that we may assume $S_1,\dots,S_m$ are linearly independent (see Theorem \ref{thm:main2}), and then perform a dyadic frequency decomposition of the surface measure $d\mu$. An important difference from \cite{RyuSeeger} is an additional decomposition in the variable that determines the size of $\La^\intercal\theta$. This yields an extra gain in the $L^1$ estimate for each piece. We then use standard arguments involving Sobolev embedding and a shear transformation to reduce matters to estimates for an oscillatory integral operator.

Section \ref{sec:l2estspf} is devoted to the proof of the key estimate, Proposition \ref{prop:locall2-2}. As in \cite{RyuSeeger}, the basic strategy is to decompose the operator into building blocks and combine the estimates for these blocks using almost orthogonality. However, to obtain the desired bounds, we must keep careful track of the sizes of the singular values of $S^\theta$. This is an additional consideration that does not arise in \cite{RyuSeeger}.

Section \ref{sec:necessary} concerns necessary conditions. We first prove Theorem \ref{thm:main-hne}. Our construction of the example builds on \cite{RS26}. However, we also need to account for the scale of the singular value $s_r(\theta)$, which makes the argument somewhat more involved. In Proposition \ref{prop:matrixconst}, we construct examples for which $\cL(\cC_r,d_\perp\,|\,\cB)>2r$. To obtain a large {\L}ojasiewicz exponent, we begin with a polynomial that is a variant of an example of Masser and Philippon (see \cite[Example 15]{JKS92}). We then construct the matrices $\{S_k\}_{k=1}^m$ so that $\cC_r(1,\theta')$ has this polynomial as a factor. For this construction, we invoke a result from \cite{BDDHP17} on determinantal representations by affine linear matrix pencils, stated below as Theorem \ref{thm:determinantal}.

The details of two preliminary reductions are deferred to the appendices. Appendix \ref{app:basic-reductions} derives the sufficiency parts of Theorems \ref{thm:main-m=1} and \ref{thm:main-hsu} from Theorem \ref{thm:main2} below. Appendix \ref{app:osc-reduction} elaborates on the arguments that reduce Proposition \ref{prop:locall2} to the oscillatory integral estimates stated in the main text.

\subsection*{Notation}
We often use subscripts to denote partial derivatives. For nonnegative quantities $a$ and $b$, the notation $a\lesssim b$ means that $a\le Cb$ for some constant $C>0$. We write $a\sim b$ when both $a\lesssim b$ and $b\lesssim a$ hold. For a real quantity $c$ and a nonnegative quantity $a$, the notation $c=O(a)$ means that $|c|\le Ca$ for some constant $C>0$.

\section{Dyadic decomposition and reduction to $L^2$ estimates}\label{sec:reductions}
In this section, we make several reductions of the sufficiency statements in Theorems \ref{thm:main-m=1} and \ref{thm:main-hsu} to Proposition \ref{prop:locall2-2}, the key estimate of the proof. The necessity of the condition on $p$ in Theorems \ref{thm:main-m=1} and \ref{thm:main-hsu} follows from \cite{RS26}. Henceforth, when we refer to either theorem, we mean its sufficiency part.

\subsection{A graph formulation of the maximal estimate} In this subsection, we state the maximal estimate in a setting where the $(d+1)\times d$ matrices $S_1,\dots,S_m$ are linearly independent and $\Sigma$ is parametrized as a graph. Appendix \ref{app:basic-reductions} shows that this formulation implies Theorems \ref{thm:main-m=1} and \ref{thm:main-hsu}.

\begin{thm}\label{thm:main2}
    Let $\{S_i\}_{i=1}^m$ be $(d+1)\times d$ matrices as in \eqref{def:Si} and suppose that $S_1,\dots ,S_m$ are linearly independent.
    Let $r$ be as in \eqref{i:defr} and assume that $2\le r\le d-2$. If $m\ge 2$, we further assume that $\cL(\cC_r, d_\perp\,|\,\cB) = 2r$.    
    Let $g$ be a $C^\infty$ function  defined in an open domain near $0\in \bbR^{d-1}$ such that
    \begin{align}\label{r:propertyofg}
        g(0) = 1,\quad g'(0) = 0,\quad \text{$g''(0)$ is negative-definite}.
    \end{align}
    Define $\Gamma_g(\om') = g(\om')e_1 + \sum_{i=2}^d \om_i e_i$ and $\mathcal A_t$ by
    \begin{align*}
        \mathcal A_t f(x) = \int f(\ux - t\Gamma_g(\om'), \ox - t\sum_{i=1}^m (\ux, t)^\intercal S_i \Gamma_g(\om') e_{d+i}) \varphi_0(\om') d\om',
    \end{align*}
    where $\varphi_0$ is a $C^\infty_c$ function supported in the domain of $g$. Then there is a neighborhood $\cU$ of $0$ such that for $\varphi_0$ compactly supported in $\cU$ and $p> \frac{r+1}{r} $,
    \begin{align}\label{e:ineqforAt}
        \big\|\sup_{t>0}\big|\cA_t f\big|\big\|_{L^p(\R^{d+m})}\le C\|f\|_{L^p(\R^{d+m})}
    \end{align}
    holds for all $f\in C_c^\infty(\R^{d+m})$ with a constant $C$ independent of $f$. 
    \end{thm} 
    
    \begin{rem}\label{rem:endpoint}
    If  $I\subset (0,\infty)$ is a compact interval and if $m=1$ and $r\le d-3$ one can also show  for the endpoint $p=\frac{r+1}{r}$ the restricted weak type estimate
    \begin{align}\label{e:ineqlocal}
        \big\|\sup_{t\in I}\big|\cA_t f\big|\big\|_{L^{\frac{r+1}{r},\infty}(\R^{d+1})}\le C\|f\|_{L^{\frac{r+1}{r},1}(\R^{d+1})}.
    \end{align}
\end{rem}

The reduction of Theorems \ref{thm:main-m=1} and \ref{thm:main-hsu} to the theorem above closely follows the argument in \cite{RyuSeeger} that uses singular value decomposition, localization, and a shear transformation. In these reductions, both the rank $r$ and the {\L}ojasiewicz exponent $\mathcal L(C_r,d_\perp\mid\mathcal B)$ stay unchanged. In view of the similarity of the arguments, we present the details of this reduction in Appendix \ref{app:basic-reductions} for readability.

%and then discuss how this result can be applied to deduce Theorem \ref{thm:main-m=1} and \ref{thm:main-hsu}. The proof of Theorem \ref{thm:main2} will start in  \S\ref{sec:dyadicdecomp}.

%Next, let $g$ be a $C^\infty$ function on $\R^{d-1}$ such that
%\begin{align}\label{i:relsg}
%    g(0) = 1,\quad g'(0) = 0,\quad \text{$g''(0)$ is positive-definite}. 
%\end{align}
%Let $P$ be the $(d-1)\times d$ matrix given by $P = (0\ \ I_{d-1})$.
%Setting $\Gamma_g(\om') = g(\om')e_1 + P^\intercal \om'$, we define
%\begin{align}\label{def:cA}
%    \cA f(x,t) = \int f(\ux - t\,\Gamma_g(\om'), \ox - t(\ux, t)^\intercal \vec S \Gamma_g(\om')) \chi_\circ(\om') d\om',
%\end{align}
%where $x = (\ux,\ox)\in \R^d\times \R^n$ and $\chi_\circ\in C_c^\infty(\R^{d-1})$ such that $\chi_\circ$ is supported in a $\delta$-neighborhood of the origin. By Proposition \ref{prop:red1}, we may assume that $S_1,\dots, S_n$ are linearly independent throughout the rest of the article. Let $\cM f$ be the associated maximal function defined by
%\begin{align}\label{def:cM}
%    \cM f(x) = \sup_{t>0} |\cA f(x,t)|.
%\end{align}
%We discuss how Theorem \ref{thm:main} (with $S_1,\dots,S_n$ being linearly independent) reduces to estimating the $L^p(\R^{d+n})$ norm of $\cM f$. In the case $\La = 0$, the corresponding reduction step is thoroughly explained in \cite{RyuSeeger}, where the authors used a change of variables involving a rotation and a shear transformation. Our argument is quite similar to that of \cite{RyuSeeger}.

\subsection{Dyadic decomposition and $L^2$ estimates}\label{sec:dyadicdecomp}
For a given collection of $d\times d$ skew-symmetric matrices $J_1,\dots,J_m$ and vectors $\{\la_i\}_{i=1}^m\subset \R^d$, assume that $S_i$, $1\le i\le m$, are linearly independent. We define the quantity $r^\sharp$ as
\[
r^\sharp = \min\{r(\theta):\theta\in \bbS^{m-1}\}.
\]
Note that
\begin{align}\label{cond:rsharp}
    2\le r^\sharp \le r.
\end{align}
In fact, if $\La^\intercal \theta\neq 0$, then $r(\theta)\ge 2$ by the assumption $r\ge 2$. Otherwise, linear independence implies $J^\theta\neq 0$, and a nonzero skew-symmetric matrix has rank at least two.
We also note that $r^\sharp = r$ when $m=1$.
To prove Theorem \ref{thm:main2}, we follow the approach taken in \cite{RyuSeeger} (see also \cite{MuellerSeeger2004}). Let $\nu$ be the distribution on $\R^{d+m}$ defined by 
\[
\inn{\nu}{f} = \int f(g(\om'),\om',\La \Gamma_g(\om'))\varphi_0(\om')d\om',
\]
where $g, \varphi_0, \Gamma_g(\om')$ are defined as in Theorem \ref{thm:main2}. Observe that $\cA_t f = f\ast \nu_t$.

Note that $\nu$ is expressed as the pairing of the measure \[
d\nu(x) = \varphi_1(x_1,\ox) \varphi_0(x') \delta(x_1-g(x'), \ox - \La \ux),
\]
where $\delta$ is the Dirac measure in $\bbR^{m+1}$, $\varphi_1\in C_c^\infty(\R^{m+1})$ such that
\begin{align*}
    \varphi_1(g(\om'),\La\Gamma_g(\om')) = 1 \quad &\text{for}\ \om'\in \supp(\varphi_0) \\
    \varphi_1(x_1,\ox) = 0 \quad &\text{unless}\ |x_1 - 1| + |\ox - \La e_1|\le c
\end{align*}
with a small constant $c = c(g,\La)>0$. Using the Fourier inversion formula, we  write $\nu$ as an oscillatory integral distribution, 
\begin{align*}
    \nu(x) = \frac{\varphi_1(x_1,\ox) \varphi_0(x')}{(2\pi)^{m+1}} \int e^{i\sigma(x_1 - g(x')) + i\inn{\tau}{\overline x - \La \ux}} d\sigma d\tau,\quad \sigma\in \R,\ \tau\in \R^m.
\end{align*}
We dyadically decompose the integral in the frequency variables $\sigma, \tau$. Let $\zeta_0\in C_c^\infty(\R)$ such that $\supp(\zeta_0)\subset (-1,1)$, $\zeta_0 \equiv 1$ on $(-1/2,1/2)$. For $k\ge 1$, we define $\zeta_k(s) = \zeta_0(2^{-k}s) - \zeta_0(2^{1-k}s)$, so that $\supp(\zeta_k)\subset (-2^k, -2^{k-2})\cup(2^{k-2},2^k)$ for $k\ge 1$ and $\sum_{k\ge 0} \zeta_k \equiv 1$. Let
\[
\nu^k(x) = \frac{\varphi_1(x_1,\ox) \varphi_0(x')}{(2\pi)^{m+1}} \int \zeta_k(\sqrt{\sigma^2 + |\tau|^2}) e^{i\sigma(x_1 - g(x')) + i\inn{\tau}{\ox - \La \underline{x}}} d\sigma d\tau
\]
for $k\ge 0$. We make an extra dyadic decomposition in $\sigma$. For $k\ge 1$, $l\le k$, define $\zeta_{k,l}(\sigma,\tau) = \zeta_k((\sigma^2 + |\tau|^2)^{1/2})\zeta_{k-l}(\sigma)$ if $l<k$ and $\zeta_{k,l}(\sigma,\tau) = \zeta_k((\sigma^2 + |\tau|^2)^{1/2})\zeta_0(\sigma)$ if $l=k$. Note that $\zeta_{k,l}$ identically vanishes if $l\le -3$. Then, for $k\ge 1,\ l\le k$, we set
\[
\nu^{k,l}(x) = \frac{\varphi_1(x_1,\ox) \varphi_0(x')}{(2\pi)^{m+1}} \int \zeta_{k,l}(\sigma,\tau) e^{i\sigma(x_1 - g(x')) + i\inn{\tau}{\overline  x - \La \underline x}} d\sigma d\tau.
\]
Then, $\mathcal A_t f$ is decomposed as $\mathcal A_t f = f\ast \nu^0_t + \sum_{k\ge 1}\sum_{l\le k}f\ast\nu^{k,l}_t$, with convergence in the sense of distributions.

We distinguish between the two cases $i)$ $r^\sharp = r$, $ii)$ $r^\sharp < r$. If $r^\sharp = r$, we estimate $\nu^0$ and $\nu^{k,l}$ directly. If $r^\sharp < r$, we introduce an additional dyadic decomposition adapted to the size of $\La^\intercal \tau$.

To do so, we denote the rank of $\La$ by $m^\sharp$. Note that $m^\sharp = m = 1$ when $m=1$. If $m^\sharp=m$, then $\ker(\Lambda^\intercal)=\{0\}$, and hence $r^\sharp = r$. Thus, $r^\sharp < r$ implies $m^\sharp < m$. After linear changes of variables in $\ox$ and in $\oy$ we may assume that
\Be\label{eq:la-split} 
\begin{aligned} &\la_1,\dots, \la_{m-m^\sharp} = 0,\\  &\lambda_{m-m^\sharp+1},\dots, \la_m {\text{ are linearly independent.}}
\end{aligned}
\Ee
Now we write $\tau = (\tau^{(1)}, \tau^{(2)})\in \R^{m-m^\sharp}\times \R^{m^\sharp}$ and denote by $P^{(1)}$, $P^{(2)}$ the corresponding projections $P^{(1)}: \bbR^m\to \bbR^{m-m^\sharp}$, $P^{(2)}: \bbR^m\to \bbR^{m^\sharp}$, respectively. We perform an additional dyadic decomposition in the $\tau^{(2)}$-variables as follows. For $-2\le j< l$, let
\[
\nu^{k,l,j}(x) = \frac{\varphi_1(x_1,\ox) \varphi_0(x')}{(2\pi)^{m+1}} \int \zeta_{k,l}(\sigma,\tau)\zeta_{k-j}(|\tau^{(2)}|) e^{i\sigma(x_1 - g(x')) + i\inn{\tau}{\ox - \La \ux}} d\sigma d\tau.
\]
If $j = l$, we define $\nu^{k,l,l}(x)=$ 
\[ 
\frac{\varphi_1(x_1,\ox) \varphi_0(x')}{(2\pi)^{m+1}} \int \zeta_{k,l}(\sigma,\tau)\sum_{j=l}^k\zeta_{k-j}(|\tau^{(2)}|) e^{i\sigma(x_1 - g(x')) + i\inn{\tau}{\ox - \La \ux}} d\sigma d\tau.
\]
%\end{multline*} 
Then, $\nu^{k,l} = \sum_{-2\le j\le l}\nu^{k,l,j}$ for all $k\ge 1$ and $-2\le l\le k$.
%This ends the dyadic decomposition process for $\nu$, in case \textit{ii)}.

We now state our main local estimates for the decomposed pieces. For $-2\le j\le l$, let
\begin{align*}
    \mathbf B_{j,l}(r^\sharp,r) = \min(2^{jr}, 2^{l(r-r^\sharp)}).
\end{align*}

\begin{prop}\label{prop:locall2}
Let $k,l,j\in \Z$ such that $k\ge 1$, $-2\le l\le k$, and $-2\le j\le l$. Then there exists a constant $C>0$ independent of $k,l,j$ such that the following a priori estimates hold for all $f\in C_c^\infty(\R^{d+m})$. 

(i) If $r^\sharp = r$ then  
    \begin{multline*}
        \Big(\int_1^2 \|f\ast \nu_t^{k,l}\|_2^2 \, dt\Big)^{1/2}+ 2^{-k} \Big(\int_1^2\|\partial_t(f\ast \nu_t^{k,l})\|_{2}^2 \,dt\Big)^{1/2} 
        \le C 2^{-\frac{d-2}{2}k}2^{\frac{d-r-2}{2}l}\|f\|_{2}.
\end{multline*}

(ii) If $r^\sharp < r$ then 
\begin{multline*}   \Big(\int_1^2\|f\ast \nu_t^{k,l,j}\|_2^2\, dt\Big)^{1/2} + 
2^{j-k}\Big(\int_1^2 \|\partial_t(f\ast \nu_t^{k,l,j})\|_2^2\,  dt\Big)^{\frac 12} \\
        \le C 2^{-k\frac{d-2}{2}} 2^{\frac 12 l(d-r-2)} (\mathbf B_{j,l}(r^\sharp,r))^{\frac12}\|f\|_{2}.
\end{multline*} 
% \begin{align*}
%         \|f\ast \nu_t^{k,l,j}\|_{L^2(\R^{d+m}\times [1,2])} + &2^{j-k}\|\partial_t(f\ast \nu_t^{k,l,j})\|_{L^2(\R^{d+m}\times [1,2])}\\
%         &\le C 2^{-\frac{d-2}{2}k} 2^{\frac12(d-r-2)l} 2^{\frac12(r-r^\sharp-m^\sharp)j}\|f\|_{L^2(\R^{d+m})}.
% \end{align*}
%for all $f\in C_c^\infty(\R^{d+m})$.
%\end{itemize}    
\end{prop}

We establish these local $L^2$ estimates in Sections \ref{sec:red} and \ref{sec:l2estspf}. For the remainder of this section, we explain how the  proposition implies Theorem \ref{thm:main2}.

\subsection{From dyadic estimates to maximal bounds} We begin by proving \eqref{e:ineqforAt} with $\sup_{t>0}|\cA_t f|$ replaced by a local maximal operator. This result will follow from $L^p$ estimates for $f\ast \nu_t^{k,l,j}$ ($f\ast \nu_t^{k,l}$ if $r = r^\sharp$) which we now deduce from Proposition \ref{prop:locall2}.

We first analyze the case $r^\sharp <r$. Note that
\[
\nu_t^{k,l,j}(x) = \int a^\nu_{k,l,j}(x,t,\sigma,\tau) e^{i\phi(x,t,\sigma,\tau)}d\sigma d\tau,
\]
where $a^\nu_{k,l,j}$ and $\phi$ are defined by
\begin{align*}
    a^\nu_{k,l,j}(x,t,\sigma,\tau) &=(2\pi)^{-m-1}t^{-d-2m} \varphi_1(t^{-1}x_1, t^{-2}\ox)\varphi_0(t^{-1}x')\zeta_{k,l}(\sigma,\tau)\zeta_{k-j}(|\tau^{(2)}|), \\
\phi(x,t,\sigma,\tau) &= t^{-1}\sigma(x_1 - tg(t^{-1}x')) + t^{-2}\inn{\tau}{\ox-t\La\ux}
\end{align*}
for $j<l$. $a^\nu_{k,l,l}$ is defined similarly but $\zeta_{k-j}$ is replaced by $\sum_{j=l}^k\zeta_{k-j}$.
The $t$-derivative of the phase function is 
%$\partial_t\phi$  is 
%expressed as
$\partial_t\phi= -u^{[1]} - u^{[2]}- u^{[3]}$ where 
% \begin{align*}
%     \partial_t\phi(x,t,\sigma,\tau) &= -\big(t^{-2} \sigma(x_1 - tg(t^{-1}x')) - 2t^{-3}\inn{\tau}{\ox-t\La\ux}\big) \\
%     &\hspace{30mm} - t^{-1}\sigma \partial_t(t g(t^{-1}x')) - t^{-2}\inn{\tau}{\La \ux} \\
%     & =: \phi^{[1]}(x,t,\sigma,\tau) + \phi^{[2]}(x,t,\sigma,\tau) + \phi^{[3]}(x,t,\sigma,\tau).
% \end{align*}
\begin{align*}
    u^{[1]} (x,t,\sigma,\tau) &= t^{-2} \sigma(x_1 - tg(t^{-1}x')) + 2t^{-3}\inn{\tau}{\ox-t\La\ux} 
    \\
        u^{[2]} (x,t,\sigma,\tau) &= t^{-1}\sigma \partial_t\big(t g(t^{-1}x')\big) 
        \\
            u^{[3]} (x,t,\sigma,\tau) &= t^{-2}\inn{\tau}{\La\ux}\,. 
\end{align*}
Since $|\tau^{(2)}|\lc  2^{k-j+2}$, and  $|\sigma|\lc 2^{k-l}$ we get $u^{[3]}=O(2^{k-j})$ and $u^{[2]}=O( 2^{k-l})$ on $\supp(a^\nu_{k,l,j})$.
In view of the derivative  bounds \[\big |\partial_\sigma^n \partial_{\tau^{(1)}} ^\alpha 
\partial_{\tau^{(2)}} ^\beta 
\big (\zeta_{k,l}(\sigma,\tau^{(1)}, \tau^{(2)})\zeta_{k-j}(|\tau^{(2)}|) \big) \big|\lesssim_{n,\alpha,\beta} 2^{-n(k-l)}2^{-|\alpha|k} 2^{-|\beta| (k-j)}\]
we may integrate by parts and get the kernel bound 
\begin{multline}
%\begin{aligned}
\label{e:bdsnuklj}
    |\nu^{k,l,j}_t(x)|, \,\, \,  |2^{j-k}\partial_t\nu^{k,l,j}_t(x)|, \,\,\, 
    |2^{-k}\nabla_x\nu^{k,l,j}_t(x)|,\,\, \, |2^{j-2k}\nabla_x\partial_t\nu^{k,l,j}_t(x)| 
    \\
    % &{\color{green} \quad\lesssim_N \frac{2^{k-l}}{(1 + 2^{k-l}|x_1 - g(x')|)^N}\frac{2^{km}}{(1 + 2^k|\ox - \La\ux|)^N},}
   % \\
        %\quad\lesssim_N \frac{2^{k-l}}{(1 + 2^{k-l}|x_1 - tg(t^{-1}x')|)^N}\frac{2^{km}}%{(1 + 2^k|\ox - t\La\ux|)^N},
        \quad\lesssim_N \frac{2^{k-l}}{(1 + 2^{k-l}|x_1 - tg(t^{-1}x')|)^N}
        \frac{2^{k(m-m^\sharp)}}{(1 + 2^k|P^{(1)}\ox |)^N}
\frac{2^{(k-j)m^\sharp}}{(1 + 2^{k-j}|P^{(2)}(\ox - t\La\ux)|)^N}  
%\end{aligned}
\end{multline}
provided that $|t^{-1}x_1-1| + |t^{-1}x'| + |t^{-2}\ux-\La e_1|\le c$ with small $c>0$ and $t\in [1,2]$. 
% To deduce it, we also use the bounds $|\partial_\sigma^n \partial_\tau^\alpha (\zeta_{k,l}(\sigma,\tau)\zeta_{k-j}(|\tau^2|))|\lesssim_{n,\alpha} 2^{n(l-k)}2^{|\alpha|k}$ where $(n,\alpha)\in \N\times \N^{m}$.
% {\color{green}For later use, we note another pointwise bounds
% \begin{align*}
%     |2^{-k}\nabla_x\nu^{k,l,j}_t(x)|,\ &|2^{j-2k}\nabla_x\partial_t\nu^{k,l,j}_t(x)| \\
%     &\quad\lesssim_N \frac{2^{k-l}}{(1 + 2^{k-l}|x_1 - g(x')|)^N}\frac{2^{km}}{(1 + 2^k|\ox - \La\ux|)^N},
% \end{align*}
% }
These pointwise estimates yield  in turn %\as{Here we do not need the $t$-integral}
\begin{subequations} 
\begin{align} 
\label{e:bdsnuklj-2a}
    \|\nu_t^{k,l,j}\|_{L^1(\R^{d+m})} + 2^{j-k}\|\partial_t \nu_t^{k,l,j}\|_{L^1(\R^{d+m})}\lesssim 1.
    \\ 
\label{e:bdsnuklj-2}
    \|\nabla_x\nu_t^{k,l,j}\|_{L^1(\R^{d+m})} + 2^{j-k}\|\nabla_x\partial_t \nu_t^{k,l,j}\|_{L^1(\R^{d+m})}\lesssim 2^k.
\end{align}
\end{subequations}

We thus get  $L^1$ estimates for the convolution of $f$ with 
$\nu_t^{k,l,j}$  and $\partial_t \nu_t^{k,l,j}$. An interpolation  argument invoking the  $L^2$ estimate from Proposition \ref{prop:locall2} $(ii)$ yields,  for $1\le p\le 2$,
\begin{align}
    \begin{aligned}\label{e:locallpklj}
        \Big(\int_1^2 \|f\ast \nu_t^{k,l,j}\|_p^p\,dt\Big)^{1/p} &%{p}
        + 2^{j-k}\Big(\int_1^2\|\partial_t(f\ast \nu_t^{k,l,j})\|_{p}^p dt\Big)^{1/p} 
        %\times [1,2])}\
        \\ &\le C 2^{-k\frac{d-2}{p'}}2^{l\frac{d-r-2}{p'}}(\mathbf B_{j,l}(r^\sharp,r))^{\frac{1}{p'}}\|f\|_{p}
    \end{aligned}
\end{align}

We turn to the other case $r^\sharp = r$, in which we can obtain the $L^p$ estimates
\begin{align}
    \begin{aligned}\label{e:locallpkl}
        \Big(\int_1^2 \|f\ast \nu_t^{k,l}\|_{p}^p\,dt\Big)^{\frac 1p}+ 2^{-k}\Big(\int_1^2 \|\partial_t(f\ast &\nu_t^{k,l})\|_{p}^p\, dt\Big)^{\frac 1p} \\
        &\le C 2^{-k\frac{d-2}{p'}}2^{ l\frac{d-r-2}{p'}}\|f\|_{p}
    \end{aligned}
\end{align}
after using $L^2$ bounds stated in Proposition \ref{prop:locall2} $(i)$ and the $L^1$ bounds
\begin{align*}%\label{e:fnul1}
    \|f\ast \nu_t^{k,l}\|_{1} + 2^{-k}\|\partial_t(f\ast \nu_t^{k,l})\|_{1}  \lesssim \|f\|_{1}.
\end{align*}
The $L^1$ bounds follow as above, with $u^{[3]}=O(2^k)$ and without the $j$-decomposition. In this case $r=r^\sharp$ we get the pointwise  bounds 
%\begin{align}\begin{aligned}
\begin{multline}\label{e:bdsnukl}
  \  |\nu^{k,l}_t(x)|,\ |2^{-k}\nabla\nu^{k,l}_t(x)|,\ |2^{-k}\partial_t\nu^{k,l}_t(x)|,\ |2^{-2k}\partial_t\nabla\nu^{k,l}_t(x)| \\
    \quad\lesssim_N \frac{2^{k-l}}{(1 + 2^{k-l}|x_1 - tg(t^{-1}x')|)^N}\frac{2^{km}}{(1 + 2^k|\ox - t\La\ux|)^N}.
%\end{aligned}
\end{multline}
%\end{align}

Next, we will deduce maximal inequalities in Theorem \ref{thm:main2} using $L^p$ estimates \eqref{e:locallpklj} and \eqref{e:locallpkl}. 

\begin{proof}[$L^p$ estimates for the local maximal operator]
%The restricted weak type inequality \eqref{e:ineqlocal}] \as{revise later}
Since $\nu^0$ is compactly supported and bounded, the local maximal function $f\mapsto \sup_{1\le t\le 2} |f\ast \nu_t^0|$ maps $L^p\to L^p$ for $1\le p\le \infty$. 
% This implies
% \begin{align}\label{e:lpest-nu0}
%     \big\|\sup_{t>0}\big|f\ast \nu_t^0\big|\big\|_{L^p(\R^{d+m})}\le C\|f\|_{L^p(\R^{d+m})}
% \end{align}
% for all $p\in (1,\infty]$.
%Hence the maximum of $f\ast \nu_t^0$ satisfies a much better boundedness property than we need.

We now estimate the maximal function $\sup_{t\in [1,2]} |f\ast  (\nu_t - \nu_t^0)|$ and distinguish the  cases $r^\sharp < r$ and $r^\sharp = r$. By the Sobolev inequality $\|F\|_\infty \lesssim \|F\|_p + \|F\|_p^{1/p'}\|F'\|_p^{1/p}$, the $L^p$ estimates \eqref{e:locallpklj} imply
\begin{align}\label{e:sup-nuklj}
        \|\sup_{t\in [1,2]}|f\ast \nu_t^{k,l,j}|\|_{p} 
        \lc 2^{k(\frac1p-\frac{d-2}{p'})}2^{l\frac{d-r-2}{p'}}(\mathbf B_{j,l}(r^\sharp,r))^{\frac{1}{p'}}2^{-\frac jp}\|f\|_{p}
\end{align}
where $1\le p\le 2$. We note that
\[
\mathbf B_{j,l}(r^\sharp,r) \le \begin{cases}
    2^{jr} & -2\le j\le j_\circ, \\
    2^{l(r-r^\sharp)} & j_\circ < j\le l,
\end{cases}\quad j_\circ = \frac{l(r-r^\sharp)}{r}.
\]
Substituting this into \eqref{e:sup-nuklj} and summing over $j$ yields
\begin{align}\label{e:sup-nukl}
        \|\sup_{t\in [1,2]}|f\ast \nu_t^{k,l}|\|_{p} 
        \lc 2^{k(\frac1p-\frac{d-2}{p'})}2^{l\alpha(d,p,r,r^\sharp)}\|f\|_{p}
\end{align}
for $\frac{r+1}{r} < p \le 2$, where 
\[\alpha(d,p,r,r^\sharp) = \frac{d-2}{p'}-\frac1p-\frac{r^\sharp}{r}(\frac{r}{p'} - \frac1p).
\]
Since $r\le d-2$, $r^\sharp < r$, and $\frac{r+1}{r} < p \le 2$, we have $\alpha(d,p,r,r^\sharp)> 0$. Thus, summing \eqref{e:sup-nukl} over $l:-2\le l\le k$ gives 
\[
\|\sup_{t\in [1,2]}|f\ast \nu_t^{k}|\|_p 
        \lc  2^{-k\frac{r^\sharp}{r}(\frac{r}{p'} - \frac1p)}\|f\|_{p}
\] 
for every $\frac{r+1}{r} < p \le 2$. Note that $\frac{r}{p'} - \frac1p > 0$ if $\frac{r+1}{r} < p \le 2$. Also, by \eqref{cond:rsharp}, $r^\sharp/r>0$. Hence, taking the last summation over $k$, we obtain the desired estimate for the local maximal operator.

We now address the case $r^\sharp = r$, where the proof is simpler. Similarly, using the Sobolev inequality and \eqref{e:locallpkl} gives
\begin{align*}
    \|\sup_{t\in [1,2]}|f\ast \nu_t^{k,l}|\|_{p} \lc 2^{k(\frac1p-\frac{d-2}{p'})}2^{l\frac{d-r-2}{p'}}\|f\|_{p}
\end{align*}
which implies
\begin{align*}
    \|\sup_{t\in [1,2]}|f\ast \nu_t^{k}|\|_p \lc  \begin{cases}
        2^{k(\frac1p-\frac{r}{p'})}\|f\|_{p} & \text{if}\ r<d-2, \\
        (k+3) 2^{k(\frac1p-\frac{r}{p'})}\|f\|_{p} & \text{if}\ r=d-2.
    \end{cases}
\end{align*}
Summing the estimate over $k$, we get the local version of \eqref{e:ineqforAt} for $\frac{r+1}{r}<p\le 2$, as desired. Interpolation with the trivial $L^\infty$ bound gives the remaining range $2<p\le\infty$. We also note that if $r\le d-3$, Bourgain's interpolation trick can be applied to get the restricted weak type bound \eqref{e:ineqlocal} for $p=\frac{r+1}{r}$, for the local maximal operator.
\end{proof}

\begin{rem}\label{rem:lpest-hrank}
    From the proof above, we can also establish $L^p$ estimates for the local maximal operator when $d-1\le r\le d$. We only consider the case $r^\sharp<r$, as the other case is simpler. Let $d-1\le r\le d$. Then, we note that $\alpha(d,\frac{d-1}{d-2},r,r^\sharp) < 0$. Choose $p^\sharp>\frac{d-1}{d-2}$ sufficiently close to $\frac{d-1}{d-2}$ so that $\alpha(d,p,r,r^\sharp) < 0$ for $\frac{d-1}{d-2}<p<p^\sharp$. Summing \eqref{e:sup-nukl} first over $l$ and then over $k$, we conclude that $\fM_I f$ is bounded on $L^p(G)$ for $\frac{d-1}{d-2}<p<p^\sharp$ and a compact interval $I\subset (0,\infty)$. Interpolating with the $L^\infty$ bound gives the $L^p$ boundedness for the local maximal operator in the range $\frac{d-1}{d-2}<p\le \infty$.
\end{rem}

Now we prove \eqref{e:ineqforAt}. The proof essentially follows the same argument as in \cite[Section 4]{RyuSeeger} (or \cite{MuellerSeeger2004}), thus we shall briefly outline it. Since the case $r = r^\sharp$ can be treated similarly to the case $r^\sharp < r$, we include the proof only for $r^\sharp < r$. It suffices to prove the desired maximal estimate for $\frac{r+1}{r}<p\le2$, since the range $2<p\le\infty$ then follows by interpolation with the trivial $L^\infty$ bound.

We have $\sup_{t>0} |f\ast \nu_t^0|\lc CM_{HL}f$ with a constant $C>0$, where $M_{HL}$ denotes the Hardy-Littlewood maximal operator associated with the Carnot balls. Hence, this maximal function is bounded on $L^p$ for $p>1$.

We now consider the global maximal function associated with $\nu^{k,l,j}$  for $-2\le j\le l\le k$. Integration by parts shows that $|\int \nu^{k,l,j} dx|\lc 2^{-kN}$. 
Thus, following \cite{MuellerSeeger2004}, we can write $\nu^{k,l,j}(x) = \fK^{k,l,j}(x) + \rho^{k,l,j}(x)$ where $\fK^{k,l,j}, \rho^{k,l,j}\in C_c^\infty(\R^{d+m})$ such that $\fK^{k,l,j}$ has the cancellation property $\int \fK^{k,l,j}(x) dx = 0$ and the $L^1$ norms of $\rho^{k,l,j}$ and its  derivatives up to order $2(d+m)$ are $O(2^{-kN})$. 
Define $\fK^{k,l,j}_t(\ux,\ox) = t^{-d-2m}\fK^{k,l,j}(t^{-1}\ux, t^{-2}\ox)$ and $\rho^{k,l,j}_t(\ux,\ox) = t^{-d-2m}\rho^{k,l,j}(t^{-1}\ux, t^{-2}\ox)$.
By the above we dominate $\sup_{t>0} |f*\rho^{k,l,j}_t(x)|\lc_N 2^{-kN}  M_{HL} f(x)$. 
Hence the task is to obtain an inequality for  $\sup_{t>0}|f\ast \fK_t^{k,l,j}|$ that is summable over $-2\le j\le l \le  k$. For that, it suffices to prove
\begin{multline}\label{e:supl2Lp}
    \Big\|\Big(\sum_{n\in\Z}\big(\sup_{s\in [1,2]}|f\ast \fK_{2^ns}^{k,l,j}|\big)^2\Big)^{1/2}\Big\|_{L^p(\R^{d+m})} \\
    \le C (1+k)^{\frac 1p} 2^{k(\frac1p-\frac{d-2}{p'})}2^{l\frac{d-r-2}{p'}}(\mathbf B_{j,l}(r^\sharp,r))^{\frac{1}{p'}}2^{-\frac jp}\|f\|_{L^p(\R^{d+m})}
\end{multline}
for $1< p\le 2$. 
By a standard Sobolev-type argument involving the fundamental theorem of calculus, the inequality  \eqref{e:supl2Lp} is a consequence of 
\begin{align}
   \notag \Big(\int_1^2 \Big\|\Big(\sum_{n\in\Z}|f\ast \fK_{2^ns}^{k,l,j}|^2\Big)^{1/2}\Big\|_{p}^p&\,ds\Big)^{\frac 1p}  +\Big(\int_1^2
    \Big\|\Big(\sum_{n\in\Z} \Big|f\ast \frac{\frac{d}{ds}\fK_{2^ns}^{k,l,j}}{2^{k-j}}\Big |^2\Big)^{1/2}\Big\|_{p} ^p\,ds\Big)^{\frac 1p} 
    \\ \label{e:l2Lp}
    &\le C (1+k)^{\frac 1p} 2^{-k\frac{d-2}{p'}}2^{l\frac{d-r-2}{p'}}(\mathbf B_{j,l}(r^\sharp,r))^{\frac{1}{p'}}\|f\|_{p}.
\end{align}
We give the estimate for the first term on the left hand side of 
\eqref{e:l2Lp} as the argument for the second term is similar. We first prove the $L^2$ estimate and then interpolate it with the following weak type $(1,1)$ bound
\begin{align}
\begin{aligned}\label{e:weakL1}
    \Big|\Big\{(x,s):\Big(\sum_{n\in\Z}|f\ast \fK_{2^ns}^{k,l,j}(x)|^2\Big)^{1/2} > \alpha\Big\}\Big|\le C  (1+k)\frac{\|f\|_{1}}\alpha.
\end{aligned}
\end{align}
We first consider the $L^2$ estimate. 
Via a  Cotlar-Stein type lemma as in 
\cite[Lemma 3.2]{MuellerSeeger2004}, the $L^2$ bound (for the first term)  follows once we prove
\begin{align*}
   \Big(\int_1^2  \big\|f\ast \fK_{2^ns}^{k,l,j}\big\|_{2}^2 ds\Big)^{1/2} 
    \lc 2^{-k\frac{d-2}{2}}2^{l\frac{d-r-2}{2}}(\mathbf B_{j,l}(r^\sharp,r))^{\frac12}\|f\|_{2}
\end{align*}
for $n\in \Z$, and 
\begin{align*}
    \Big(\int_1^2 \big\|f\ast \fK_{2^ns}^{k,l,j} \ast \wt\fK_{2^{n'}s}^{k,l,j}\big\|_2^2 ds  \Big)^{1/2}
    \lc  2^{10mk}2^{-|n-n'|/10}\|f\|_{2}
\end{align*}
for $n,n'\in \Z$, where $\wt\fK_{2^{n'}s}^{k,l,j}(x) = \overline{\fK_{2^{n'}s}^{k,l,j}(-x)}$.
The first inequality is a consequence of Proposition \ref{prop:locall2} and scaling. To prove the second inequality, we use the argument in \cite{MuellerSeeger2004, RyuSeeger} which relies on the straightforward non-optimal bound 
%This  have this work, $\fK^{k,l,j}$ should satisfy
\begin{align}\label{e:goodboundcK}
    |\fK^{k,l,j}(x)| + 2^{-k}|\nabla\fK^{k,l,j}(x)|\lesssim 2^{k(m+1)}(1+|\ux|)^{-N}(1+|\ox|)^{-N},\quad N\in\N.
\end{align}
%This bound holds  by \eqref{e:bdsnukl} and the rapid decay $\rho^{k,l,j} = O(2^{-kN})$. As similar arguments apply to the second term in \eqref{e:l2Lp}  we get  inequality for $p=2$. 

To show the weak $L^1$ bound \eqref{e:weakL1}, we use the Calder\'on-Zygmund decomposition for the Carnot-Carath\'eodory balls on  a step two nilpotent group, as established in \cite[Ch.3A]{FollandStein1982}. Again, the arguments  in \cite{MuellerSeeger2004} or \cite[pp.16--20]{RyuSeeger} work with only notational modifications, in view of  \eqref{e:goodboundcK} 
(and corresponding estimates where $\fK^{k,l,j}_{2^n s} $ is replaced with  $2^{j-k} \frac{d}{ds} \fK_{2^n s}^{k,l,j}$).
This completes the sketch of the proof of Theorem \ref{thm:main2}. Our task is now to establish Proposition \ref{prop:locall2}.

\subsection{Oscillatory integral reduction and singular value bounds}\label{sec:red}
We reduce Proposition \ref{prop:locall2} to Proposition \ref{prop:locall2-2}, whose proof is given in Section \ref{sec:l2estspf}.

\begin{prop}\label{prop:locall2-2}
    Let $d\ge 2$ and let $k,l,j\in \Z$ be as in Proposition \ref{prop:locall2}. Let $\mathbf C\in [1,\infty)$. For a given $d\times d$ skew-symmetric matrix $J$ and $\lambda\in \R^d$, let $S = (J^\intercal\ \lambda)^\intercal$. Assume $\|S\|\le \mathbf C$. Let $s_1\ge s_2\ge \dots\ge s_d\ge 0$ be singular values of $S$. Let $P=\big(0\ I_{d-1}\big)$ be a $(d-1)\times d$ matrix, and let $g$ satisfy \eqref{r:propertyofg}. Suppose that there exist natural numbers $2\le r_1\le r_2\le d$ and a positive constant $C>0$ such that
    \begin{align}\label{b:lowerboundsi}
            s_i \ge C\quad \text{for}\ \ 1\le i\le r_1,\quad \Big|\prod_{i=1}^{r_2}s_i\Big|\ge C 2^{-jr_2}
    \end{align}
    hold. If $r_1 = r_2$ or $j=l$, only the first set of lower bounds, $s_i\ge C$ for $1\le i\le r_1$, is required.
    Define
    \begin{align*}
        T^{k,l} f(x,t) &= \int_{\R^d} e^{i2^k\Psi(x,t,y)} \chi(x,t,y) \zeta(2^lt \sigma(x',t,y_1)) f(y) dy,
    \end{align*}
    where
    \begin{align*}
        \Psi(x,t,y)&= y_1(x_1-tg(\tfrac{x'-y'}{t}) )+ (x,t)^\intercal S (P^\intercal y'-tg(\tfrac{x'-y'}{t})e_1), \\
    \sigma(x',t,y_1) & = y_1 + (x,t)^\intercal S e_1,
    \end{align*}
    and $\chi\in C_c^\infty(\R^{2d+1})$ such that
    \begin{align}\label{r:chisupp}
        \supp(\chi)\subset \{(x,t,y): t\sim 1,\ |x'|,\,|y'|\le \eps,\ |x_1|,\,|y_1|\lesssim 1\}
    \end{align}
    with a small number $\eps > 0$, and lastly, $\zeta$ stands for $s\mapsto \zeta_1(2s)$ if $l < k$ and $\zeta_0$ if $l = k$.
    Then we have
    \begin{align}\label{e:locall2-2}
        \|T^{k,l}\|_{L^2(\R^d)\to L^2(\R^d\times [1,2])}\lesssim
            2^{-\frac{d-1}{2}k}2^{\frac{d-r_2-2}{2}l} (\mathbf B_{j,l}(r_1,r_2))^{\frac12}
    \end{align}
    where the implicit constant depends on $d, \chi, \zeta, g,\mathbf C$, and $C$ in \eqref{b:lowerboundsi}, but is independent of $k,l,j$.
\end{prop}

Aside from establishing the singular value lower bounds \eqref{b:lowerboundsi}, the reduction of Proposition \ref{prop:locall2} to Proposition \ref{prop:locall2-2} closely follows the arguments in \cite{MuellerSeeger2004, RyuSeeger}, using a shear transformation and the Fourier transform. For readability, we present these standard reductions in Appendix \ref{app:osc-reduction} and establish \eqref{b:lowerboundsi} in this subsection.

For each $\overline u\in\R^m$, let $\Psi_{\overline u}$ and $\sigma_{\overline u}$ be defined as in Proposition \ref{prop:locall2-2}, with $S$ replaced by $S^{\overline u}$. We consider the oscillatory integral operator
\[
    \cG^{k,l}_{\overline u}f(x,t) = \int_{\R^d} e^{i2^k\Psi_{\overline u}(x,t,y)} \chi_{\overline u}(x,t,y) \zeta\big(2^l t \sigma_{\overline u}(x',t,y_1)\big) f(y)\,dy.
\]
Here, $\zeta$ is the cutoff function stated in Proposition \ref{prop:locall2-2}. The amplitudes $\chi_{\overline u}\in C_c^\infty(\R^{2d+1})$ arise from the reductions in Appendix \ref{app:osc-reduction}. They satisfy \eqref{r:chisupp}, with $\chi$ replaced by $\chi_{\overline u}$, and their derivatives are bounded uniformly in $\overline u$ for $|\overline u|\lesssim1$.

%We treat the case $2^{-l}\sim 1$ separately below. For sufficiently large $l$, the frequency cutoffs restrict the relevant parameters to $|\overline u|\sim1$. 
As shown in Appendix \ref{app:osc-reduction}, the desired estimates for $f\ast\nu_t^{k,l}$ and $f\ast\nu_t^{k,l,j}$ then follow once we establish the following uniform operator bounds.
If $r^\sharp<r$, we require, for $k\ge1$, $l\le k$, and $-2\le j\le l$,
\begin{align}\label{e:claim3}
    \sup_{\substack{\overline u\in\R^m \\ |\overline u^{(2)}|\sim2^{-j}}} \big\|\cG^{k,l}_{\overline u}\big\|_{L^2(\R^d)\to L^2(\R^{d+1})} &\lesssim
    2^{-k\frac{d-1}{2}} 2^{\frac12 l(d-r-2)} \big(\mathbf B_{j,l}(r^\sharp,r)\big)^{\frac12}, && j<l,\\\label{e:claim3-b}
    \sup_{\substack{\overline u\in\R^m \\ |\overline u^{(2)}|\lesssim2^{-l}}} \big\|\cG^{k,l}_{\overline u}\big\|_{L^2(\R^d)\to L^2(\R^{d+1})} &\lesssim 2^{-k\frac{d-1}{2}} 2^{\frac12 l(d-r^\sharp-2)}, && j=l.
\end{align}
If $r^\sharp=r$ and $l$ is large, it suffices to prove
\begin{equation}\label{e:claim4}
    \sup_{\substack{\overline u\in\R^m \\ |\overline u|\sim1}} \big\|\cG^{k,l}_{\overline u}\big\|_{L^2(\R^d)\to L^2(\R^{d+1})}\lesssim2^{-k\frac{d-1}{2}}2^{\frac12 l(d-r-2)}
\end{equation}
When $2^{-l}\sim 1$, proving \eqref{e:claim4} with the range $|\overline u|\sim 1$ replaced by $|\overline u|\lesssim 1$ is sufficient. See \eqref{e:Gu2jlarge} below. The implicit constants in these estimates are independent of $k$, $l$, and $j$.

To establish \eqref{e:claim3}--\eqref{e:claim4}, we apply Proposition \ref{prop:locall2-2} with $S=S^{\overline u}$. The following lemma provides the required singular value lower bounds under the assumption $\mathcal L(C_r,d_\perp\mid\mathcal B)=2r$.

\begin{lem}\label{lem:claim3}
    Let $\overline u\in \R^m$ such that $10^{-1} \le |\overline u|\le 10$. Let 
    $s_i=s_i(\overline u)$ be the ordered singular values of $S^{\overline u}$ (such that $s_k\ge s_{k+1}$ for $k=1,\dots, d-1$). Assume that $\cL(\cC_r, d_\perp\,|\,\cB) = 2r$.
    Then there exists a constant $C = C(\vec S)>0$ such that the following estimates hold.

    (i) $s_i(\overline u) \ge C$ for $1\le i\le r^\sharp$.

    (ii) We have
        \[
        \prod_{i=1}^r s_i(\overline u) \ge C |\overline u^{(2)}|^r.
        \]
    %\begin{itemize}
    %    \item [(i)] $s_i(\overline u) \ge C_S$ for $1\le i\le r^\sharp$.
    %    \item [(ii)] We have
    %    \[
    %    \prod_{i=1}^r s_i(\overline u) \ge C |\overline u^{(2)}|^r.
    %    \]
    %\end{itemize}
\end{lem}

\begin{proof}
    By definition of $r^\sharp$, the first $r^\sharp$ singular values $s_i = s_i(\theta)$, $1\le i\le r^\sharp$, of $S^{\theta}$ are nonzero, for all $10^{-1}\le |\theta|\le 10$. By \cite[Appendix D]{HornJohnson} or \cite{BhatiaElsnerKrause}, the singular values (i.e. the square roots of the eigenvalues of $S^\intercal S$) are continuous as functions of the matrix coefficients. Thus the $s_i$ are continuous as functions of $\theta$ and hence each $s_i$ has a minimum on $\{\theta:10^{-1}\le |\theta|\le 10\}$ which is positive for $1\le i\le r^\sharp$. Thus we have shown that there is $C = C(\vec S)>0$ such that $s_i(\theta)\ge C$ for $i=1,\dots r^\sharp$, for all $10^{-1}\le |\theta|\le 10$. This proves part $(i)$.

    Part $(ii)$ is a simple application of Vieta's formula. We consider the characteristic polynomial $p_{S,\theta}(\lambda)$ of $(S^\theta)^\intercal S^\theta$. By Vieta's formula, 
\[
\cC_r(\theta) = c\sum_{1\le i_1 < \dots < i_r\le d} \prod_{j=1}^r (s_{i_j}(\theta))^2
\]
where $c$ is either $-1$ or $1$. This implies the inequality
   \begin{align}\label{e:prod-cd-r}
       \prod_{i=1}^r (s_i(\theta))^2\le  |\cC_r(\theta)|\le \begin{pmatrix}
           d \\ r
       \end{pmatrix} \prod_{i=1}^r (s_i(\theta))^2.
   \end{align}
   The hypothesis $\mathcal L(C_r,d_\perp\,|\,\mathcal B)=2r$ gives $|\cC_r(\overline u)|\gtrsim (d_\perp(\overline u))^{2r}=|\overline u^{(2)}|^{2r}$. Combining this with the upper bound in \eqref{e:prod-cd-r} proves $(ii)$.
\end{proof}

\begin{rem}\label{rem:loja-2r}
We note that the condition $\cL(\cC_r, d_\perp\,|\,\cB) = 2r$ is fulfilled if $m=2$. Indeed, when $m=m^\sharp=2$, the statement is straightforward by the identity $\cC_r(\theta) = |\theta|^{2r}\cC_r(\theta/|\theta|)$. To see that $\cL(\cC_r, d_\perp\,|\,\cB) = 2r$ holds if $m=2$ and $m^\sharp = 1$, we denote $\theta = (\theta_1,\theta_2)\in \R\times \R$. Since $\cC_r$ is a homogeneous polynomial of degree $2r$,
\begin{align*}
    \cC_r(\theta_1,\theta_2) = \theta_2^{2r} \cC_r(\theta_1/\theta_2, 1)
\end{align*}
for $\theta\in \R^2$ with $\theta_2 \neq 0$. By our assumption, $\cC_r(\vartheta_1, 1)$ is either a positive or a negative polynomial in $\vartheta_1$. Thus,
\begin{align*}
    |\cC_r(\theta_1,\theta_2)| = \theta_2^{2r}|\cC_r(\theta_1/\theta_2, 1)|\ge C\theta_2^{2r}
\end{align*}
for $(\theta_1,\theta_2)\in \R\times (\R\setminus\{0\})$ with a constant $C>0$, which yields $\cL(\cC_r, d_\perp\,|\,\cB) = 2r$.
\end{rem}

We now explain how \eqref{e:claim3}--\eqref{e:claim4} follow from Proposition \ref{prop:locall2-2} and Lemma \ref{lem:claim3}. We first consider $2^{-l}\sim 1$. In this range $2^{-j}\sim 1$ as well, so it suffices to prove
\begin{align}\label{e:Gu2jlarge}
    \sup_{\overline u: |\overline u|\lesssim 1} \big\|\cG^{k,l}_{\overline u}\big\|_{L^2(\R^d)\to L^2(\R^{d+1})}\lesssim 2^{-\frac{d-1}{2}k}.
\end{align}
If $k=l$, then the estimate above is straightforward from Schur's test. Thus, we may assume $l<k$. We use a standard $TT^*$-argument. Recall that $P = \big(0\ I_{d-1}\big)$. Let $J^{\overline u}$, $\la^{\overline u}$ be a $d\times d$ skew-symmetric matrix and an element of $\R^d$, respectively, such that $S^{\overline u} = \big((J^{\overline u})^\intercal\ \la^{\overline u}\big)^\intercal$.
Then, by calculation, we get
\begin{align*}
\partial_{y'}^\intercal \partial_{x'}  \Psi^{\overline u}\big|_{(x,t,y)} = t^{-1}\sigma_{\overline u}(x',t,y_1) g''(\frac{x'-y'}{t}) + PJ^{\overline u} P^\intercal + PJ^{\overline u}e_1(g'(\frac{x'-y'}{t}))^\intercal.
\end{align*}
Since $|x'-y'|\le 2\eps$ for $(x,t,y)\in \supp(\chi_{\overline u})$, the last summand on the right-hand side is $O(\eps)$. On the other hand, the condition $2^{-l}\sim 1$ implies $|\sigma_{\overline u}(x',t,y_1)|\sim 1$ for all $(x,t,y)$ in the support of the amplitude function $\chi_{\overline u}(x,t,y)\zeta(2^lt\sigma_{\overline u}(x',t,y_1))$. Combining with the assumption that $g''(0)$ is negative-definite, one can show that $t^{-1}\sigma_{\overline u} g'' + PJ^{\overline u}P^\intercal$ is invertible by \cite[Lemma 5.3]{MuellerSeeger2004} for its proof. Choosing $\eps$ sufficiently small, we conclude that  $|\det(\partial_{y'}^\intercal \partial_{x'} \Psi^{\overline u})|\ge c$ on the support of $\chi_{\overline u}(x,t,y)\zeta(2^lt\sigma_{\overline u}(x',t,y_1))$, with a uniform constant $c>0$ in $\overline u$. From this, \eqref{e:Gu2jlarge} follows.

%From \eqref{i:mixedhessian}, one can derive that $|\det(\partial_y^\intercal \partial_{x,t}  \psi)|\ge c$ for $c>0$. Hence a standard argument for estimating oscillatory integral operator yields the bound $\|T^{k,l}\|_{L^2\to L^2}\lesssim \lambda^{-(d-1)/2}$.

%Therefore, Proposition \ref{e:locall2-2} implies \eqref{e:claim3} and \eqref{e:claim4} whenever $S^{\overline u}$ satisfies the lower bound stated in the proposition uniformly in $\overline u$. This is guaranteed by the following lemma.

Now we assume that $2^{-l}$ is small. Precisely, $2^{-l}$ should be so small that $\zeta_1(2(t^2\sigma_{\overline u}(x',t,y_1)^2 + t^4|\overline u|^2)^{1/2})$ vanishes unless $10^{-1}\le |\overline u|\le 10$, provided that $|\sigma_{\overline u}(x',t,y_1)|\le 2^{-l+3}$. Observe that the kernel of $\cG^{k,l}_{\overline u}$ is expressed identically to that of $T^{k,l}$, with $S$ and $\chi$ replaced by $S^{\overline u}$ and $\chi_{\overline u}$, respectively. Therefore, by the above discussion, we see that \eqref{e:claim3}--\eqref{e:claim4} are consequences of Proposition \ref{prop:locall2-2} and Lemma \ref{lem:claim3}. This finishes the proof of Proposition \ref{prop:locall2}.

\section{Proof of the oscillatory integral estimate}\label{sec:l2estspf} 
In this section, we prove Proposition \ref{prop:locall2-2}.
It turns out that under the assumption $r\le d-2$ we can derive optimal (or nearly optimal) bounds by freezing the variable $y_1$ and considering operators $\widetilde T^{k,l}_{y_1}$ mapping functions on $\bbR^{d-1}$ to functions on $\bbR^{d+1}$.
%The first step of the proof is to freeze the variable $y_1$. This is purposed to avoid dealing with the situation where $\Psi$ exhibits more complicated singularity behavior, by disregarding the contribution the variable $y_1$ can give. 
By the H\"older and Minkowski inequalities,
\begin{align}\label{e:freezing}
\|T^{k,l}f\|_{L^2(\R^d\times [1,2])}^2 \lesssim 2^{-l} \int \|\wt T^{k,l}_{y_1}f_{y_1}\|_{L^2(\R^d\times [1,2])}^2 dy_1    
\end{align}
for $f\in L^2(\bbR^{d})$, where $f_{y_1}(y') = f(y_1,y')$ and $\wt T^{k,l}_{y_1}$ is defined by
%an operator sending functions in $\R^{d-1}$ to functions in $\R^{d+1}$ given by
\[
\wt T^{k,l}_{y_1}g(x,t)  = \int_{\R^{d-1}} e^{i2^k\Psi(x,t,y)} \chi(x,t,y) \zeta(2^l t \sigma(x',t,y_1)) g(y') dy'
\]
for $g\in L^2(\R^{d-1})$.
From now on, we fix $y_1$ and suppress the dependence on $y_1$. That is, we write $\wt T^{k,l}_{y_1} = \wt T^{k,l}$, $\sigma(x',t,y_1) = \sigma(x',t)$, $\Psi(x,t,y) = \Psi(x,t,y')$, and $\chi(x,t,y) = \chi(x,t,y')$. By \eqref{e:freezing}, the proof would be completed once we show that for all $f\in \cS(\R^{d-1})$,
\begin{align}\label{e:lrtkl*}
    \|\wt T^{k,l}f\|_{L^2(\R^d\times [1,2])} \lesssim 2^{-k\frac{d-1}{2}} 2^{\frac12 l(d-r_2-1)} (\mathbf B_{j,l}(r_1,r_2))^{\frac12}\|f\|_{L^2(\R^{d-1})}
\end{align}
holds with the implicit constant independent of $y_1$ and $f$.

Establishing the estimate \eqref{e:lrtkl*} is the most technical part of the paper. We begin with a heuristic overview of the main ideas. For $w^\tau=(1-\tau)y'+\tau z'$, $0\le\tau\le1$, a direct calculation gives
\begin{align}
    \begin{aligned}\label{i:lrmh*}
    \partial_{x,t}\partial_{y'}^\intercal  \Psi(x,t, w^\tau) &= S P^\intercal + t^{-1}\sigma(x',t)\, \begin{pmatrix}
    0 \\ g''(\frac{x'-w^\tau}{t}) \\ 0
\end{pmatrix} \\
&\  - \begin{pmatrix}
    0 \\ 0 \\ t^{-2}\sigma(x',t)(x'-w^\tau)^\intercal g''(\frac{x'-w^\tau}{t})
\end{pmatrix} + Se_1 g'(\frac{x'-w^\tau}{t})^\intercal.
\end{aligned}
\end{align}
Here, the last two terms will be treated as error terms.

The proof relies on almost orthogonality and proceeds through two successive decompositions of the $y'$-space. The initial decomposition exploits the negative definiteness of $g''$ in \eqref{i:lrmh*}. At this stage, the skew-symmetry of $J$ allows us to obtain the relevant estimates without addressing the singular value structure of $SP^\intercal$. Heuristically, for the purpose of this decomposition, we may regard the mixed Hessian as $2^{-l} (0\ g''\ 0)^\intercal$. This leads to a decomposition of the $y'$-space into
cubes of side length $2^{-k+l}$.

We then perform a finer decomposition for each cube. Having accounted for the curvature of $g$ in the initial decomposition, we now focus on the singular values of $SP^\intercal$. Only singular values $\gtrsim 2^{-l}$ are relevant to this decomposition; see the definitions of $\mathcal N_1$ and $\mathcal N_2$ below.

An additional difficulty arises because we have frozen the coordinate $y_1$ and reduced the analysis to the $y'$-space. In this regard, we should understand how the direction $e_1$ interacts with $S$. This leads us to distinguish between the cases $|Se_1|\le2^{-l}$ and $|Se_1|>2^{-l}$. In the first case, the frozen direction does not significantly affect the decomposition. In the second case, however, we must take into consideration the position of $e_1$ relative to the subspace spanned by the right singular vectors of $S$ corresponding to singular values $\gtrsim 2^{-l}$. When the projection of $e_1$ onto this subspace is small, we first rescale the output variables $x,t$ so that the transformed vector $LSe_1$ has size $O(2^{-l})$. We then choose the subdivisions using the singular directions and singular values of the rescaled operator.

\subsection{Initial block decomposition}
Let $P$ denote the $(d-1)\times d$ matrix $\big(0\ \ \ I_{d-1}\big)$. Choose an orthonormal basis $\{\fe_i\}_{i=1}^{d-1}\subset \R^{d-1}$ so that $\fe_1$ is parallel to $PJe_1$. We choose $\fe_1$ arbitrarily on the unit sphere if $PJe_1 = 0$. We set, for $\fz\in \Z^{d-1}$, 
\begin{align*}
    \fq(\fz) := &\{y': \inn{c(\fz)}{\fe_i}\le \inn{y'}{\fe_i}<\inn{c(\fz_1+1,\dots,\fz_{d-1}+1)}{\fe_i},\ 1\le i\le d-1\}
\end{align*}
where $c(\fz)$ denotes a lower corner of $\fq(\fz)$ defined by
\[
c(\fz) = \sum_{i=1}^{d-1} 2^{-k+l} \fz_i \fe_i.
\]
Obviously, $\{\fq(\fz)\}_{\fz\in \Z^{d-1}}$ is a tiling of $\R^{d-1}$. Using this, we decompose
\[
\wt T^{k,l} f(x,t) = \sum_{\fz\in \Z^{d-1}} \wt T^{k,l}_{\fz} f(x,t) := \sum_{\fz\in \Z^{d-1}} \wt T^{k,l}[f\chi_{\fq(\fz)}](x,t).
\]
To exploit the almost orthogonality among $\wt T^{k,l}_\fz$, we use the Cotlar-Stein lemma of the form
\begin{align}\label{e:cotlar-stein}
    \big\|\sum_{\fz} \wt T^{k,l}_\fz\big\|_{L^2\to L^2} \lesssim \sup_\fz \sum_{\fz'} \|(\wt T^{k,l}_\fz)^*\wt T^{k,l}_{\fz'}\|_{L^2\to L^2}^{1/2} + \sup_{\fz} \sum_{\fz'} \|\wt T^{k,l}_\fz(\wt T^{k,l}_{\fz'})^*\|_{L^2\to L^2}^{1/2},
\end{align}
This inequality can be proved with a small modification of the proof in \cite[ch. VII.2]{Stein-harmonic}. 
By definition, $\wt T^{k,l}_\fz(\wt T^{k,l}_{\fz'})^* = 0$ if $\fz \neq \fz'$. The following proposition gives rapid decay of this operator norm in $|\fz-\fz'|$.

\begin{prop}\label{prop:estTzz'}
    Let $S$ be as in Proposition \ref{prop:locall2-2}. Let $\fz,\fz'\in \Z^{d-1}$ and $N\in \N$. Then we have
    \[
    \|(\wt T^{k,l}_\fz)^*\wt T^{k,l}_{\fz'}\|_{L^2\to L^2} \lesssim_N 2^{-k(d-1)}2^{l(d-r_2-1)} \mathbf B_{j,l}(r_1,r_2)(1+|\fz-\fz'|)^{-N}.
    \]
\end{prop}
Once the proposition is established, \eqref{e:lrtkl*} follows immediately from \eqref{e:cotlar-stein}. To prove Proposition \ref{prop:estTzz'}, we begin by showing that the kernel of $(\wt T^{k,l}_\fz)^*\wt T^{k,l}_{\fz'}$ can be expressed as an oscillatory integral where the amplitude decays rapidly as $|\fz-\fz'|$ grows.

\begin{lem}\label{lem:decayfzfz'}
    Let $\fz,\fz'\in \Z^{d-1}$. Then we have
    \begin{align}
    \begin{aligned}\label{i:Tz*Tz'}
        &(\wt T^{k,l}_\fz)^*\wt T^{k,l}_{\fz'}(y',z') \\
        &\quad= \chi_{q(\mathfrak z)}(y')\chi_{q(\mathfrak z')}(z')\int e^{i2^k(\Psi(x,t,z') - \Psi(x,t,y'))} \cA_{N,k,l}(x,t,y',z') dxdt,
    \end{aligned}
    \end{align}
    where $\cA_{N,k,l}$ is a function such that for $(y',z')\in \fq(\fz)\times \fq(\fz')$, $\cA_{N,k,l}$ is smooth in $x,t$,
    \begin{align*}
        &\supp(\cA_{N,k,l}) \\
        &\quad\subset \{(x,t,y',z'): t\sim 1,\ |x_1|\lesssim 1,\ |\sigma(x',t)|\lesssim 2^{-l},\  |x'|,\,|y'|,\,|z'|\le \eps\}
    \end{align*}
    with the small number $\eps>0$ in \eqref{r:chisupp}, and $\cA_{N,k,l}$ satisfies the following for all $N\in \N$:
    \begin{itemize}[topsep=2pt, itemsep=2pt]
        \item [i)] If $|Se_1|\le 2^{-l}$, then 
        \begin{align*}
        |\partial_{x,t}^\alpha\cA_{N,k,l}(x,t,y',z')|\lesssim_{\alpha,N} (1+|\fz-\fz'|)^{-N}
    \end{align*}
    for all $\alpha\in \N_0^{d+1}$. Here, $\N_0:=\N\cup\{0\}$.
    \item [ii)] If $|Se_1|\ge 2^{-l}$, then
    \begin{align*}
        |\inn{v}{\partial_{x,t}}^{M_1}\inn{\frac{Se_1}{|Se_1|}}{\partial_{x,t}}^{M_2}&\cA_{N,k,l}(x,t,y',z')| \\
        &\lesssim_{N,M_1,M_2} (2^l|Se_1|)^{M_2}(1+|\fz-\fz'|)^{-N}
    \end{align*}
    for all $M_1,M_2\in \N_0$ and unit vectors $v$ in an orthogonal complement of $\text{span}(Se_1)\subset \R^{d+1}$.
    \end{itemize}
\end{lem}

\begin{proof}
    Without loss of generality, we assume that the cutoff function $\chi$ is real-valued. We first consider the case $|\fz-\fz'|\le 10\sqrt d$ or $l=k$. In this case, we define 
    \begin{align*}
        \cA_{N,k,l}(x,t,y',z') = \chi(x,t,y')\chi(x,t,z') \big(\zeta(2^lt\sigma(x',t))\big)^2.
    \end{align*} 
    This function satisfies all the properties stated in the lemma. Hence, in the remainder of the proof, we assume that $|\fz-\fz'|>10\sqrt d$ and $l<k$. For $(y',z')\in \fq(\fz)\times \fq(\fz')$, we write
    \begin{align*}
        (\wt T^{k,l}_\fz)^*&\wt T^{k,l}_{\fz'}(y',z') \\
        &= \int e^{i2^k(\Psi(x,t,z') - \Psi(x,t,y'))} \chi(x,t,y')\chi(x,t,z') \big(\zeta(2^lt\sigma(x',t))\big)^2 dxdt.
    \end{align*}
    As before, we omit the dependence on $y_1$.

    We address the case $|Se_1|\le 2^{-l}$ first. In this case, the amplitude satisfies the bound
    \begin{align}\label{e:ampderivs1}
        \big|\partial_{x,t}^\alpha\big(\chi(x,t,y')\chi(x,t,z') \big(\zeta(2^lt\sigma(x',t))\big)^2\big)\big|\lesssim_\alpha 1
    \end{align}
    for all $\alpha\in \N_0^{d+1}$ because 
    \[
    \nabla_{x,t}\sigma(x',t) = Se_1
    \]
    and $2^l|Se_1|\le 1$ by assumption. To analyze the behavior of the phase function, we write
    \begin{align}\label{i:diffPsi}
        \Psi(x,t,y') - \Psi(x,t,z') = \int_0^1 \partial_{y'}^\intercal \Psi(x,t, w^\tau(y',z'))(y'-z') d\tau,
    \end{align}
where we denote $w^\tau(y',z') = w^\tau := (1-\tau) y' + \tau z'$.
Recall that $g'(0) = 0$, $g''(0)$ is negative-definite, and $|x'-y'|$, $|x'-z'|$ are sufficiently small. Thus, from the above identity, we notice that if
\[
v = \frac{1}{|y'-z'|}\begin{pmatrix}
    0 \\
    y'-z' \\
    0
\end{pmatrix}\in \R^{d+1},
\]
then
\begin{align}
    \begin{aligned}\label{e:lbub-vPsi}
        |\inn{v}{\partial_{x,t}\partial_{y'}^\intercal  \Psi(x,t,w^\tau)(y'-z')}|&\ge c2^{-l}|y'-z'|, \\
        |\inn{v}{\partial_{x,t}}^N \inn{\partial_{y'}  \Psi(x,t,w^\tau)}{y'-z'}|&\le C_N2^{-l}|y'-z'|
    \end{aligned}
\end{align}
for all $N\in \N$ with constants $c, C_N>0$. Here, we have used
\[
\inn{y'-z'}{PJP^\intercal(y'-z')} = 0.
\]
Now, let $L_1$ be a differential operator on $\R^{d+1}$ given by
\[
L_1 f(x,t) = \frac{\inn{v}{\partial_{x,t} f(x,t)}}{i2^k \inn{v}{\partial_{x,t}(\Psi(x,t,z') - \Psi(x,t,y'))}}.
\]
We define $\cA_{N,k,l}$ for the case $|Se_1|\le 2^{-l}$ by
\begin{align}\label{i:defANkl}
    \cA_{N,k,l}(x,t,y',z') = (L_1^\intercal)^N\big(\chi(x,t,y')\chi(x,t,z') \big(\zeta(2^lt\sigma(x',t))\big)^2\big).
\end{align}
By repeated integration by parts, it is immediate that $(\wt T^{k,l}_\fz)^*\wt T^{k,l}_{\fz'}(y',z')$ is expressed as in \eqref{i:Tz*Tz'}. Moreover, by definition, $\cA_{N,k,l}$ satisfies the support property stated in the lemma. Lastly, combining \eqref{e:lbub-vPsi} with the fact that $||y'-z'| - 2^{-k+l}|\fz-\fz'||\le \sqrt{d-1} 2^{-k+l}$, we derive the bound for $\cA_{N,k,l}$ stated in the assertion $i)$.

Next, we assume $|Se_1|>2^{-l}$. The proof of the desired statement is similar to that of \cite[Lemma 6.2]{RyuSeeger}. We divide the argument into two subcases and consider each separately:
\[
\textbf{I}.\ |\inn{\frac{y'-z'}{|y'-z'|}}{\fe_1}| > 2^{-l}|PJe_1|^{-1},\quad \textbf{I\!I}.\ |\inn{\frac{y'-z'}{|y'-z'|}}{\fe_1}| \le 2^{-l}|PJe_1|^{-1}.
\]

If $PJe_1 = 0$, then only the case \textbf{I\!I} applies. To address \textbf{I}, we compute $\partial_{x_1}\Psi(x,t,z') - \partial_{x_1}\Psi(x,t,y') = \inn{PJe_1}{y'-z'}$. Thus,
\begin{align}\label{e:lbdiff1}
    \big|\partial_{x_1}\big(\Psi(x,t,y') - \Psi(x,t,z')\big)\big| \ge C2^{-l}|y'-z'|.
\end{align}
Note that all the higher $x_1$-derivatives of $\Psi(x,t,z') - \Psi(x,t,y')$ vanish. Since $\inn{Je_1}{e_1} = 0$, the $x_1$-derivative of $\sigma(x',t)$ also vanishes, which yields
\begin{align}\label{e:ampderivs2}
    \big|\partial_{x_1}^N\big(\chi(x,t,y')\chi(x,t,z') \big(\zeta(2^lt\sigma(x',t))\big)^2\big)\big|\lesssim_N 1
\end{align}
for all $N\in \N$. Now we define
\[
\cA_{N,k,l}(x,t,y',z') = (L_2^\intercal)^N\big(\chi(x,t,y')\chi(x,t,z') \big(\zeta(2^lt\sigma(x',t))\big)^2\big),
\]
where $L_2 = (i2^k \partial_{x_1}(\Psi(x,t,z') - \Psi(x,t,y')))^{-1}\partial_{x_1}$. Using \eqref{e:lbdiff1} and \eqref{e:ampderivs2}, we see that $|\cA_{N,k,l}(x,t,y',z')|\lesssim (1+ 2^{k-l} |y'-z'|)^{-N}$ for all $N\in \N$. The derivative estimates in part $ii)$ follow in the same way, using $\inn{v}{\partial_{x,t}\sigma(x',t)} = 0$ for $v\perp Se_1$. This finishes the proof for the case \textbf{I}.

We now turn to the case \textbf{I\!I}. As before, we choose $v = |y'-z'|^{-1}(0, y'-z', 0)^\intercal\in \R^{d+1}$. Using \eqref{i:lrmh*}, we note
\begin{align*}
    \inn{v}{\partial_{x,t}\partial_{y'}^\intercal\Psi(x,t,w^\tau) (y'-z')} &= \frac{t^{-1}\sigma}{|y'-z'|} \inn{y'-z'}{g''(\frac{x'-w^\tau}{t})(y'-z')}\\
    &\qquad + \frac{\inn{PJe_1}{y'-z'}}{|y'-z'|}\inn{g'(\frac{x'-w^\tau}{t})}{y'-z'}.
\end{align*}
Since $g'(0) = 0$ and $|x'-w^\tau|$ is sufficiently small, the condition in \textbf{I\!I} implies that the second summand on the right-hand side is bounded by $\varepsilon 2^{-l} |y'-z'|$ for small $\varepsilon>0$. Meanwhile, the negative-definiteness of $g''$ implies that the absolute value of the first summand is $\ge C 2^{-l}|y'-z'|$. From these bounds and \eqref{i:diffPsi}, it follows that
\begin{align}\label{e:psidiff1}
    |\inn{v}{\partial_{x,t}(\Psi(x,t,y') - \Psi(x,t,z'))}|\gtrsim 2^{-l}|y'-z'|.
\end{align}

On the other hand, note that $|\inn{v}{\partial_{x,t}\sigma(x',t)}| \lesssim 2^{-l}$, a consequence of the assumption \textbf{I\!I}. This inequality implies
\[
|\inn{v}{\partial_{x,t}}^N\inn{\partial_{y'}\Psi(x,t,w^\tau)}{y'-z'}|\lesssim_N 2^{-l}|y'-z'|
\]
for $N\in \N$, which in turn yields
\begin{align}\label{e:psidiffub1}
    |\inn{v}{\partial_{x,t}}^N(\Psi(x,t,y') - \Psi(x,t,z'))|\lesssim_N 2^{-l}|y'-z'|.
\end{align}
Additionally, the bound for $\sigma$ gives that, for $N\in\N$,
\begin{align}\label{e:ampderivs3}
    \big|\inn{v}{\partial_{x,t}}^N\big(\chi(x,t,y')\chi(x,t,z') \big(\zeta(2^lt\sigma(x',t))\big)^2\big)\big|\lesssim_N 1.
\end{align}

We now define $A_{N,k,l}$ as in \eqref{i:defANkl}. Then, by \eqref{e:psidiff1}--\eqref{e:ampderivs3}, we have 
\[
|\cA_{N,k,l}(x,t,y',z')|\lesssim (1+ |\fz-\fz'|)^{-N}.
\]
Moreover, the derivatives 
\[
\inn{v}{\partial_{x,t}}^{M_1}\inn{\frac{Se_1}{|Se_1|}}{\partial_{x,t}}^{M_2}\cA_{N,k,l},\quad v\in \mathbb S^d \cap (\text{span}(Se_1))^\perp,\ M_1,M_2\in \N_0
\]
are controlled by $(2^l|Se_1|)^{M_2}(1+ |\fz-\fz'|)^{-N}$ up to a constant $C_{M_1,M_2,N}>0$ because $\inn{v}{\partial_{x,t} \sigma} = 0$ if $v\perp Se_1$. This proves $ii)$ for the case \textbf{I\!I}. 
\end{proof}

\subsection{A finer decomposition}
While the lemma validates a rapid decay of the kernel $(\wt T_\fz^{k,l})^*\wt T_{\fz'}^{k,l}(y',z')$ with respect to $|\fz-\fz'|$, it is insufficient to establish the desired estimate \eqref{e:lrtkl*} only using those decay estimates. We point out that the proof of Lemma \ref{lem:decayfzfz'} does not utilize the size condition on the singular values of $S$ in Lemma \ref{lem:claim3}. To exploit this, a finer decomposition is required. 

Throughout the decomposition step, we fix $\fz, \fz'$. Recall that $\{s_i\}_{i=1}^d$ denote singular values of $S$. Let $\{w_i\}_{i=1}^d$ be an orthonormal eigenbasis of $S^\intercal S$ such that $S^\intercal S w_i = s_i^2 w_i$ for $1\le i\le d$. The vectors $Sw_i$ are mutually orthogonal, with $|S w_i| = s_i$. %where $K$ is a sufficiently large number such that $K^{10}2^{-l}\le 1$.

To proceed, we consider the following two cases separately:
\[
a).\ |Se_1|\le 2^{-l},\quad b).\ |Se_1|> 2^{-l}.
\]
Let $\cN_1$ be a subspace of $\R^d$ generated by the vectors $w_i$ corresponding to singular values $s_i\ge 10^{3}2^{-l}$. We define $\cN_2$ in a similar fashion to $\cN_1$, but $s_i\ge 10^{3}2^{-l}$ is replaced by $s_i> 10^{-2}2^{-l}$. For a given subspace $V\subset \R^d$, let $\Pi_V$ be an orthogonal projection onto $V$. Fix a sufficiently small constant $\varrho>0$. Since the construction of the finer decomposition depends on the size of $|\Pi_{\cN_2}e_1|$, we further divide case $b)$ into the following two subcases:
\[
A.\ |\Pi_{\cN_2}e_1|\ge \varrho, \quad B.\ |\Pi_{\cN_2}e_1|< \varrho.
\]

\subsubsection{Decomposition when $|Se_1|\le 2^{-l}$}
Let $\mathbf M_1$ be a $(d+1)\times d$ matrix given by
\[
\mathbf M_1 = S\Pi_{(\text{span}(e_1))^\perp}.
\]
Let $\wt s_i$, $1\le i\le d$, be singular values of $\mathbf M_1$ such that $\wt s_1\ge \dots\ge \wt s_d\ge 0$. When $|Se_1|$ is smaller than $2^{-l}$, it is reasonable to expect that if $s_i\ge K_1 2^{-l}$ for some $i$ and large $K_1>0$, then $\wt s_i\ge K_2 2^{-l}$ for some $K_2>0$ with $K_1>K_2$. This expectation is indeed true.

\begin{lem}\label{lem:wtsibds}
    Assume $|Se_1|\le 2^{-l}$. Then, $\wt s_i\ge 10^2 2^{-l}$ for $1\le i\le \dim({\cN_1})$.
\end{lem}
\begin{proof}
Choose $1\le k\le \dim({\cN_1})$ so that $s_k\ge K_1 2^{-l}$ for some $K_1\ge 10^3$. Fix $k$. It is sufficient to show $\wt s_k\ge K_22^{-l}$ with a constant $K_2$ satisfying $K_2\ge 10^{-1}K_1$. Let $V_\circ = \text{span}(\{w_1,\dots,w_k\})$. Now we let $\wt V_\circ$ be a subspace of $\R^{d-1}$ given by $\wt V_\circ = \{Pz: z\in V_\circ\}$.
Note that $e_1 \notin {\cN_1}$, following from the definition of ${\cN_1}$. Thus, $\dim(\wt V_\circ) = \dim(V_\circ) = k$. By the assumption $|Se_1|\le 2^{-l}$, we obtain that, for $v\in V_\circ$ such that $|v| = 1$,
\begin{align*}
    K_12^{-l}\le |Sv|\le |Se_1||\inn{v}{e_1}| + |SP^\intercal P v| \le 2^{-l} + \frac{|SP^\intercal Pv|}{|Pv|},
\end{align*}
which implies that $|SP^\intercal w|\ge (K_1-1)2^{-l}$ for all unit vectors $w$ in $\wt V_\circ$. Thus, by the Courant--Fischer min-max theorem (see \cite{HornJohnson}), $\wt s_k\ge K_22^{-l}$ with $K_2 = K_1 -1$, which justifies our claim.
\end{proof}

%Combining the above observation with our assumptions \eqref{b:silbds} and $\vr\ge 10^32^{-l}$, we conclude
%\begin{align}\label{e:lbwtsi1}
%    \wt s_i\ge \begin{cases}
%        c & 1\le i\le r_1, \\
%        \vr/2, & r_1< i \le r_2,
%    \end{cases}
%\end{align}
%for some $c>0$. 
%, provided that $2^{-j}$ is sufficiently large compared to $2^{-l}$, precisely, when $C_S2^{-j}\ge 10^32^{-l}$, where $C_S$ is a constant appearing in \eqref{b:silbds}. If $C_S2^{-l}\le C_S2^{-j}<10^32^{-l}$, then the only lower bound we need is $\wt s_i\ge c > 10^2 2^{-l}$ with $c>0$ for $1\le i\le r$, which follows from the assumption that $2^{-l}$ is sufficiently small.

We now start decomposing the cubes $\fq(\fz)$, $\fq(\fz')$. % Denote $\pi_{e_1^\perp} = \pi_{(\text{span}(e_1))^\perp}$. 
Let $\wt \cN_1$ be a subspace of $\R^d$ generated by eigenvectors of $\mathbf M_1^\intercal \mathbf M_1$ of which eigenvalues exceed $10^42^{-2l}$.
Denote $r^{[1]} = \dim(\wt \cN_1)$. By Lemma \ref{lem:wtsibds}, $r^{[1]}\ge \dim({\cN_1})$. Let $(\wt w_i^{[1]})_{i=1}^{r^{[1]}}$ be orthonormal eigenvectors of $\mathbf M_1^\intercal \mathbf M_1$ such that $|S\wt w_i^{[1]}| = \wt s_i$ for $1\le i\le r^{[1]}$. For $\fn\in \Z^{r^{[1]}}$, we define
\begin{align*}
    \fr^{[1]}(\fz,\fz';\fn) = &\{Py\in \R^{d-1}:\ \inn{c^{[1]}(\fz,\fz';\fn)}{\wt w^{[1]}_i}\le \inn{y}{\wt w^{[1]}_i} \\
&\qquad\quad  < \inn{c^{[1]}(\fz,\fz';\fn_1+1,\dots,\fn_{r^{[1]}}+1)}{\wt w^{[1]}_i},\ 1 \le i\le r^{[1]}\},
\end{align*}
where
\[
c^{[1]}(\fz,\fz';\fn) = \sum_{i=1}^{r^{[1]}} 2^{-k} K (1+|\fz-\fz'|)  \wt s_i^{-1}\fn_i\wt w_i^{[1]}.
\]
Here, $K$ is a sufficiently large number to be chosen later.
Define
\[
G^{[1]}_{k,l,\fn,\fz,\fz'} f(x,t) = \wt T^{k,l}[f\chi_{\fq(\fz)\cap\fr^{[1]}(\fz,\fz';\fn)}].
\]
Throughout this argument, $\fz$ and $\fz'$ are fixed. We abbreviate $G^{[1]}_{k,l,\fn,\fz,\fz'}$ as $G^{[1]}_{\fn,\fz}$ and use analogous notation with $\fz$ and $\fz'$ interchanged. Then, $(\wt T^{k,l}_\fz)^* \wt T^{k,l}_{\fz'}$ is decomposed as
\[
(\wt T^{k,l}_\fz)^* \wt T^{k,l}_{\fz'} f(y') = \sum_{\fn,\fn'\in \Z^{r^{[1]}}} (G^{[1]}_{\fn,\fz})^*G^{[1]}_{\fn',\fz'} f(y').
\]
\begin{prop}\label{prop:estGnn'1}
    Let $S, \fz,\fz'$ be as in Proposition \ref{prop:estTzz'}. Assume that $|Se_1|\le 2^{-l}$. Let $\fn,\fn'\in \Z^{r^{[1]}}$ and $N, M\in \N$. Then
    \begin{align*}
        &\|(G^{[1]}_{\fn,\fz})^*G^{[1]}_{\fn',\fz'}\|_{L^2\to L^2} \\
        &\qquad \lesssim_{N,M} 2^{-k(d-1)}2^{l(d-1-r^{[1]})}\Big(\prod_{i=1}^{r^{[1]}}\wt s_i\Big)^{-1}(1+|\fz-\fz'|)^{-N}(1+|\fn-\fn'|)^{-M}.
    \end{align*}
\end{prop}
We postpone proving the proposition to the next subsection. Assuming it for now, we prove Proposition \ref{prop:estTzz'} in case a), then describe the decomposition for case A.

\begin{proof}[Proof of Proposition \ref{prop:estTzz'} when $|Se_1|\le 2^{-l}$]
    We use the Cotlar-Stein lemma \eqref{e:cotlar-stein} by which the problem is reduced to that of establishing
    \begin{align}
    \begin{aligned}\label{e:G4est}
        \sup_{\fn,\fn'\in\Z^{r^{[1]}}} \sum_{\wt\fn, \wt\fn'\in \Z^{r^{[1]}}}&\|(G^{[1]}_{\fn',\fz'})^*G^{[1]}_{\fn,\fz}(G^{[1]}_{\wt\fn,\fz})^*G^{[1]}_{\wt\fn',\fz'}\|_{L^2\to L^2}^{1/2} \\
        &\quad\lesssim_N 2^{-k(d-1)}2^{l(d-r_2-1)}\mathbf B_{j,l}(r_1,r_2)(1+|\fz-\fz'|)^{-N}
    \end{aligned}
    \end{align}
    for all $\fz,\fz'\in \Z^{d-1}$ and $N\in \N$, where the implicit constant depends on $N$. Fix $\fz,\fz'$ and $\fn,\fn'$. By definition, $(G^{[1]}_{\fn',\fz'})^*G^{[1]}_{\fn,\fz}(G^{[1]}_{\wt\fn,\fz})^*G^{[1]}_{\wt\fn',\fz'} = 0$ unless $\fn = \wt\fn$. Applying Proposition \ref{prop:estGnn'1} twice, we get
    \begin{align}
    \begin{aligned}\label{e:G4est-nsum}
        &\sum_{\wt\fn\in \Z^{r^{[1]}}}\|(G^{[1]}_{\fn',\fz'})^*G^{[1]}_{\fn,\fz}(G^{[1]}_{\wt\fn,\fz})^*G^{[1]}_{\wt\fn',\fz'}\|_{L^2\to L^2}^{1/2} \\
        &\qquad\quad\lesssim_{N,M} 2^{-k(d-1)}2^{l(d-1-r^{[1]})}\Big(\prod_{i=1}^{r^{[1]}}\wt s_i\Big)^{-1}\\
        &\hspace{24mm}\times (1+|\fz-\fz'|)^{-N}(1+|\fn-\fn'|)^{-M}(1+|\fn-\wt\fn'|)^{-M}.
    \end{aligned}
    \end{align}
    Recall $r^{[1]}\ge \dim({\cN_1})$. Moreover, by the argument in the proof of Lemma \ref{lem:wtsibds}, $\wt s_i\gtrsim s_i$ for $1\le i\le \dim({\cN_1})$. Thus
    \begin{align}\label{e:singsizecomp}
        \Big(\prod_{i=1}^{r^{[1]}}\wt s_i\Big)^{-1}\lesssim 2^{l(r^{[1]}-\dim({\cN_1}))} \Big(\prod_{i=1}^{\dim({\cN_1})} s_i\Big)^{-1}.
    \end{align}
    We now claim that
    \begin{align}\label{e:clm-pdsingsize}
        2^{l(d-1-\dim({\cN_1}))} \Big(\prod_{i=1}^{\dim({\cN_1})} s_i\Big)^{-1}\lesssim 2^{l(d-1-r_2)} \mathbf B_{j,l}(r_1,r_2).
    \end{align}
    To see this, we consider two cases $2^{l(r_2-r_1)} < 2^{jr_2}$ and $2^{jr_2}\le 2^{l(r_2-r_1)}$ separately. If $2^{l(r_2-r_1)} < 2^{jr_2}$, we use the assumption in \eqref{b:lowerboundsi} that $s_i\ge C$ for $1\le i\le r_1$. Using this, we get $(\prod_{i=1}^{\dim({\cN_1})} s_i)^{-1}\lesssim 2^{l(\dim({\cN_1})-r_1)}$, yielding \eqref{e:clm-pdsingsize}. When $2^{jr_2}\le 2^{l(r_2-r_1)}$, we show
    \[
    \Big(\prod_{i=1}^{\dim({\cN_1})} s_i\Big)^{-1}\lesssim 2^{l(\dim({\cN_1})-r_2)}\Big(\prod_{i=1}^{r_2} s_i\Big)^{-1} \lesssim 2^{l(\dim({\cN_1})-r_2)} 2^{jr_2}.
    \]
    Note that the second inequality is immediate from the assumption \eqref{b:lowerboundsi}. To prove the first inequality, we further distinguish two subcases $\dim({\cN_1})\ge r_2$ and $\dim({\cN_1})<r_2$. Both subcases can be addressed by using that $s_i\gtrsim 2^{-l}$ for $1\le i\le \dim({\cN_1})$ and $s_i\lesssim 2^{-l}$ for $\dim({\cN_1})< i\le d$. We omit the details. Hence, the claim follows.
    
    Substituting \eqref{e:singsizecomp} and \eqref{e:clm-pdsingsize} into \eqref{e:G4est-nsum} and summing over $\wt\fn'\in \Z^{r^{[1]}}$ and then taking the supremum over $\fn,\fn'\in \Z^{r^{[1]}}$, we obtain \eqref{e:G4est}.
\end{proof}

\subsubsection{Decomposition when $|\Pi_{\cN_2}e_1|\ge\varrho$} Recall that $\cN_2$ denotes a subspace of $\R^d$ generated by the vectors $w_i$ with singular values $s_i> 10^{-2}2^{-l}$. We assume $|Se_1|>2^{-l}$ and $|\Pi_{\cN_2}e_1|\ge\varrho$. 

Since $|Sv|\le 10^{-2}2^{-l}|v|$ for $v\in ({\cN_2})^\perp$, by the assumption $|Se_1|>2^{-l}$, $\Pi_{\cN_2} e_1$ does not vanish. Define $e_1^\cN$ to be a unit vector parallel to $\Pi_{\cN_2} e_1$. We let $\wt {\cN_2} = {\cN_2} \cap (\text{span}(e_1^\cN))^\perp$.
Denote $r^{[2]} = \dim(\wt {\cN_2})$ and $\Pi_{(e_1^\cN)^\perp} = \Pi_{(\text{span}(e_1^\cN))^\perp}$ for simplicity. By definition, $r^{[2]} = \dim(\cN_2) - 1$. Let $\mathbf M_2$, $\mathring{\mathbf M}_2$ be $(d+1)\times d$ matrices defined by
\[
\mathbf M_2 = \Pi_{S\cN_2}S\Pi_{\cN_2},\quad\mathring{\mathbf M}_2 = \Pi_{(Se_1^\cN)^\perp}\Pi_{S{\cN_2}}S\Pi_{\cN_2}\Pi_{(e_1^\cN)^\perp},
\]
where $S{\cN_2}$ is the image of ${\cN_2}$ under $S$. Let $(\bar s_i)_{i=1}^d$ be singular values of $\mathring{\mathbf M}_2$, ordered so that $\bar s_1\ge \dots\ge \bar s_d\ge 0$.
Let $\{\wt w_i^{[2]}\}_{i=1}^{r^{[2]}}$ denote a collection consisting of orthonormal eigenvectors of $(\mathring{\mathbf M}_2)^\intercal \mathring{\mathbf M}_2$ with nonzero eigenvalues. We note that $\wt w_1^{[2]},\dots, \wt w_{r^{[2]}}^{[2]}$ span $\wt\cN_2\subset \{(0,x')\in \R^d:x'\in\R^{d-1}\}$. For $\fn\in \Z^{r^{[2]}}$, let
\begin{align*}
    c^{[2]}(\fz,\fz';\fn) &= \sum_{i=1}^{r^{[2]}} 2^{-k} K (1+|\fz-\fz'|) \bar s_i^{-1}\fn_i\wt w_i^{[2]}.
\end{align*}
As above, $K$ is a sufficiently large number. Then, we define
\begin{align*}
    \fr^{[2]}(\fz,\fz';\fn) := &\{Py\in \R^{d-1}:\ \inn{c^{[2]}(\fz,\fz';\fn)}{\wt w^{[2]}_i}\le \inn{y}{\wt w^{[2]}_i} \\
&\qquad\quad  < \inn{c^{[2]}(\fz,\fz';\fn_1+1,\dots,\fn_{r^{[2]}}+1)}{\wt w^{[2]}_i},\ 1\le i\le r^{[2]}\}
\end{align*}
and
\[
G^{[2]}_{\fn,\fz} f(x,t) \equiv G^{[2]}_{k,l,\fn,\fz,\fz'} f(x,t) := \wt T^{k,l}[f\chi_{\fq(\fz)\cap\fr^{[2]}(\fz,\fz';\fn)}].
\]
Using this setup, we obtain the decomposition
\[
(\wt T^{k,l}_\fz)^* \wt T^{k,l}_{\fz'} f(y') = \sum_{\fn,\fn'\in \Z^{r^{[2]}}} (G^{[2]}_{\fn,\fz})^*G^{[2]}_{\fn',\fz'} f(y').
\]

The following is a counterpart to Proposition \ref{prop:estGnn'1}.

\begin{prop}\label{prop:estGnn'2}
    Let $S, \fz,\fz'$ be as in Proposition \ref{prop:estTzz'}. Let $\fn,\fn'\in \Z^{r^{[2]}}$ and $N, M\in \N$. Assume that $|Se_1| > 2^{-l}$ and $|\fn-\fn'|>10\sqrt d$. Then
    \begin{align*}
        &\|(G^{[2]}_{\fn,\fz})^*G^{[2]}_{\fn',\fz'}\|_{L^2\to L^2} \\
        &\  \lesssim_{N,M} 2^{-k(d-1)}2^{l(d-2-r^{[2]})} \Big(|Se_1|\prod_{i=1}^{r^{[2]}}\bar s_i\Big)^{-1}(1+|\fz-\fz'|)^{-N}|\fn-\fn'|^{-M}.
    \end{align*}
\end{prop}

Postponing the proof of this proposition, we complete the proof of Proposition \ref{prop:estTzz'} after assuming that Proposition \ref{prop:estGnn'2} is valid. For that, we need an elementary result as follows.

\begin{lem}\label{lem:elem}
    We have the identity
    \[
    |Se_1^\cN| \prod_{i=1}^{r^{[2]}} \bar s_i = \prod_{i=1}^{\dim(\cN_2)} s_i.
    \]
\end{lem}

\begin{proof}
    Decompose 
    \[
    \R^{d} = \text{span}(e_1^\cN) + \wt{\cN_2} + ({\cN_2})^\perp,\quad \R^{d+1} = \text{span}(Se_1^\cN) + \wt\cN_{2,S} + (S{\cN_2})^\perp
    \]
    where $\wt\cN_{2,S} = S{\cN_2}\cap (\text{span}(Se_1^\cN))^\perp$. Since the restriction of $S$ to ${\cN_2}$ is injective, $\dim(\wt\cN_{2,S}) = r^{[2]}$. Additionally, note that $S({(\cN_2})^\perp)\subset (S{\cN_2})^\perp$, because $(\cN_2)^\perp$ is an invariant subspace for $S^\intercal S$. Let $(\bar w_i)_{i=1}^{r^{[2]}}$ be an orthonormal basis of $\wt\cN_{2,S}$. Now define the orthogonal matrices $R_1$, $R_2$ on $\R^{d+1}$, $\R^d$ respectively, such that $R_1 e_1 = Se_1^\cN/|Se_1^\cN|$, $R_1 e_i = \bar w_{i-1}$ for $2\le i\le r^{[2]}+1$, $R_1(\text{span}(\{e_i\}_{i=r^{[2]}+2}^{d+1})) = (S{\cN_2})^\perp$, and $R_2 e_1 = e_1^\cN$, $R_2 e_i = \wt w^{[2]}_{i-1}$ for $2\le i\le r^{[2]}+1$, $R_2(\text{span}(\{e_i\}_{i=r^{[2]}+2}^{d})) = ({\cN_2})^\perp$.

Then, a calculation gives
\begin{align*}
    R_1^\intercal \mathbf M_2 R_2 = \begin{pmatrix}
        \mathbf N & 0 \\
        0 & 0
    \end{pmatrix} \quad \text{where}\quad \mathbf N = \begin{pmatrix}
        |Se_1^\cN| &\ast \\
        0 & \mathring{\mathbf N}
    \end{pmatrix},
\end{align*}
and $\mathring{\mathbf N}$ is the $r^{[2]}\times r^{[2]}$ minor of $R_1^\intercal \mathring{\mathbf M}_2 R_2$ obtained by deleting all zero rows and columns of $R_1^\intercal \mathring{\mathbf M}_2 R_2$. From the above identity, it follows that $\det(\mathbf N) = |Se_1^\cN|\det(\mathring{\mathbf N})$. By construction of $R_1$ and $R_2$, it follows that
\[
\prod_{i=1}^{\dim(\cN_2)} s_i = |\det(\mathbf N)|,\quad \prod_{i=1}^{r^{[2]}} \bar s_i = |\det(\mathring{\mathbf N})|,
\]
implying the identity in the lemma, as desired.
\end{proof}

\begin{proof}[Proof of Proposition \ref{prop:estTzz'} for case A]
    For $|\fn-\fn'|\le 10\sqrt d$, the required estimate follows from Lemma \ref{lem:decayfzfz'} and Schur's test. Hence, we assume $|\fn-\fn'|> 10\sqrt d$. Since $|\Pi_{\cN_2}e_1|\ge\varrho$, the estimate $|Se_1|\sim |Se_1^\cN|$ holds.
    Combining this with Proposition \ref{prop:estGnn'2} and Lemma \ref{lem:elem}, we obtain
    \begin{align*}
        &\|(G^{[2]}_{\fn,\fz})^*G^{[2]}_{\fn',\fz'}\|_{L^2\to L^2} \\
        &\ \lesssim_{N,M} 2^{-k(d-1)}2^{l(d-1-\dim({\cN_2}))}\Big(\prod_{i=1}^{\dim({\cN_2})} s_i\Big)^{-1}(1+|\fz-\fz'|)^{-N}(1+|\fn-\fn'|)^{-M}.
    \end{align*}
    We also use that $r^{[2]} = \dim(\cN_2) - 1$. By the same argument as in the proof of \eqref{e:clm-pdsingsize}, we get
    \[
    2^{l(d-1-\dim({\cN_2}))}\Big(\prod_{i=1}^{\dim({\cN_2})} s_i\Big)^{-1}\lc 2^{l(d-r_2-1)}\mathbf B_{j,l}(r_1,r_2).
    \]
    Combining this with the $L^2$ estimate above, we get
    \begin{align*}
        &\|(G^{[2]}_{\fn,\fz})^*G^{[2]}_{\fn',\fz'}\|_{L^2\to L^2} \\
        &\quad \lesssim_{N,M} 2^{-k(d-1)}2^{l(d-r_2 -1)}\mathbf B_{j,l}(r_1,r_2)(1+|\fz-\fz'|)^{-N}(1+|\fn-\fn'|)^{-M}.
    \end{align*}
    Using the Cotlar-Stein lemma in the same manner as in the proof of Proposition \ref{prop:estTzz'} for case $a)$ gives the desired estimate for $(\wt T^{k,l}_\fz)^*\wt T^{k,l}_{\fz'}$. We omit the details.
\end{proof}

\subsection{Proofs of Propositions \ref{prop:estGnn'1} and \ref{prop:estGnn'2}} Our aim is to prove
\begin{align}\label{e:kerg1g1'}
    |(G^{[1]}_{\fn,\fz})^*G^{[1]}_{\fn',\fz'}(y',z')|&\lesssim_{N,M} (1+|\fz-\fz'|)^{-N}(1+|\fn-\fn'|)^{-M}, \\\label{e:kerg2g2'}
    |(G^{[2]}_{\fn,\fz})^*G^{[2]}_{\fn',\fz'}(y',z')|&\lesssim_{N,M} 2^{-l}|Se_1|^{-1}(1+|\fz-\fz'|)^{-N}(1+|\fn-\fn'|)^{-M}.
\end{align}
These pointwise kernel estimates imply Propositions \ref{prop:estGnn'1} and \ref{prop:estGnn'2} by Schur's test.

To show \eqref{e:kerg1g1'} and \eqref{e:kerg2g2'}, we may assume that $|\fn-\fn'|>10\sqrt d$. Otherwise, these estimates immediately follow from Lemma \ref{lem:decayfzfz'}. We first prove \eqref{e:kerg1g1'}. Assume that $y'\in \fq(\fz)\cap\fr^{[1]}(\fz,\fz';\fn)$, $z'\in \fq(\fz')\cap\fr^{[1]}(\fz,\fz';\fn')$, and $|Se_1|\le 2^{-l}$. Then,
\begin{align}\label{e:disty'z'}
    |y'-z'|\le C_d2^{-k+l}(1+|\fz-\fz'|)
\end{align}
with a constant $C_d>0$. We now choose
\begin{align*}
    v^{[1]} = \dfrac{S\Pi_{\wt\cN_1}P^\intercal(y'-z')}{|S\Pi_{\wt\cN_1}P^\intercal(y'-z')|}.
\end{align*}
Recalling \eqref{i:lrmh*}, we compute
\begin{align*}
    \inn{v^{[1]}}{\partial_{x,t}\partial_{y'}^\intercal \Psi(x,t,w^\tau)&(y'-z')} \\
    &= |S\Pi_{\wt\cN_1}P^\intercal(y'-z')| + \cR_1(y',z') + \cQ_1(x,t,y',z',\tau),
\end{align*}
where
\begin{align*}
    \cR_1(y',z') &= \inn{v^{[1]}}{S\Pi_{(\wt\cN_1)^\perp}P^\intercal(y'-z')},
\end{align*}
and
\begin{align*}
    &\cQ_1(x,t,y',z',\tau) \\
    &\quad\qquad= t^{-1}\sigma(x',t)\inn{v^{[1]}}{\fg(x,t,y',z',\tau)} + \inn{Se_1}{v^{[1]}}\inn{g'(\frac{x'-w^\tau}{t})}{y'-z'}.
\end{align*}
Here, the function $\fg$ is defined by
\[
\fg(x,t,y',z',\tau) = \begin{pmatrix}
    0 \\
    g''(\frac{x'-w^\tau}{t})(y'-z') \\
    -t^{-1}(x'-w^\tau)^\intercal g''(\frac{x'-w^\tau}{t})(y'-z')
\end{pmatrix}.
\]
The geometry of the rectangles gives
\begin{align}
\begin{aligned}\label{e:main-gigi'-c1}
    |S\Pi_{\wt\cN_1}P^\intercal(y'-z')| \ge 2^{-k} C_\circ^{[1]} K (1+|\fz-\fz'|)|\fn-\fn'|,
\end{aligned}
\end{align}
with $C_\circ^{[1]}>0$ and a sufficiently large $K>0$. On the other hand, by definition,
\[
    |S\Pi_{(\wt\cN_1)^\perp}P^\intercal(y'-z')| = |\mathbf M_1\Pi_{(\wt\cN_1)^\perp}P^\intercal(y'-z')|\le 2^{-l+10}|y'-z'|.
\]
Combining this with \eqref{e:disty'z'}, we get $|S\Pi_{(\wt\cN_1)^\perp}P^\intercal(y'-z')| \lesssim C_d2^{-k+10}(1+|\fz-\fz'|)$, implying $|\cR_1(y',z')|\le C_d2^{-k+10}(1+|\fz-\fz'|)$. Moreover,
\[
|\cQ_1(x,t,y',z',\tau)|\le C_g^{[1]} 2^{-k}(1+|\fz-\fz'|)
\]
with a constant $C_g^{[1]}>0$ depending on $g$. Here, we use the bounds $|Se_1|\le 2^{-l}$ and $|\sigma(x',t)|\lesssim 2^{-l}$. Now choose $K>0$ so large that 
\[
C_\circ^{[1]}K>2\max(2^{10}C_d, C_g^{[1]}).
\]
Then, combining all the bounds for $\cR_1, \cQ_1$, and \eqref{e:main-gigi'-c1} together with \eqref{i:diffPsi} gives
\begin{align}\label{e:lbsPsidiff3-c1}
    \big|\inn{v^{[1]}}{\partial_{x,t}}(\Psi(x,t,y') - \Psi(x,t,z'))\big|\gtrsim 2^{-k} K (1+|\fz-\fz'|)|\fn-\fn'|.
\end{align}
On the other hand, using $|Se_1|\le 2^{-l}$, \eqref{i:lrmh*}, and \eqref{i:diffPsi} together with $|y'-z'|\lesssim 2^{-k}2^l(1+|\fz-\fz'|)$ yields
\begin{align}\label{e:ubsPsidiff3}
    \big|\inn{v^{[1]}}{\partial_{x,t}}^M(\Psi(x,t,y') - \Psi(x,t,z'))\big|\lesssim_M 2^{-k} (1+|\fz-\fz'|),\quad M\ge 2.
\end{align}
Lastly, combining Lemma \ref{lem:decayfzfz'} with the assumption $|Se_1|\le 2^{-l}$, we get
\begin{align}\label{e:ubsamp3}
    \big|\inn{v^{[1]}}{\partial_{x,t}}^M\cA_{N,k,l}(x,t,y',z')\big|\lesssim_{N,M} (1+|\fz-\fz'|)^{-N},\quad N,M\ge 1.
\end{align}
Putting together the estimates \eqref{e:lbsPsidiff3-c1}--\eqref{e:ubsamp3} and applying integration by parts gives \eqref{e:kerg1g1'}. This completes the proof of Proposition \ref{prop:estGnn'1}.

We proceed to the proof of Proposition \ref{prop:estGnn'2} assuming $|Se_1|>2^{-l}$. 
%To this end, we investigate the following two cases separately: $(\fn,\fn')\in \fI_1$, $(\fn,\fn')\in \fI_2$, where
%\begin{align*}
%    \fI_1 &= \{(\fn,\fn'):  |\fn_{\dim({\cN_2})}-\fn'_{\dim({\cN_2})}| \le \mathbf C (1+|\fn-\fn'|_{r^{[2]}})\}, \\
%    \fI_2 &= \{(\fn,\fn'):  |\fn_{\dim({\cN_2})}-\fn'_{\dim({\cN_2})}| > \mathbf C (1+|\fn-\fn'|_{r^{[2]}})\}.
%\end{align*}
%Here, $|\fn|_{r^{[2]}}$ denotes the size of the vector $(\fn_1,\dots,\fn_{r^{[2]}},0)\in \Z^{\dim({\cN_2})}$, and $\mathbf C$ is a large constant which will be determined later. 
To prove \eqref{e:kerg2g2'}, it suffices to show that
\begin{align}\label{e:fGbd}
    \fG_{N,\fz,\fz',\fn,\fn'} \lesssim_M |\fn-\fn'|^{-M}
\end{align}
for any $M\in \N$, where
\begin{align*}
    \fG_{N,\fz,\fz',\fn,\fn'} = 2^l|Se_1|(1+|\fz-\fz'|)^{N}\sup_{\substack{y'\in \fq(\fz)\cap\fr^{[2]}(\fz,\fz';\fn)\\ z'\in \fq(\fz')\cap\fr^{[2]}(\fz,\fz';\fn')}} |(G^{[2]}_{\fn,\fz})^*G^{[2]}_{\fn',\fz'}(y',z')|.
\end{align*}
Fix $y'\in \fq(\fz)\cap\fr^{[2]}(\fz,\fz';\fn), z'\in \fq(\fz')\cap\fr^{[2]}(\fz,\fz';\fn')$. Let
\begin{align*}
    v^{[2]} = \frac{\Pi_{(Se_1^\cN)^\perp}S\Pi_{\wt\cN_2}P^\intercal(y'-z')}{|\Pi_{(Se_1^\cN)^\perp}S\Pi_{\wt\cN_2}P^\intercal(y'-z')|}.
\end{align*}
We note that
\begin{align}
    \begin{aligned}\label{e:y'-z'-c21}
        |y'-z'|&\le C_d2^{-k+l} (1+|\fz-\fz'|), \\
        |\Pi_{(Se_1^\cN)^\perp}S\Pi_{\wt\cN_2}P^\intercal(y'-z')|&\ge 2^{-k} C_\circ^{[2]} K (1+|\fz-\fz'|)|\fn-\fn'|
    \end{aligned}
\end{align}
with constants $C_d, C_\circ^{[2]}>0$.
As in the previous argument, we write
\begin{align*}
    \inn{v^{[2]}}{\partial_{x,t}\partial_{y'}^\intercal &\Psi(x,t,w^\tau)(y'-z')} \\
    &= |\Pi_{(Se_1^\cN)^\perp}S\Pi_{\wt\cN_2}P^\intercal(y'-z')| + \cR_2(y',z') + \cQ_2(x,t,y',z',\tau),
\end{align*}
where $\cR_2$ and $\cQ_2$ are defined as in $\cR_1$ and $\cQ_1$, but $v^{[2]}, \wt\cN_2$ instead of $v^{[1]}, \wt\cN_1$.
We note that every $u\in (\wt\cN_2)^\perp$ can be written as $u = \wt u + a e_1^\cN$. This implies 
\[
|\Pi_{(Se_1^\cN)^\perp}S\Pi_{(\wt\cN_2)^\perp}P^\intercal(y'-z')| \le 2^{-l}|y'-z'|.
\]
From this and \eqref{e:y'-z'-c21}, it follows that $|\cR_2(y',z')|\le C_d2^{-k}(1+|\fz-\fz'|)$. We also note that
\[
\inn{Se_1}{v^{[2]}} = \inn{S(e_1-\Pi_{\cN_2}e_1)}{v^{[2]}} + \inn{S\Pi_{\cN_2}e_1}{v^{[2]}},
\]
which consequently yields
\begin{align}\label{e:innprodse1v}
    |\inn{Se_1}{v^{[2]}}|\le 2^{-l}.
\end{align}
By \eqref{e:y'-z'-c21}, \eqref{e:innprodse1v}, and $|\sigma(x',t)|\lesssim 2^{-l}$, we get 
\[
|\cQ_2(x,t,y',z',\tau)|\le C_g^{[2]} 2^{-k}(1+|\fz-\fz'|)
\]
with a constant $C_g^{[2]}>0$. Choose $K>0$ such that $C_\circ^{[2]}K>2\max(C_d,C_g^{[2]})$.

Using the second estimate in \eqref{e:y'-z'-c21} together with the estimates for $\cR_2$, $\cQ_2$ above, we have
\begin{align}\label{e:lbsPsidiff3}
    \big|\inn{v^{[2]}}{\partial_{x,t}}(\Psi(x,t,y') - \Psi(x,t,z'))\big|\gtrsim 2^{-k} K (1+|\fz-\fz'|)|\fn-\fn'|.
\end{align}
We also use \eqref{i:diffPsi}. Next, note that $\inn{v}{\partial_{x,t}}\sigma(x',t) = \inn{Se_1}{v}$. Combining this identity, \eqref{i:diffPsi}, \eqref{e:innprodse1v} along with $|y'-z'|\lesssim 2^{-k+l}(1+|\fz-\fz'|)$, we obtain
\begin{align}\label{e:ubsPsidiff3-c21}
    \big|\inn{v^{[2]}}{\partial_{x,t}}^M(\Psi(x,t,y') - \Psi(x,t,z'))\big|\lesssim_M 2^{-k} (1+|\fz-\fz'|),\quad M\ge 2.
\end{align}
To control the derivatives of the amplitude function $\cA_{N,k,l}$, we write
\begin{align}\label{i:v2-se1}
    v^{[2]} = \Pi_{(Se_1)^\perp} v^{[2]} + \frac{1}{|Se_1|^2}\inn{Se_1}{v^{[2]}}Se_1.
\end{align}
By \eqref{e:innprodse1v} and Lemma \ref{lem:decayfzfz'}, we get
\begin{align}\label{e:ubsamp3-c21}
    \big|\inn{v^{[2]}}{\partial_{x,t}}^M\cA_{N,k,l}(x,t,y',z')\big|\lesssim_{N,M} (1+|\fz-\fz'|)^{-N}, \quad N,M\ge 1.
\end{align}
Combining \eqref{e:lbsPsidiff3}--\eqref{e:ubsamp3-c21} and integrating by parts yields \eqref{e:kerg2g2'}, as desired.

\subsection{Rescaling and decomposition in case B} We complete the proof of Proposition \ref{prop:estTzz'} in this subsection by considering the remaining case where $|\Pi_{\cN_2}e_1|<\varrho$ and $|Se_1|>2^{-l}$.
Let $\{\bar\fe_i\}_{i=1}^{d+1}$ be a set of orthonormal vectors in $\R^{d+1}$ such that $\bar\fe_1$ is parallel to $Se_1^\cN$. Consider a $(d+1)\times (d+1)$ rotation matrix $R$ such that $R\bar\fe_i = e_i$ for $1\le i\le d+1$. Let $D$ be a $(d+1)\times (d+1)$ diagonal matrix whose entries are $2^{-l}|S\Pi_{\cN_2}e_1|^{-1}$, $1,\dots, 1$. We then define
\[
 L =  R^\intercal D R.
\]
Let $\{\breve w_i\}_{i=1}^{\dim(\mathcal N_2)}$ be an orthonormal collection of right singular vectors of $\mathbf M_3 = LS\Pi_{\mathcal N_2}$ with nonzero singular values. Denote these singular values by $\breve s_i$, $1\le i\le \dim(\mathcal N_2)$. Note that
\begin{align}\label{i:prodcomp}
    2^{-l}|S\Pi_{\cN_2}e_1|^{-1}\prod_{i=1}^{\dim(\mathcal N_2)} s_i =  \prod_{i=1}^{\dim(\mathcal N_2)} \breve s_i.
\end{align}
Indeed, $L$ acts as the identity on $(\cN_2)^\perp$, so its volume factor is exactly $2^{-l}|S\Pi_{\mathcal N_2}e_1|^{-1}$. This gives \eqref{i:prodcomp}.
For $\fn\in \Z^{\dim(\mathcal N_2)}$, setting
\begin{align*}
    c^{[3]}(\fz,\fz';\fn) &= \sum_{i=1}^{\dim(\cN_2)} 2^{-k} K (1+|\fz-\fz'|) \breve s_i^{-1}\fn_i\breve w_i
\end{align*}
with $K$ a sufficiently large number, we define
\begin{align*}
    \fr^{[3]}(\fz,\fz';\fn) := &\{y'\in \R^{d-1}:\ \inn{c^{[3]}(\fz,\fz';\fn)}{\breve w_i}\le \inn{\Pi_{\mathcal N_2}P^\intercal y'}{\breve w_i} \\
&\   < \inn{c^{[3]}(\fz,\fz';\fn_1+1,\dots,\fn_{\dim({\cN_2})}+1)}{\breve w_i},\ 1\le i\le \dim(\cN_2)\}
\end{align*}
and
\[
G^{[3]}_{\fn,\fz} f(x,t) \equiv G^{[3]}_{k,l,\fn,\fz,\fz'} f(x,t) := \wt T^{k,l}[f\chi_{\fq(\fz)\cap\fr^{[3]}(\fz,\fz';\fn)}].
\]
Similarly, we have
\[
(\wt T^{k,l}_\fz)^* \wt T^{k,l}_{\fz'} f(y') = \sum_{\fn,\fn'\in \Z^{\dim({\cN_2})}} (G^{[3]}_{\fn,\fz})^*G^{[3]}_{\fn',\fz'} f(y').
\]

\begin{prop}\label{prop:estGnn'3}
    Let $\fn,\fn'\in \Z^{\dim(\cN_2)}$. Assume that $|Se_1|>2^{-l}$, $|\Pi_{\mathcal N_2}e_1| < \varrho$, and $|\fn-\fn'|>10\sqrt d$. Let $N,M\in \N$. Then
    \begin{align*}
        &\|(G^{[3]}_{\fn,\fz})^*G^{[3]}_{\fn',\fz'}\|_{L^2\to L^2} \\
        &\  \lesssim_{N,M} 2^{-k(d-1)}2^{l(d-2-\dim(\cN_2))} |Se_1|^{-1}\Big(\prod_{i=1}^{\dim(\cN_2)}\breve s_i\Big)^{-1}(1+|\fz-\fz'|)^{-N}|\fn-\fn'|^{-M}.
    \end{align*}
\end{prop}

From the assumption $|Se_1|\ge 2^{-l}$, the inequality 
\begin{align}\label{e:Se1Picomp}
    2^{-1}|Se_1|\le |S\Pi_{\cN_2}e_1|\le 2|Se_1|
\end{align} 
follows. Combining this with \eqref{i:prodcomp} and Proposition \ref{prop:estGnn'3}, we have
\begin{align}
\begin{aligned}\label{e:G*G3-ref}
    &\|(G^{[3]}_{\fn,\fz})^*G^{[3]}_{\fn',\fz'}\|_{L^2\to L^2} \\
        &\  \lesssim_{N,M} 2^{-k(d-1)}2^{l(d-1-\dim(\cN_2))} \Big(\prod_{i=1}^{\dim(\cN_2)} s_i\Big)^{-1}(1+|\fz-\fz'|)^{-N}|\fn-\fn'|^{-M}.
\end{aligned}
\end{align}
for $\fn,\fn'\in \Z^{\dim(\cN_2)}$ such that $|\fn-\fn'|>10\sqrt d$. The remaining case $|\fn-\fn'|\le 10\sqrt d$ can be addressed via Schur's test. To do so, we need to use the bound
\begin{align}
\begin{aligned}\label{e:ubfqfr3}
    &|\fq(\fz)\cap \fr^{[3]}(\fz,\fz';\fn)| \\
    &\quad\lesssim 2^{-k(d-1)} 2^{l(d-\dim(\cN_2))} (1+|\fz-\fz'|)^{\dim(\cN_2)} |S\Pi_{\cN_2}e_1| \big(\prod_{i=1}^{\dim(\cN_2)} s_i\big)^{-1}
\end{aligned}
\end{align}
which follows from \eqref{i:prodcomp} and \eqref{e:Se1Picomp}. In applying Schur's test, we use the kernel estimate with $N+\dim(\cN_2)$ in place of $N$. Consequently, for $|\fn-\fn'|\le 10\sqrt d$, we get a similar estimate to \eqref{e:G*G3-ref} without the factor $|\fn-\fn'|^{-M}$.
Applying the Cotlar–Stein lemma, arguing as in the previous two cases, we obtain the desired bound for $(\wt T^{k,l}_\fz)^*\wt T^{k,l}_{\fz'}$ in Proposition \ref{prop:estTzz'}. We omit the details.

We proceed to the proof of Proposition \ref{prop:estGnn'3}. By \eqref{e:ubfqfr3} and Schur's test, it is sufficient to prove
\begin{align}\label{e:G*G3}
    |(G^{[3]}_{\fn,\fz})^*G^{[3]}_{\fn',\fz'}(y',z')|\lesssim_{N,M} 2^{-l}|Se_1|^{-1}(1+|\fz-\fz'|)^{-N}(1+|\fn-\fn'|)^{-M}
\end{align}
for $N,M\in \N$. Fix $y'\in \fq(\fz)\cap\fr^{[3]}(\fz,\fz';\fn)$, $z'\in \fq(\fz')\cap\fr^{[3]}(\fz,\fz';\fn')$. Note that $L\bar\fe_1 = 2^{-l}|S\Pi_{\cN_2}e_1|^{-1}\bar\fe_1$ and $L\bar\fe_i = \bar\fe_i$ for $2\le i\le d+1$. Thus, by \eqref{i:Tz*Tz'} and a change of variables 
\[
\begin{pmatrix}
    x \\ t
\end{pmatrix} =  L\begin{pmatrix}
    u \\ s
\end{pmatrix}
\]
the kernel $(G^{[3]}_{\fn,\fz})^*G^{[3]}_{\fn',\fz'}(y',z')$ is expressed as
\begin{align*}
    (G^{[3]}_{\fn,\fz})^*G^{[3]}_{\fn',\fz'}(y',z')= 2^{-l}|S\Pi_{\cN_2}e_1|^{-1}\int e^{i2^k\Phi_L(u,s,y',z')} \mathring\cA^L_{N,k,l}(u,s,y',z') duds
\end{align*}
where
\begin{align*}
    \Phi_L(u,s,y',z') &= \Psi(L(u,s),z') - \Psi(L(u,s),y'),\\
    \mathring\cA^L_{N,k,l}(u,s,y',z') &= \cA_{N,k,l}(L(u,s),y',z') \chi_{\fq(\fz)\cap\fr^{[3]}(\fz,\fz';\fn)}(y')\chi_{\fq(\fz')\cap\fr^{[3]}(\fz',\fz;\fn')}(z').
\end{align*}
Note that
\begin{align*}
    &|\inn{\bar\fe_1}{\partial_{u,s}}^N\mathring\cA^L_{N,k,l}(u,s,y',z')|\\
    &\qquad\lesssim (2^{-l}|S\Pi_{\cN_2}e_1|^{-1})^N\big|\inn{\frac{S\Pi_{\cN_2}e_1}{|S\Pi_{\cN_2}e_1|}}{\partial_{x,t}}^N\cA_{N,k,l}(x,t,y',z')\big|_{(x,t) = L(u,s)}\big|
\end{align*}
for $N\in\N$.
Writing $S\Pi_{\cN_2}e_1 = Se_1 - S\Pi_{(\cN_2)^\perp}e_1$ and using the estimate above together with \eqref{e:Se1Picomp}, $|S\Pi_{(\cN_2)^\perp} e_1|\le 10^{-2}2^{-l}$, and $ii)$ in Lemma \ref{lem:decayfzfz'}, we see that the derivatives of the amplitude function $\mathring\cA^L_{N,k,l}$ in $u,s$ satisfy the bounds
\begin{align*}
    |\partial_{u,s}^\alpha \mathring\cA^L_{N,k,l}(u,s,y',z')|\lesssim_{N,\alpha} (1+|\fz-\fz'|)^{-N}
\end{align*}
for $N\in \N$ and $\alpha \in \N_0^{d+1}$. Moreover, by \eqref{e:Se1Picomp}, the size of the support of $(u,s)\mapsto \mathring\cA^L_{N,k,l}(u,s,y',z')$ is uniformly bounded in $y',z'$, hence
\begin{align}\label{e:intringA}
    \int |\partial_{u,s}^\alpha\mathring\cA^L_{N,k,l}(u,s,y',z')| duds \lesssim_{N,\alpha} (1+|\fz-\fz'|)^{-N}.
\end{align}
On the other hand,
\begin{align}
\begin{aligned}\label{i:pusPsiL}
    \partial_{u,s}\partial_{y'}^\intercal &\Psi(L(u,s),w^\tau)(y'-z') \\
    & = LS\Pi_{\mathcal N_2}P^\intercal(y'-z') + \cR_3(u,s,y',z',\tau) + \cQ_3(u,s,y',z',\tau)
\end{aligned}
\end{align}
where
\begin{align*}
    \cR_3(u,s,y',z',\tau) &= LSe_1\inn{g'(\frac{x'-w^\tau}{t})}{y'-z'},\\
    \cQ_3(u,s,y',z',\tau) & = t^{-1}\sigma(L(u,s)) L\fg(L(u,s)) + LS\Pi_{(\mathcal N_2)^\perp}P^\intercal(y'-z').
\end{align*}
We use that $L$ is a symmetric matrix, that is, $L^\intercal = L$. We suppress the dependence of $\fg$ on $y',z',\tau$ for simplicity. As before, we aim to obtain a lower bound for a directional derivative of $\Phi_L$ and use repeated integration by parts, which will yield \eqref{e:G*G3}. In light of this, we estimate each summand on the right-hand side of \eqref{i:pusPsiL}. Since $(y',z')\in \fr^{[3]}(\fz,\fz';\fn)\times \fr^{[3]}(\fz',\fz;\fn')$ and $|\fn-\fn'| > 10\sqrt d$, we see that there exists a constant $C_\circ^{[3]}>0$ such that
\begin{align}\label{e:lb-lead3}
    |LS\Pi_{\mathcal N_2}P^\intercal(y'-z')|\ge C_\circ^{[3]}K2^{-k}(1+|\fz-\fz'|)|\fn-\fn'|.
\end{align}
We proceed to estimate $\cQ_3$. By definition, we have
\[
|\sigma(L(u,s))|\lesssim 2^{-l},\quad |S\Pi_{(\mathcal N_2)^\perp}P^\intercal(y'-z')|\lesssim 2^{-l}|y'-z'|,
\]
Combining these bounds with $|y'-z'|\lesssim 2^{-k+l}(1+|\fz-\fz'|)$ gives
\begin{align}\label{e:ub-cQ3}
    |\cQ_3(u,s,y',z',\tau)|\le C_g^{[3]} 2^{-k}(1+|\fz-\fz'|).
\end{align}
Lastly, to estimate $\cR_3$, we note that
\begin{align*}
    |\inn{g'(\frac{x'-w^\tau}{t})}{y'-z'}|\lesssim \eps |y'-z'|\lesssim \eps 2^{-k+l}(1+|\fz-\fz'|)
\end{align*}
and by \eqref{e:Se1Picomp}, $|LSe_1| \lesssim |S\Pi_{\cN_2}e_1||L\bar\fe_1| + 2^{-l} \lesssim 2^{-l}$. Putting together these bounds gives
\begin{align}\label{e:ub-cR3}
    |\cR_3(u,s,y',z',\tau)|\le C_*^{[3]} \eps  2^{-k}(1+|\fz-\fz'|)
\end{align}
with a constant $C_*^{[3]}>0$. Thus, letting
\[
v^{[3]} = \frac{LS\Pi_{\mathcal N_2}P^\intercal(y'-z')}{|LS\Pi_{\mathcal N_2}P^\intercal(y'-z')|}
\]
and using \eqref{e:lb-lead3}--\eqref{e:ub-cR3} together with \eqref{i:diffPsi}, we have
\begin{align}\label{e:lbPhiL}
    |\inn{v^{[3]}}{\partial_{u,s}\Phi_L(u,s,y',z')}|\ge CK2^{-k}(1+|\fz-\fz'|)|\fn-\fn'|
\end{align}
with a constant $C>0$ provided that $K$ is so large that 
\[
C_\circ^{[3]}K\ge \max(C_g^{[3]}, C_*^{[3]} \eps).
\]
Furthermore, by a calculation,
\begin{align*}
    |\inn{v^{[3]}}{\partial_{u,s}} \sigma(L(u,s))| = |\inn{v^{[3]}}{LSe_1}|\le|LSe_1|\lesssim 2^{-l}.
\end{align*}
Consequently, for $N\ge 2$,
\begin{align}\label{e:ubPhiL}
    |\inn{v^{[3]}}{\partial_{u,s}}^N\Phi_L(u,s,y',z')|\lesssim_N 2^{-k}(1+|\fz-\fz'|).
\end{align}
Repeated integration by parts in the $v^{[3]}$-direction, using \eqref{e:intringA}, \eqref{e:lbPhiL}, and \eqref{e:ubPhiL}, yields \eqref{e:G*G3} and completes the proof of Proposition \ref{prop:estGnn'3}.

\section{Necessary conditions and examples}\label{sec:necessary}

\subsection{Proof of Theorem \ref{thm:main-hne}}
Recall that $\Sigma$ is the boundary of a convex open domain in $\R^d$ containing $0$, with nonvanishing Gaussian curvature everywhere. Assume that
\begin{align}\label{e:assum-loja}
    \cL\equiv\cL(\cC_r,d_\perp\,|\,\cB) > 2r.
\end{align}
Recall that $m^\sharp$ denotes the rank of $\La^\intercal$. As before, after making changes of variables in $\ox$ and in $\oy$, we may assume \eqref{eq:la-split}. In this setting, we have
\[
d_\perp(\theta) = \big(\sum_{i=m-m^\sharp+1}^{m} \theta_i^2\big)^{\frac12}.
\]
By \eqref{e:assum-loja} and \eqref{e:prod-cd-r}, the product of the $r$ largest singular values $s_i(\theta)$, $1\le i\le d$, satisfies
\begin{align}\label{i:sing-L}
    \inf A_\cB = \frac{\cL}{2} >r
\end{align}
with $A_\cB :=\{\beta: \exists\,c>0 \ \text{s.t.}\ \forall\theta\in \cB\,\big(\prod_{i=1}^r s_i(\theta) \ge cd_\perp(\theta)^\beta\big)\}$.
Let
\begin{align*}
    A_{\bbS} = \big\{\beta: \exists\,c>0 \ \text{s.t.}\ \forall\theta\in \bbS^{m-1}\,\big(\prod_{i=1}^r s_i(\theta) \ge cd_\perp(\theta)^\beta\big)\big\}.
\end{align*}
Then, note that
\begin{align}\label{i:Aeps-L}
    \inf A_{\bbS} = \frac{\cL}{2}.
\end{align}
Indeed, it follows immediately from the definition that $A_\cB \subset A_\mathbb S$. On the other hand, since $\cC_r$ is homogeneous of degree $2r$, by \eqref{e:prod-cd-r}, there exist $c_\beta, c_{\beta,\cB}'>0$ for $\beta\in A_{\bbS}\cap [r,\infty)$ such that
\begin{align*}
    |\cC_r(\theta)| = |\theta|^{2r}|\cC_r(\theta/|\theta|)| \ge c_\beta |\theta|^{2r} d_\perp(\theta/|\theta|)^{2\beta} \ge c_{\beta,\cB}'d_\perp(\theta)^{2\beta}
\end{align*}
for all $\theta\in \cB$. Moreover, $(-\infty,r)\cap A_\bbS = \emptyset$ because by homogeneity of $\cC_r$, $r$ belongs to $A_\cB$ otherwise. Using \eqref{e:prod-cd-r} again, we deduce $A_{\bbS} = A_\cB$. Thus, $\inf A_{\bbS} = \cL/2$. %Now choose $0<\eps<1$. It is immediate from the definition that $A_1\subset A_\eps$. On the other hand, let $\beta\in A_\eps$. Then $\prod_{i=1}^r s_i(\theta)\ge c_\beta d_\perp(\theta)^\beta$ for $\theta\in \bbS^{m-1}$ satisfying $d_\perp(\theta)\in (0,\eps]$. By our assumption, $\prod_{i=1}^r s_i(\theta)\neq 0$ whenever $d_\perp(\theta)>0$, thus there exists $C_\eps>0$ such that $\prod_{i=1}^r s_i(\theta)\ge C_\eps$ for $\theta\in \bbS^{m-1}$ such that $d_\perp(\theta) \ge \eps$. Hence,
%\begin{align*}
%    \prod_{i=1}^r s_i(\theta)\ge \max(c_\beta, C_\eps) d_\perp(\theta)^\beta
%\end{align*}
%for all $\theta\in \bbS^{m-1}$, which implies $A_\eps = A_1$. Thus, $\inf A_\eps = \inf A_1$. Combining this with $\inf A_1 = \cL$ gives \eqref{i:Aeps-L}.

By \eqref{i:Aeps-L}, we can choose $\alpha\in (r, \cL/2)$ so that the following holds: there exist a sequence $\{\eps_k\}_{k\in \N}$ of positive numbers converging to $0$ and a sequence $\{\vartheta^k\}_{k\in \N}$ of unit vectors in $\R^m$ such that $d_\perp(\vartheta^k) = \eps_k$ and
\begin{align}\label{e:prodvarth}
    \prod_{i=1}^r s_i(\vartheta^k) \le d_\perp(\vartheta^k)^\alpha = \eps_k^\alpha.
\end{align}
For $k\in \N$, let $\wt A_t^k$ be an operator on $\R^{d+1}$ given by
\begin{align*}
     \wt A_t^k f(\ubar x,x_{d+1})= \int_\Sigma f(\ubar x-t\om, x_{d+1} -t(x,t)^\intercal S^{\vartheta^k}\om) d\mu(\om).
\end{align*}
Throughout the proof, we assume that $I=[1/4,4]$. We define $\wt\fM^k f(x) = \sup_{t\in I} |\wt A_t^k f(x)|$ for $k\in \N$. Then, we claim that
\begin{align}\label{e:inf-wtfMk}
    \|\fM\|_{L^p\to L^p} \gtrsim \sup_{k} \|\bbone_{\wt B}\wt\fM^k\|_{L^p\to L^p},
\end{align}
where $\wt B$ denotes the unit ball in $\R^{d+1}$.

The proof of \eqref{e:inf-wtfMk} follows the same argument as in \cite{RyuSeeger}. Choose $k\in \N$ and fix it. There exist unit vectors $\{b^{k,i}\}_{2\le i\le m}\subset \R^m$ such that $\{\vartheta^k,b^{k,2},\dots,b^{k,m}\}$ forms an orthonormal basis of $\R^m$. Let $R_k$ be an $m\times m$ rotation matrix with columns $\vartheta^k,b^{k,2}\dots b^{k,m}$. A calculation gives
\begin{align*}
    A_t f(\ux, \ox) = \int f\big(\ux-t\om, R_k \big(R_k ^\intercal \ox - t R_k ^\intercal \big((\ux,t)^\intercal \vec S \om\big)\big)\big) d\mu(\om),
\end{align*}
where
\begin{align*}
    \inn{R_k ^\intercal ((\ux,t)^\intercal \vec S \om)}{e_1} = (\ux,t)^\intercal (\sum_{i=1}^m \inn{\vartheta^k}{e_i} S_i)\om.
\end{align*}
By \eqref{eq:la-split} and the definition of $\vartheta^k$, we have 
\[
C_{\La,1}\eps_k\le |\sum_{i=1}^m \inn{\vartheta^k}{e_i} \la_i|\le C_{\La,2}\eps_k
\]
with constants $0<C_{\La,1}<C_{\La,2}$ depending only on $\La$. Thus, after conjugating with $R_k$, we may identify $S^{\vartheta^k}$ with $S_1$. Moreover, we may assume that
\begin{align}\label{cond:S1la1}
    \rank(S_1) \ge r,\quad C_{\La,1}\eps_k\le |\la_1|\le C_{\La,2}\eps_k,\quad \prod_{i=1}^r s_i \le \eps_k^\alpha
\end{align}
where $\{s_i\}$ denote singular values of $S_1$. We now consider a function $f$ given by $f(\uy, \oy_1,\oy') = \wt f(\uy, \oy_1) h(\oy')$ with $\wt f\in C_c^\infty(\R^{d+1})$, $h\in C_c^\infty(\R^{m-1})$. Choose $h\in C^\infty_c(\bbR^{m-1})$ so that $h\equiv 1$ on a large compact set. Note that we may choose $h$ independently of $k$. %specifically
%\[h(\sum_{i=2}^m (\ox_i - t (\ux, t)^\intercal S_i \om) e_{d+i}) = 1\] for every $\om \in \supp(\mu)$, $t\in I$, $x\in \widetilde B\times B'$.
Then 
\begin{align*}
    \big\|\bbone_{\widetilde B\times B'}\sup_{t\in I}\big|A_t f\big|\big\|_{L^p(\bbR^{d+m})} = |B'|^{1/p} \big\|\bbone _{\wt B}\sup_{t\in I}\big|\wt A^k_t \wt f\big|\big\|_{L^p(\bbR^{d+1})},
\end{align*}
where $B'$ denotes the unit ball in $\R^{m-1}$. Recall that we have identified $S^{\vartheta^k}$ with $S_1$. Hence
\begin{align*}
    \|\fM\|_{L^p\to L^p} \ge C \|\bbone_{\wt B} \wt \fM^k\|_{L^p \to L^p}
\end{align*}
with a constant $C>0$ independent of $k$. Taking the supremum over $k\in \N$ proves the claim.

To conclude the proof of Theorem \ref{thm:main-hne}, we will prove that there exists $a>0$ and $p_\circ \in (\frac{r+1}{r},\infty]$ such that
\begin{align}\label{e:lb-wtfMk}
    \|\bbone_{\wt B}\wt\fM^k\|_{L^p\to L^p}\gtrsim \eps_k^{a(\frac{1}{p_\circ}-\frac1p)}.
\end{align}
By \eqref{e:inf-wtfMk}, it is immediate that \eqref{e:lb-wtfMk} implies Theorem \ref{thm:main-hne}.

To establish \eqref{e:lb-wtfMk}, let $\delta_k = s_r(\vartheta^k)$ for $k\in \N$. Fix $k\in \N$ and simply write $s_i = s_i(\vartheta^k)$ for $1\le i\le d$. As before, we use the notation $S_1$ instead of $S^{\vartheta^k}$. We write $S_1 = \begin{pmatrix}
    J_1 \\ \la_1^\intercal
\end{pmatrix}$. Choose an orthonormal basis $\{u_1,\dots, u_d\}$ of $\R^d$ such that $S_1^\intercal S_1 u_i = s_i^2 u_i$ for all $1\le i\le d$. Recall that $\Sigma = \partial\Omega$, where $\Omega$ is a bounded convex domain containing the origin. Choose $\om^\circ\in \Sigma$ on the ray generated by $u_d$. Consider the tangent hyperplane $\om^\circ + T_{\om^\circ}\Sigma$. Since $0\in \Omega$, $0\notin \om^\circ + T_{\om^\circ}\Sigma$. Thus, $\Sigma$ can be parametrized near $\om^\circ$ by
\[
\Gamma(\om) = \sum_{i=1}^{d-1} \om_i u_i + G(\om) u_d,\quad \om\in B(0,\delta_*)\subset \R^{d-1}
\]
with a small $\delta_*>0$ and $G\in C^\infty(\R^{d-1})$ such that $G(0) u_d = \om^\circ$. Let $\Sigma_*$ denote the image of $B(0,\delta_*)$ under $\Gamma$. Since the test function $f_*$ is nonnegative, restricting the integral to $\Sigma_*$ gives a lower bound for $\widetilde A_t^k f_*$.

%For simplicity, let $\Pi = \Pi_{\vartheta^k}$ throughout the proof. 
We set
\begin{align*}
    U_* &= \{x = (\ux,x_{d+1})\in \R^{d+1}: |\ux|\le c_*\delta_*,\ \la_1^\intercal \ux,\,x_{d+1}\in [a_*\eps_k/2, 2a_*\eps_k]\}, \\
    V_* &= \{y = (\uy,y_{d+1}) \in\R^{d+1}: |\inn{\uy}{u_i}|\le |\inn{u_i}{\fu}|\ \text{for}\ 1\le i\le d,\ |y_{d+1}|\le C_*\delta_k\},
\end{align*}
where $c_*, C_*$ are sufficiently small and large constants, respectively, and
\[
a_* = \frac{C_{\La,1}c_*\delta_*}{4},\quad \fu = c_*\delta_*\big(\sum_{i=1}^r \delta_k s_i^{-1}u_i + \sum_{i=r+1}^{d-1} u_i\big) + C_*|\om^\circ|u_d.
\]
By \eqref{cond:S1la1}, we have
\begin{align}\label{e:sizeUV*}
    |U_*| \sim (c_*\delta_*)^d\eps_k,\quad |V_*|\sim \delta_k^{r+1} \big(\prod_{i=1}^r s_i\big)^{-1}.
\end{align}
Lastly, let $f_*(y) = \chi_{V_*}(y)$. Note that
\begin{align*}
    \wt A_t^k f_*(\ux,x_{d+1}) \ge \int f_*(\ux-t\om, x_{d+1} - t\la_1^\intercal\ux + (\ux,t)^\intercal S_1 (\ux - t\om)) \chi_{\Sigma_*}(\om)d\mu(\om).
\end{align*}
For $(\ux,x_{d+1})\in U_*$, we choose $t_x = x_{d+1}/(\la_1^\intercal \ux)$. By the definition of $U_*$, $t_x\in [1/4, 4]$. By the parametrization of $\Sigma_*$ by $\Gamma$, the surface measure of
\begin{align*}
    \cW(\ux) = \big\{\om\in\Sigma_*: |\inn{\ux - t_x\om}{u_i}|\le |\inn{u_i}{\fu}|\ \forall 1\le i\le d\big\}
\end{align*}
satisfies
\begin{align*}
    \big|\cW(\ux)\big| \gtrsim \delta_k^r \big(\prod_{i=1}^r s_i\big)^{-1} \quad \forall x=(\ux,x_{d+1})\in U_*.
\end{align*}
Moreover, since $s_i(\vartheta^k)\le \delta_k$ for $r+1\le i\le d$, we have that, for every $x\in U_*$,
\[
|x_{d+1} - t_x\la_1^\intercal\ux + (\ux,t_x)^\intercal S_1 (\ux - t_x\om)| \le C_*\delta_k
\]
for all $\om \in \cW(\ux)$, if $C_*$ is sufficiently large. Thus, for $x\in U_*$,
\begin{align*}
    |\wt \fM^k f_*(x)| \ge |\wt A_{t_x}^k f_*(x)| \gtrsim |\cW(\ux)| \gtrsim \delta_k^r \big(\prod_{i=1}^r s_i\big)^{-1}.
\end{align*}
Combining this with \eqref{e:sizeUV*}, we obtain
\begin{align}\label{e:wtfM-f*Lp}
    \|\wt\fM^k f_*\|_{L^p(U_*)} \gtrsim \delta_k^{r-\frac{r+1}{p}} \big(\prod_{i=1}^r s_i\big)^{-1 + \frac 1p} \eps_k^{\frac 1p} \|f_*\|_{L^p}
\end{align}
where the implicit constant depends on $c_*, \delta_*, C_*$. We suppress the dependence on these constants because these are independent of $k$. Note that
\begin{align*}
    \delta_k \ge (\sup_k\|S^{\vartheta^k}\|^{r-1})^{-1}\prod_{i=1}^r s_i \ge c(\sup_k\|S^{\vartheta^k}\|^{r-1})^{-1}\eps_k^{\frac{\cL}{2}+\varrho}
\end{align*}
for all sufficiently small numbers $\varrho>0$ with $c>0$ by the definition of $\cL$. Note that $\|S^{\vartheta^k}\|$ is uniformly bounded in $k$, thus $(\sup_k\|S^{\vartheta^k}\|^{r-1})^{-1}$ above is acceptable. Suppose now that $p\in ((r+1)/r,\infty]$. Substituting the inequality above and $\prod s_i \le \eps_k^\alpha$ into \eqref{e:wtfM-f*Lp}, we get
\begin{align*}
    \|\wt\fM^k f_*\|_{L^p(U_*)} \gtrsim \eps_k^{\frac{\cL r-2\alpha}{2}-\frac{\cL(r+1)-2\alpha-2}{2p}+\varrho'}\|f_*\|_{L^p}
\end{align*}
for every small number $\varrho'>0$. Indeed, we can choose $\varrho'=r\varrho$.
Since $r\ge 2$ by assumption, both $\cL r-2\alpha$ and $\cL(r+1)-2\alpha-2$ are positive. Moreover, setting
\[
\frac{1}{p_\circ}:=\frac{\cL r-2\alpha+2\varrho'}{\cL(r+1)-2\alpha-2},
\]
we have $\frac{1}{p_\circ} < \frac{r}{r+1}$ provided $\varrho'$ is sufficiently small. This completes the proof of \eqref{e:lb-wtfMk}.

\subsection{Examples with large Łojasiewicz exponents} In this subsection, we construct $(d+1)\times d$ matrices $\{S_i\}_{1\le i\le m}$ for some tuples $(d,m,r)$ such that $r=\min\{r(\theta):\theta\notin \ker(\La^\intercal)\}$ and $\cL(\cC_r,d_\perp\,|\,\cB)$ is larger than $2r$. For $(d,m,r)\in \N^3$, let $n = \big[\frac r2\big]\in \Z$ and
\begin{align*}
    \fI(d,m,r) = \big\{\wt D\in \N: \wt D\ge 2,\ N^*(m-1,2\wt D)< n\big\}
\end{align*}
where $N^*$ denotes a function on $\N^2$ defined by
\[
N^*(k,l) = \begin{cases} 2\binom{k'+l-1}{k'}-1 & k=2k'\\
2\sqrt{3l+1}\binom{k'+l-1}{k'}-1 & k=2k'+1.
\end{cases}
\]
For $(d,m,r)\in \N^3$ such that $\fI(d,m,r)\neq \emptyset$, let $D\equiv D(d,m,r)$ denote the maximum of $\fI(d,m,r)$. Throughout this subsection, we assume $m\ge 3$ because the {\L}ojasiewicz exponent is always equal to $2r$ otherwise, as noted in $(a)$ of \textit{Remarks} in Section \ref{sec:intro}.

\begin{prop}\label{prop:matrixconst}
    Let $d\ge 2$, $m\ge 3$, and $2\le r\le d$. Assume that $\fI(d,m,r)\neq \emptyset$. Then, there exists a collection $\{S_i\}_{i=1}^m$ of $(d+1)\times d$ matrices such that $\min\{r(\theta):\theta\notin \ker(\La^\intercal)\} = r$ and 
    \begin{align}\label{e:lb-lojarD}
        \cL(\cC_r,d_\perp\,|\,\cB) \ge 2r+8D^{m-2}(D-1). %2^{2-m}n^2
    \end{align}
\end{prop}

To prove the proposition, we use the known result that every multivariate polynomial with real coefficients can be represented as a determinant of an affine linear matrix pencil. We first introduce notation for linear matrix pencils. Let $k,N\in \N$. An $N\times N$ linear pencil $L_M$ in the $k$ indeterminates $x_1,\dots, x_k$ is an expression of the form $L_M(x) = x_1M_1 + \dots + x_kM_k$ where $M_1,\dots,M_k$ are $N\times N$ matrices with real entries. Likewise, an $N\times N$ affine linear pencil $\wt L_M$ denotes an expression of the form $\wt L_M = M_0 +x_1M_1+\dots+x_kM_k$ where $M_0,\dots,M_k$ are $N\times N$ real matrices.

\begin{thm}[\cite{BDDHP17}]\label{thm:determinantal}
    Let $k\ge 2$ and $l\ge 1$. Let $p(x)$ be a polynomial of degree $l$ in $k$ variables over $\R$ such that $p(0)\neq 0$. Then, there is an $N\times N$ linear pencil $L_M$ such that
    \[
    p(x) = p(0) \det(I + L_M(x))
    \]
    and $N$ is the greatest integer less than or equal to $N^*(k,l)$.
\end{thm}

This follows from the construction in \cite[p. 426]{BDDHP17}. See also \cite{HMV06, Qu12} for related determinantal representations.

\begin{proof}[Proof of Proposition \ref{prop:matrixconst}]
    Throughout the proof, write $\theta = (\theta_1, \theta')\in \R\times\R^{m-1}$. Note that $D\ge 2$. Let $p_D$ be a polynomial of degree $2D$ in $m-1$ variables defined by
    \begin{align}\label{i:def-pD}
        p_D(\theta') = (\theta_2^{D-1}\theta_3-1)^2+(\theta_4-\theta_3^D)^2+\dots+(\theta_m-\theta_{m-1}^D)^2+\theta_m^{2D}.
    \end{align}
    We claim that there are $(d+1)\times d$ matrices $\{S_k\}_{k=1}^m$ such that the following holds:
    \begin{itemize}[topsep=2pt, itemsep=2pt]
        \item [a)] $S_k = \begin{pmatrix}
            J_k \\ \la_k^\intercal
        \end{pmatrix}$ where $J_k$ is a $d\times d$ skew-symmetric matrix. In particular, $|\la_1| = 1$ and $\la_k = 0$ for $2\le k\le m$.
        \item [b)] $\min\{r(\theta):\theta_1\neq 0\} = r$.
        \item [c)] The polynomial $\cC_r(1,\theta')$ is equal to $a(p_D(\theta'))^4$ with $a\neq0$.
    \end{itemize}
    We first show that $a)$--$c)$ imply \eqref{e:lb-lojarD}, and then construct the matrices. It suffices to show \eqref{e:lb-lojarD} with $\{S_k\}_{k=1}^m$ satisfying $a)$--$c)$ above. By $a)$, $d_\perp(\theta) = |\theta_1|$. Since $\cC_r$ is a homogeneous polynomial of degree $2r$, by $c)$,
    \[
    |\cC_r(\theta_1,\theta')| = |a|\theta_1^{2r}(p_D(\theta'/\theta_1))^4.
    \]
    Now consider a curve $\gamma_D:(0,1]\to \R^m$ given by 
    \[
    \gamma_D(t) = (t, 1, t^D,t^{D(D-1)+1},\dots,t^{D^{m-3}(D-1)+1}).
    \]
    Clearly, $\gamma_D(t)\in \cB$ for $t\in (0,1]$ and
    \begin{align*}
        |\cC_r(\gamma_D(t))| = |a| t^{2r+8D^{m-2}(D-1)}.
    \end{align*}
    Thus, we have that
    \[
    2r+8D^{m-2}(D-1)-\eps\notin \big\{\beta:\exists c>0\ \text{s.t.}\ \forall \theta\in \cB\,\big(|\cC_r(\theta)|\ge c|\theta_1|^\beta\big)\big\}
    \]
    for every $\eps>0$. Therefore, \eqref{e:lb-lojarD} holds.
    
    We now prove the claim. Let $N$ be the largest integer at most $N^*(m-1,2D)$. By the definition of $D$, $N< n$. By Theorem \ref{thm:determinantal}, there are $N\times N$ matrices $\{M_k\}_{k=2}^m$ such that $p_D(\theta') = \det(I_N+L_M(\theta'))$ where $L_M(\theta')$ denotes the $N\times N$ linear pencil $\theta_2M_2+\dots+\theta_m M_m$. Let $\{\wt M_k\}$ be $n\times n$ matrices defined by
    \begin{align*}
        \wt M_1 = I_n,\quad \wt M_k =  \begin{pmatrix}
            0 & 0 \\
            0 & M_k
        \end{pmatrix}\quad \text{for}\ 2\le k\le m.
    \end{align*}
    By definition, the determinant of the linear pencil $L_{\wt M}(\theta) = \theta_1\wt M_1 +\dots +\theta_m \wt M_m$ is a homogeneous polynomial $g_{\wt M}(\theta)$ of degree $n$ such that $g_{\wt M}(1,\theta') = p_D(\theta')$. We have to extend $\{\wt M_k\}$ to $[d/2]\times [d/2]$ matrices so as to construct $S_k$. We define $[d/2]\times [d/2]$ matrices $\{A_k\}_{k=1}^m$ by adjoining $[d/2]-n$ zero rows above $\wt M_k$ and $[d/2]-n$ zero columns to its left. Moreover, define $\la\in \R^d$ by
    \begin{align*}
        \la = \begin{cases}
            e_{d-[\frac d2]-n+1} & r=2n, \\
            e_1 & r=2n+1.
        \end{cases}
    \end{align*}
    Now let $\{S_k\}_{k=1}^m$ be $(d+1)\times d$ matrices defined by
    \begin{align*}
        S_1 = \begin{pmatrix}
            0 &  -A_1^\intercal \\
            A_1 &  0 \\
            \multicolumn{2}{c}{\la^\intercal}
        \end{pmatrix},\quad S_k = \begin{pmatrix}
            0 &  -A_k^\intercal \\
            A_k&  0 \\
            \multicolumn{2}{c}{0}
        \end{pmatrix}
    \end{align*}
    for $2\le k\le m$ when $d$ is even. For odd $d$, we define
    \begin{align*}
        S_1 = \begin{pmatrix}
            0 & 0 & 0\\
            0 & 0 &-A_1^\intercal \\
            0 & A_1 &  0 \\
            \multicolumn{3}{c}{\la^\intercal}
        \end{pmatrix},\quad S_k = \begin{pmatrix}
            0 & 0 & 0\\
            0 & 0 &  -A_k^\intercal \\
            0 & A_k &  0 \\
            \multicolumn{3}{c}{0}
        \end{pmatrix},\quad 2\le k\le m.
    \end{align*}
    %first define $(2[d/2]+1)\times 2[d/2]$ matrices $\wt S_k$ in the same manner as above. Then, we define $S_k$ by inserting one additional zero rows and columns to the top and left of $\wt S_k$, respectively. 
    It is immediate that $a)$ holds for the matrices $\{S_k\}$. To verify that $b)$ and $c)$ also hold, we compute $\cC_r(\theta)$ for $S^\theta$. By direct computation, $(S^\theta)^\intercal S^\theta$ can be written in block form as
    \begin{align*}
        (S^\theta)^\intercal S^\theta = \begin{pmatrix}
        0 & 0 & 0 & 0 \\
        0 & (\wt M^\theta)^\intercal \wt M^\theta & 0 & 0 \\
        0 & 0 & 0 & 0 \\
        0 & 0 & 0 & \wt M^\theta(\wt M^\theta)^\intercal
    \end{pmatrix} + \theta_1^2\la\la^\intercal
    \end{align*}
    where $\wt M^\theta = \theta_1\wt M_1 + \dots +\theta_m \wt M_m$. From this, it follows that 
    \[
    \cC_r(\theta) = (-1)^{d-r} \det(\mathbf M_\theta)
    \]
    where $\mathbf M_\theta$ is an $r\times r$ matrix given by
    \begin{align*}
        \mathbf M_\theta = \begin{pmatrix}
            \theta_1^2 & 0 & 0 \\
            0 & (\wt M^\theta)^\intercal \wt M^\theta & 0 \\
            0 & 0 & \wt M^\theta (\wt M^\theta)^\intercal
        \end{pmatrix}
    \end{align*}
    if $r$ is odd, and
    \begin{align*}
        \mathbf M_\theta = \begin{pmatrix}
            (\wt M^\theta)^\intercal \wt M^\theta & 0 \\
            0 & \wt M^\theta (\wt M^\theta)^\intercal
        \end{pmatrix} + \theta_1^2e_1e_1^\intercal
    \end{align*}
    if $r$ is even. Here, $e_1$ denotes $(1,0,\dots,0)\in \R^r$. We first consider the case where $r$ is odd. Note that $\det(\mathbf M_\theta) = \theta_1^2\det(\wt M^\theta)^4 = \theta_1^2 (g_{\wt M}(\theta))^4$. Clearly, the condition $c)$ holds if we set $a=(-1)^{d-r}$. Furthermore, since $g_{\wt M}(\theta) = \theta_1^n p_D(\theta'/\theta_1)$ and $p_D$ is a positive polynomial, $\det(\mathbf M_\theta)$ does not vanish as long as $\theta_1\neq 0$. This proves $b)$. Thus, $\{S_k\}$ satisfy all the desired conditions $a)$--$c)$.

    Lastly, assume that $r$ is even. Recall that $N < n$. Set $M^\theta=\theta_1I_N+\sum_{k=2}^m\theta_kM_k$, so that $\widetilde M^\theta=\operatorname{diag}(\theta_1I_{n-N},M^\theta)$. Note that $\mathbf M_\theta$ is expressed as
    \begin{align*}
        \mathbf M_\theta = \begin{pmatrix}
            \theta_1^2I_{n-N} + \theta_1^2e_1e_1^\intercal & 0 & 0 \\
            0 &(M^\theta)^\intercal M^\theta & 0 \\
            0 & 0 & \wt M^\theta(\wt M^\theta)^\intercal
        \end{pmatrix}
    \end{align*}
    where $e_1$ denotes $(1,0,\dots,0)\in \R^{n-N}$. Thus, a straightforward calculation shows
    \[
    \det(\mathbf M_\theta) = 2(g_{\wt M}(\theta))^4.
    \]
    Since $g_{\wt M}(\theta)$ does not vanish if $\theta_1\neq 0$, we see that the conditions $b)$ and $c)$ hold with $a=2\cdot (-1)^{d-r}$. This concludes the construction of $\{S_k\}$.
\end{proof}

\appendix

\section{Reduction to Theorem \ref{thm:main2}}
\label{app:basic-reductions}
 
%discuss how this result can be applied to deduce Theorem \ref{thm:main-m=1} and \ref{thm:main-hsu}.
In this appendix, we explain how Theorems \ref{thm:main-m=1} and \ref{thm:main-hsu} follow from Theorem \ref{thm:main2}. We only treat $m\ge2$. For $m=1$, the linear independence reduction is unnecessary and the remaining argument is the same. The reduction to Theorem \ref{thm:main-hsu} involves a localization argument and changes of variables in $x$ and $y$,  as in \cite{RyuSeeger} for the case $\La=0$. Here we need to make sure that both $r$ and $\cL(\cC_r,d_\perp\,|\,\cB)$ are not affected by this reduction.

    We first discuss why we may assume that  $S_1,\dots, S_m$ are linearly independent. Let $(E_k)_{k=1}^{\frac{d(d-1)}{2}}$ be a basis of the vector space of $d\times d$ skew-symmetric matrices. We then define a collection of $(d+1)\times d$ matrices $(\wt E_j)_{j=1}^{\frac{d(d+1)}{2}}$ by
\begin{align*}
    \wt E_k := \begin{amatrix}{1|1}
        E_k & 0
    \end{amatrix}^\intercal, \,1\le k\le\tfrac{d(d-1)}2,\quad  \wt E_{\nu+\frac{d(d-1)}{2}} := \begin{amatrix}{1|1}
        0 & e_\nu
    \end{amatrix}^\intercal, \, 1\le \nu\le d.
\end{align*}
% for $1\le j\le d(d-1)/2$, and
% \begin{align*}
%     \wt E_{j+\frac{d(d-1)}{2}} := \begin{amatrix}{1|1}
%         0 & e_j
%     \end{amatrix}^\intercal
% \end{align*}
% for $1\le j\le d$. 
Each $S_i$ is expressed as a linear combination of $\wt E_j$'s, that is, $S_i = \sum_{j=1}^{\frac{d(d+1)}{2}} c_{ij} \wt E_j$. Denote $\cC = (c_{ij})_{1\le i\le m,\ 1\le j\le \frac{d(d+1)}{2}}$. Applying a singular value decomposition to $\cC$ we have $\cC = U D V$, where $U\in SO(m)$ and  $V\in SO(\frac{d(d+1)}{2})$. $D=(d_{ij})$ is an $m\times \frac{d(d+1)}{2}$ diagonal matrix such that its diagonal entries $d_{ii}$  are given by $d_{ii}=s_i\ge 0$ where $s_i > 0$ for $1\le i\le n$,  with $n\le \min(m,\frac{d(d+1)}{2})$, while $s_i = 0$ for $n < i\le \max(m,\frac{d(d+1)}{2})$. 
For all other entries we have $d_{ij}=0$. 

For $1\le i\le m$, let
\begin{align*}
    L_i := s_i\sum_{j=1}^{\frac{d(d+1)}{2}} V_{ij} \wt E_j\ \  \text{for $1\le i\le n$},\quad L_i = 0\ \ \text{otherwise}.
\end{align*}
Then, note that $\ux^\intercal \vec J \om + t\La \om = U\big((\ux, t)^\intercal \vec L \om\big)$ where $\vec L = (L_1,\dots, L_m)$. Thus
\begin{align*}
    f\ast \mu_t(x) = \int_\Sigma f(\ux - t\om, U(U^\intercal \ox - t (\ux, t)^\intercal \vec L\,\om))d\mu(\om),
\end{align*}
which implies 
%\begin{align*}
   $ \fM f(\ux, \ox) = \fM^{[1]} f_U(\ux, U^\intercal \ox)$ where
%\end{align*}
%Here, we denote
\[
\fM^{[1]} f(x) = \sup_{t > 0} \Big|\int_\Sigma f(\ux - t\om, \ox - t(\ux, t)^\intercal \vec L \om)d\mu(\om)\Big|,\ \ f_U(y) = f(\uy, U\oy). 
\]
The $L^p\to L^p$ operator norms of $\fM$ and $\fM^{[1]}$ are the same. 
%Since the $L^p$ operator norm remains unchanged under rotations in $\ux$ or $\uy$, the $L^p$ boundedness of $\fM$ is equivalent to the boundedness of $\fM^{[1]}$. 
Now, we decompose $\ox = (\wt x, \breve x)\in \R^n \times \R^{m-n}$. Recall that $L_i = 0$ if $n<i\le m$, and let  $\vec L=(L_1,\dots, L_n)$. Then 
\[
\fM^{[1]} f(x) = \sup_{t > 0} \Big|\int_\Sigma f(\ux - t\om, \wt x - t(\ux, t)^\intercal \vec L  \om, \breve x)d\mu(\om)\Big|.
\]
%Slightly abusing notation, we denote $\vec L = (L_1,\dots, L_n)$. 
If we define 
%
%This identity gives that, if we denote %$f^{\breve x}(\uy, \wt y) = f(\uy, \wt y, \breve x)$ and
\[
\fM^{[2]} f(\ux, \wt x) = \sup_{t > 0} \Big|\int_\Sigma f(\ux - t\om, \wt x - t(\ux, t)^\intercal \vec L \om)d\mu(\om)\Big|,
\]
then % by a slicing argument $\|\fM^{[1]} f\|_{L^p(G)}^p \le 
\[
\|\fM^{[1]} f\|_{L^p(G)}^p\le \|\fM^{[2]}\|_{L^p(\R^{d+n})\to L^p(\R^{d+n})}^p \int |f(\uy,\wt y, \breve x)|^p d\uy\,d\wt y\,d\breve x.
\]
Thus, Theorem \ref{thm:main-hsu} follows once we show that $\fM^{[2]}$ is bounded on $L^p(\R^{d+n})$ for  $p > (r+1)/r$. To verify that $r$ and the {\L}ojasiewicz exponent remain unchanged, we note that 
\begin{align}\label{i:aS=bL}
    \sum_{i=1}^{m} a_i S_i = \sum_{i=1}^m b_i L_i,\quad b = U^\intercal a \in \R^m.
\end{align}
We write $L_i = \begin{pmatrix}
    \wt J_i \\ \wt \la_i^\intercal
\end{pmatrix}$ for $1\le i\le m$. Let $\wt \La$ be an $m\times d$ matrix with rows $\wt\la_1^\intercal,\dots, \wt\la_m^\intercal$. For $1\le k\le m$, let $P_k:\R^m\to \R^k$ be the projection onto the first $k$ coordinates. From now on, $b$ denotes the first $n$ coordinates of $U^\intercal a$. By \eqref{i:aS=bL}, we have
\begin{align*}
    \dist(a, \ker(\La^\intercal)) &= \dist(U^\intercal a, \ker(\wt\La^\intercal)) = \dist (P_n U^\intercal a, P_n(\ker(\wt\La^\intercal))), \\
    \det((S^a)^\intercal S^a -\la I_d) & = \det((L^b)^\intercal L^b - \la I_d)
\end{align*}
for $a\in \R^m$, where $L^b = \sum_{i=1}^n b_i L_i$. The identities above imply that
\begin{align*}
    \cC_r(a) = \wt \cC_r(b),\quad d_\perp(a) = \wt d_\perp(b),\quad a\in \R^m
\end{align*}
with $b = P_n U^\intercal a$, where $\wt \cC_r$ denotes the coefficient of $\la^{d-r}$ in the characteristic polynomial $\det((L^b)^\intercal L^b - \la I_d)$ and $\wt d_\perp(b) = \dist (P_n U^\intercal a, P_n(\ker(\wt\La^\intercal)))$. Moreover,
\[
\min_{\substack{a \in \R^m\setminus\{0\} \\ \La^\intercal a \neq 0}}\,\rank\big(\sum_{i=1}^m a_i S_i\big) = \min_{\substack{b \in \R^n\setminus\{0\}\\\sum_{i=1}^n b_i\wt\la_i\neq 0}}\,\rank\big(\sum_{i=1}^n b_i L_i\big).
\]
Thus, $r$ and the {\L}ojasiewicz exponent are invariant under the change from $S_1,\dots, S_m$ to $L_1,\dots, L_n$.

Next, we argue that the original submanifold $\Sigma$ can be replaced with the graph of a function $g$ satisfying \eqref{r:propertyofg}.
Pick $\om_\circ\in \Sigma$. Let $\chi^{\om_\circ}$ be a smooth cutoff function supported in a $|\om_\circ|\delta$-neighborhood of $\om_\circ$. By a compactness argument, it suffices to prove $L^p$-boundedness for 
%the $L^p$ boundedness for the maximal operator generated by the convolutions  $f\ast (\chi^{\om_\circ}\mu)_t$ implies the boundedness of  $\fM^{[2]}$. Equivalently, it is sufficient to consider 
$\sup_{t>0}|f\ast (\chi^{\om_\circ}\mu)_t|$. 

Since the origin lies in the domain enclosed by $\Sigma$, we have $\om_\circ\neq 0$. Let $e_{\circ,1} = \om_\circ/|\om_\circ|$. Choose orthonormal vectors $e_{\circ,2}, \dots, e_{\circ, d}$ so that $(e_{\circ, i})_{i=1}^d$ forms a basis of $\R^d$. By the assumption that $\Sigma$ is convex with nonvanishing Gaussian curvature, $|\om_\circ|^{-1}\Sigma$ is parametrized near $e_{\circ, 1}$ by 
\[
\Gamma(\om') = G(\om') e_{\circ,1} + \sum_{i=2}^d \om_i e_{\circ, i},\quad (\om')^\intercal = (\om_2,\dots, \om_d)\in \R^{d-1},
\]
where $G\in C^\infty(\R^{d-1})$ such that $G(0) = 1$, $G'(0) = b$ with $b\in \R^{d-1}$, and $G''(0)$ is negative-definite. Using this parametrization, we write
\begin{align*}
    f\ast (\chi^{\om_\circ}\mu)_t = \int f(\ux - t|\om_\circ|\,\Gamma(\om'), \ox - t|\om_\circ|(\ux, t)^\intercal \vec L\,\Gamma(\om')) \varphi_0(\om') d\om',
\end{align*}
where $\varphi_0(\om') := A_\circ^{-1}|\om_\circ|^{d-1}(1 + |G'(\om')|^2)^{1/2}\chi^{\om_\circ}(|\om_\circ|\Gamma(\om'))$. Here, $A_\circ$ denotes the area of $\Sigma$. This constant arises because $d\mu$ is a normalized measure.

Let $R$ be a rotation such that $R e_{\circ,i} = e_i$ for $1\le i\le d$. Then $R\Gamma(\om') = \Gamma_G(\om') := G(\om')e_1 + \sum_{i=2}^d \om_i e_i$. Let $\cR$ be the $(d+1)\times (d+1)$ matrix given by 
\begin{align*}
    \cR = \begin{pmatrix}
    |\om_\circ|^2 R & 0 \\ 0 & |\om_\circ|
\end{pmatrix}.
\end{align*}
Denote $\vec S_R = (S_{R,1},\dots,S_{R,n}) = (\cR L_1 R^\intercal,\dots, \cR L_n R^\intercal)$.
Now we define
\begin{align*}
    \cA^{[3]}_t f(x) =  \int f(\ux - t\,\Gamma_G(\om'), \ox - t(\ux, t)^\intercal \vec S_R\,\Gamma_G(\om')) \varphi_0(\om') d\om'.
\end{align*}
A calculation gives
\[
f\ast (\chi^{\om_\circ}\mu)_t(x) = \cA^{[3]}_t f_R(|\om_\circ|^{-1}R\ux, \ox),\quad f_R(\uy,\oy) := f(|\om_\circ|R^\intercal \uy,\oy).
\]
By the change of variables $(\ux,\ox) \to (|\om_\circ|R^\intercal \ux,\ox)$, $(\uy,\oy)\to (|\om_\circ|^{-1}R\uy, \oy)$, it follows that the $L^p$ operator norms  of the maximal operators corresponding to $f\ast (\chi^{\om_\circ}\mu)_t$, $\cA^{[3]}_t f$ are comparable.
Hence the problem is reduced to estimating $\fM^{[3]}$ defined by $\fM^{[3]} f := \sup_{t>0}|\cA^{[3]}_t f|$.
%\as{del (formula is not necessary -  delete: equal up to a constant, that is,
%\[
%\sup_{\|f\|_2 = 1}\big\|\sup_{t>0}\big|f\ast (\chi^{\om_\circ}\mu)_t\big|\big\|_{L^p(\R^{d+n})} = C \sup_{\|f\|_2 = 1} \big\|\sup_{t>0}\big|\cA^{[3]}_t f\big|\big\|_{L^p(\R^{d+n})}
%\]
%with $C>0$.}

Moreover, we note that the rank $r$ and the {\L}ojasiewicz exponent remain invariant under the change from $\vec L$ to $\vec S_R$. Indeed, the products of the largest $r$ singular values before and after this change of coordinates are comparable. Moreover, the transpose of the tilt matrix has the same kernel. Thus both $r$ and the Łojasiewicz exponent are preserved.

To end the proof, we reduce to the situation where $G'(0) = 0$. To do so, let $g(\om') = G(\om') - b^\intercal \om'$. Note that $g(0)=1$, $g'(0)=0$ and $g''(0)$ is negative definite. % satisfies the relation \eqref{r:propertyofg}. 
Denote $S_{R,i} = \begin{pmatrix}
    J_{R,i} \\ \la_{R,i}^\intercal
\end{pmatrix}$ for $1\le i\le n$. Also, by $\bar \La$ we denote an $n\times d$ matrix with rows $\la_{R,1}^\intercal,\dots, \la_{R,n}^\intercal$. 

Let $P:\bbR^d \to \bbR^{d-1}$ be the projection with $Pe_1=0$, $Pe_j=e_{j-1}$ for $2\le j\le d$. We now set
\begin{align*}
    \mathrm J_{b,i} = J_{R,i} + P^\intercal b e_1^\intercal J_{R,i} P^\intercal P - P^\intercal P J_{R,i}^\intercal e_1 b^\intercal P,\quad \La_b = \bar\La + \bar\La e_1 b^\intercal P
\end{align*}
for $1\le i\le n$, which, using the skew-symmetry of $J_{R,i}$, can be expressed as
\Be\label{eq:changeofvarS} 
\mathrm J_{b,i} = \begin{pmatrix} 1&0\\ b& I_{d-1}\end{pmatrix} J_{R,i} \begin{pmatrix} 1&b^\intercal\\0&I_{d-1} \end{pmatrix}, \qquad  \La_b = \bar\La \begin{pmatrix} 1&b^\intercal\\0&I_{d-1} \end{pmatrix}.
\Ee

Furthermore, we set $\Gamma_g (\om') = g(\om') e_1 + \sum_{i=2}^d \om_i e_i$.
A calculation gives
\begin{align*}
    x_1 - tG(\om') &= x_1 - b^\intercal x' - tg(\om') + b^\intercal (x'-t\om'), \\
    \ux^\intercal J_{R,i} \Gamma_G(\om') &= \big((x_1 - b^\intercal x')e_1 + P^\intercal x'\big)^\intercal \mathrm J_{b,i} \Gamma_g(\om'), \\
    \bar\La \Gamma_G(\om') &= \La_b \Gamma_g(\om')
\end{align*}
for $1\le i\le n$. Here, we denote $\ux = (x_1, x')\in \R\times \R^{d-1}$. For more details (in the case $\La=0$) we refer to  \cite[Section 2]{RyuSeeger}. 

Using the above identities, we obtain
\begin{align}\label{i:comparisona1a2}
    \cA^{[3]}_t f(\ux, \ox) = \cA^{[4]}_t f_b(x_1 - b^\intercal x', x', \ox),
\end{align}
where
\begin{align*}
    \cA^{[4]}_t f(x) = \int f(\ux - t\,\Gamma_g(\om'), \ox - t(\ux, t)^\intercal \vec {\mathrm S}_b\,\Gamma_g(\om')) \varphi_0(\om') d\om'
\end{align*}
and $f_b(y_1, y', \oy) = f(y_1 + b^\intercal y', y', \oy)$. Here 
$ (\ux, t)^\intercal \vec {\mathrm S}_b\,\Gamma_g(\om') \varphi_0(\om')$ is short for $ \sum_{i=1}^n e_{d+i} 
(\ux, t)^\intercal  {\mathrm S}_{b,i} \,\Gamma_g(\om') \varphi_0(\om')
$ 
%$\vec S_{b} = (\mathrm S_{b,1},\dots, \mathrm S_{b,n})$, 
where the $\mathrm S_{b,i}$ are $(d+1)\times d$ matrices defined by 
\[
\mathrm S_{b,i} =  \begin{pmatrix}
   \mathrm  J_{b,i} \\ \la_{b,i}^\intercal
\end{pmatrix} \text{ with } \la_{b,i}^\intercal = \la_{R,i}^\intercal + (e_1^\intercal \la_{R,i}) b^\intercal P ,
 \quad 1\le i\le n.
\]
%\as{changed  from the previous  {\color {green} $\mathrm S_{b,i} = \begin{pmatrix}    \mathrm  J_{b,i} \\ e_1^\intercal \La_b^\intercal \end{pmatrix}$.}} 
%\jr{changed: $\la_i^\intercal + (e_1^\intercal \lambda) b P^\intercal \to \la_i^\intercal + (e_1^\intercal \lambda_i) b^\intercal P$ because $b^\intercal P$ is a row vector\cas} 
From \eqref{eq:changeofvarS} we get 
\begin{align*}
    \mathrm S_{b,i} = \begin{pmatrix}
        1 & 0 & 0 \\
        b & I_{d-1} & 0 \\
        0 & 0 & 1
    \end{pmatrix} S_{R,i} \begin{pmatrix}
        1 & b^\intercal \\
        0 & I_{d-1}
    \end{pmatrix}
\end{align*} for $i=1,\dots,n$,  from which we conclude that $r$ and the {\L}ojasiewicz exponent are preserved under the change from $\vec S_R$ to $ \vec {\mathrm S}_b$.
%immediately follows because the square matrices multiplied by $S_{R,i}$ on both sides are invertible.

Consequently, using  \eqref{i:comparisona1a2}, we see that the $L^p(\bbR^{d+n})$  boundedness of  $\fM^{[3]}$ 
%would be bounded on $L^p(\R^{d+n})$ 
for $p> (r+1)/r$ follows from the $L^p(\bbR^{d+n})$  boundedness of  $\fM^{[4]}$ in the same range, which in turn is equivalent to \eqref{e:ineqforAt}. Applying Theorem \ref{thm:main2} to $\cA_t^{[4]}$ proves the sufficiency assertion in Theorem \ref{thm:main-hsu}.
%if the maximal operator $\fM^{[4]}$ defined by $\fM^{[4]} f := \sup_{t>0}|\cA^{[4]}_t f|$ is bounded on $L^p(\R^{d+n})$ whenever $p> (r+1)/r$. The boundedness of $\fM^{[4]}$, however, is equivalent to \eqref{e:ineqforAt}. Hence \textit{ii)} in Theorem \ref{thm:main} follows.

\section{Reduction to oscillatory integral estimates}
\label{app:osc-reduction}

In this appendix, we give the details of the reduction of the estimates in Proposition \ref{prop:locall2} to the oscillatory integral bounds \eqref{e:claim3}--\eqref{e:claim4} using localization, a shear transformation, and Fourier transforms.

We may assume that $f$ is supported in $B(0,\delta/2)$ via a localization argument exploiting the group translation on $G$; see \cite[Section 2]{RoosSeegerSrivastava-imrn} for details. To proceed, we consider the two cases $r^\sharp = r$ and $r^\sharp < r$ separately. First, assume $r^\sharp = r$. In this case, the problem is to estimate $f\ast \nu_t^{k,l}$.
We set
\begin{align*}
    a^{[1]}(x,t,y) &= (2\pi)^{-m-1} \varphi_1(\frac{x_1 - y_1}{t}, \frac{\ox-\oy +\ux^\intercal \vec J\uy}{t^2})\varphi_0(\frac{x'-y'}{t})\chi_\delta(y), \\ 
    \psi^{[1]}(x,t,y,u_1, \overline u) &= u_1(x_1-y_1-tg(\frac{x'-y'}{t})) + \inn{\overline u}{\ox-t\La \ux-\oy+(\ux,t)^\intercal \vec S\oy}
\end{align*}
where $\chi_\delta\in C_c^\infty(\R^{d+m})$ is chosen such that $\supp(\chi_\delta)\subset B(0,\delta)$ and $\chi_\delta\equiv 1$ on $B(0,\delta/2)$.
With these notations, $f\ast \nu_t^{k,l} = \int K_t^{k,l}(x,y) f(y) dy$, where
\begin{align*}
    K_t^{k,l}(x,y) := t^{-d-2m} a^{[1]}(x,t,y) \int \zeta_{k,l}(\sigma, \tau) e^{i\psi^{[1]}(x,t,y,t^{-1}\sigma, t^{-2}\tau)} d\sigma d\tau.
\end{align*}
Making a change of variables $\upsilon_1 = 2^{-k}t^{-1}\sigma$, $\overline\upsilon = 2^{-k}t^{-2}\tau$ gives
\begin{align*}
    K_t^{k,l}(x,y) = 2^{k(m+1)} a^{[1]}(x,t,y) \int Z_{k,l}^{[1]}(t,\upsilon_1, \overline\upsilon) e^{i2^k\psi^{[1]}(x,t,y,\upsilon_1, \overline\upsilon)} d\upsilon_1 d\overline\upsilon,
\end{align*}
where
\begin{align*}
    Z_{k,l}^{[1]}(t,\upsilon_1, \overline\upsilon) = \begin{cases}
        t^{-d+1} \zeta_1(2^{l+1}t\upsilon_1)\zeta_1(2(t^2\upsilon_1^2 + t^4|\overline\upsilon|^2)^{1/2}), & l<k, \\
        t^{-d+1} \zeta_0(2^lt\upsilon_1)\zeta_1(2(t^2\upsilon_1^2 + t^4|\overline\upsilon|^2)^{1/2}), & l=k.
    \end{cases}
\end{align*}
We now make use of the shear transformation
\[
(u_1,\overline u) = (\upsilon_1 - \inn{\overline\upsilon}{(\ux,t)^\intercal \vec S e_1},\overline\upsilon)
\]
to rewrite the kernel $K_t^{k,l}$. Let
\begin{align*}
    \psi^{[2]}(x,t,y,u_1, \overline u) :=\, &u_1(x_1-y_1-tg(\frac{x'-y'}{t})) \\
    &+ \inn{\overline u}{\ox-\oy + (\ux,t)^\intercal \vec S\big( P^\intercal y'- tg(\frac{x'-y'}{t})e_1\big)}.
\end{align*}
Performing a calculation, we observe
\[
\psi^{[1]}(x,t,y,\upsilon_1, \overline\upsilon) = \psi^{[2]}(\ux, \ox - t\La \ux + x_1 (\ux,t)^\intercal \vec S e_1, t,y,\upsilon_1 - \inn{\overline\upsilon}{(\ux,t)^\intercal \vec S e_1}, \overline\upsilon).
\]
Thus, making the shear transformation $(\upsilon_1, \overline\upsilon) \mapsto (u_1, \overline u)$ yields
\begin{align}\label{i:relKandcK}
    K^{k,l}_t(x,y) = \cK_t^{k,l}(\ux, \ox - t\La \ux + x_1 (\ux,t)^\intercal \vec S e_1, y),
\end{align}
where
\begin{align*}
    \cK_t^{k,l}(x, y) &= 2^{k(m+1)} a^{[1]}(x + t\La \ux - x_1 (\ux,t)^\intercal \vec S e_1,t,y) \\
    &\qquad \times\int Z_{k,l}^{[1]}(t,u_1 + \inn{\overline u}{(\ux,t)^\intercal \vec S e_1}, \overline u) e^{i2^k\psi^{[2]}(x,t,y,u_1,\overline u)} du_1 d\overline u.
\end{align*}
Set
\begin{align*}
    a^{[2]}(x,t,y) &=  a^{[1]}(x + t\La \ux - x_1 (\ux,t)^\intercal \vec S e_1,t,y), \\
    Z_{k,l}^{[2]}(x,t,u_1,\overline u) &= Z_{k,l}^{[1]}(t,u_1 + \inn{\overline u}{(\ux,t)^\intercal \vec S e_1}, \overline u).
\end{align*}
In addition, we define $\cT^{k,l}$ by $\cT^{k,l} f(x,t) = \int \cK^{k,l}_t(x,y) f(y) dy$.
By the above discussion and \eqref{i:relKandcK}, the estimate for $f\ast \nu_t^{k,l}$ in Proposition \ref{prop:locall2} follows once we prove
\begin{align}\label{e:locall2cT}
    \|\cT^{k,l} f\|_{L^2(\R^{d+m}\times [1,2])}\lesssim 2^{-\frac{d-2}{2}k}2^{\frac{d-r-2}{2}l}\|f\|_{L^2(\R^{d+m})}.
\end{align}
Furthermore, one can also see that proving \eqref{e:locall2cT} essentially implies the estimate for the derivative (in $t$) of $f\ast \nu_t^{k,l}$ in the proposition. This is a consequence of straightforward calculation, which we omit. We refer the reader to \cite[Section 5]{RyuSeeger} for the details.

Using the same argument, in the case $r^\sharp < r$ we can deduce Proposition \ref{prop:locall2} $(ii)$ from an analogous estimate to \eqref{e:locall2cT}. Indeed, let
\begin{align*}
    Z_{k,l,j}^{[1]}(t,\upsilon_1, \overline\upsilon) = \begin{cases}
        Z_{k,l}^{[1]}(t,\upsilon_1, \overline\upsilon) \zeta_{k-j}(2^k t^2|\overline{\upsilon}^{(2)}|), & j < l, \\
        Z_{k,l}^{[1]}(t,\upsilon_1, \overline\upsilon)\sum_{j=l}^k \zeta_{k-j}(2^k t^2|\overline{\upsilon}^{(2)}|), & j=l,
    \end{cases}
\end{align*}
where we denote $\overline \upsilon = (\overline \upsilon^{(1)}, \overline \upsilon^{(2)})\in \R^{m-m^\sharp}\times \R^{m^\sharp}$.
We define $\cK_t^{k,l,j}$ and $Z_{k,l,j}^{[2]}$ as $\cK_t^{k,l}$ and $Z_{k,l}^{[2]}$ but with $Z_{k,l}^{[1]}$ in their definitions replaced by $Z_{k,l,j}^{[1]}$. Replacing $Z^{[1]}_{k,l}$ by $Z^{[1]}_{k,l,j}$ leaves the shear transformation unchanged. Thus Proposition \ref{prop:locall2} $(ii)$ reduces to
\begin{align}\label{e:locall2cT-2}
    \|\cT^{k,l,j} f\|_{L^2(\R^{d+m}\times [1,2])}\lesssim 2^{-\frac{d-2}{2}k}2^{\frac{d-r-2}{2}l} (\mathbf B_{j,l}(r^\sharp,r))^{\frac12}\|f\|_{L^2(\R^{d+m})}
\end{align}
is valid, where $\cT^{k,l,j}$ is an operator on $L^2(\R^{d+m})$ given by $\cT^{k,l,j} f(x,t) = \int \cK_t^{k,l,j}(x,y) f(y) dy$.

To reduce \eqref{e:locall2cT} and \eqref{e:locall2cT-2} to \eqref{e:locall2-2}, we claim that $a^{[2]}$ can be assumed to be a tensor product expressed as
\begin{align}\label{i:tensorpd}
    a^{[2]}(x,t,y) = \gamma(t)\prod_{i=1}^{d+m} \alpha_i(x_i) \prod_{j=1}^{d+m} \beta_j(y_j).
\end{align}
To see this, using the Fourier series expansion, we write
\begin{align*}
    a^{[2]}(x,t,y) = \sum_{n\in \Z}&\sum_{(\iota^1,\iota^2)\in \Z^{d+m}\times \Z^{d+m}} c_{n,\iota_1, \iota_2} e^{i(nt + \inn{\iota_1}{x} + \inn{\iota_2}{y})} \\
&\hspace{40mm} \times \gamma_n(t)\prod_{i=1}^{d+m} \alpha_{\iota_1,i}(x_i) \prod_{j=1}^{d+m} \beta_{\iota_2, j}(y_j).
\end{align*}
Here, $c_{n,\iota_1,\iota_2}\in \C$ are constants such that $c_{n,\iota_1,\iota_2} = O((1+|n| + |\iota_1| + |\iota_2|)^{-N})$, and $\gamma_n$, $\alpha_{\iota_1,i}$, $\beta_{\iota_2, j}$ are compactly supported smooth functions. Via the Fourier series expansion of $a^{[2]}$ above, we split $\cT^{k,l}$ and $\cT^{k,l,j}$ into the sums
\begin{align*}
    \cT^{k,l} f(x,t) &= \sum_{n\in \Z}\sum_{(\iota^1,\iota^2)\in \Z^{d+m}\times \Z^{d+m}} \cT^{k,l}_{n,\iota_1,\iota_2} f(x,t), \\
    \cT^{k,l,j} f(x,t) &= \sum_{n\in \Z}\sum_{(\iota^1,\iota^2)\in \Z^{d+m}\times \Z^{d+m}} \cT^{k,l,j}_{n,\iota_1,\iota_2} f(x,t),
\end{align*}
where $\cT^{k,l}_{n,\iota_1,\iota_2}$ and $\cT^{k,l,j}_{n,\iota_1,\iota_2}$ are defined similarly to $\cT^{k,l}$, $\cT^{k,l,j}$, respectively, but with $a^{[2]}$ replaced by the summands from the Fourier expansion of $a^{[2]}$ above. Note that inserting the modulation term $e^{i(nt + \inn{\iota_1}{x} + \inn{\iota_2}{y})}$ does not affect an operator norm. Thus, by the triangle inequality and bounds of $c_{n,\iota_1,\iota_2}$ given above, \eqref{e:locall2cT} (\textit{resp.,} \eqref{e:locall2cT-2}) follows once we prove \eqref{e:locall2cT} with $\cT^{k,l}$ (\textit{resp.,} \eqref{e:locall2cT-2} with $\cT^{k,l,j}$) replaced by a single component $\cT^{k,l}_{n,\iota_1,\iota_2}$ (\textit{resp.,} $\cT^{k,l,j}_{n,\iota_1,\iota_2}$). This justifies our claim.

We focus on the estimation of $\cT^{k,l,j}$ in \eqref{e:locall2cT-2} corresponding to the case $r^\sharp < r$. The remaining case $r^\sharp = r$ will be discussed later. Keeping \eqref{i:tensorpd} in mind, we define
\begin{align*}
    \cF^m_{k,j} G(\ux,\ox,t) &= \prod_{i=d+1}^{d+m} \alpha_i(x_i) \int_{\R^m} \zeta_{k-j}(2^k t^2|\overline{u}^{(2)}|) G(\ux,\overline u,t) e^{i2^k\inn{\ox}{\overline u}} d\overline u, \\
    \cF^{m+1}_{k} f(u_1,y',\overline u) &= \int_{\R^{m+1}} e^{-i2^k(y_1 u_1 + \inn{\oy}{\overline u})} \beta_1(y_1)\prod_{i=d+1}^{d+m} \beta_i(y_i) f(y_1,y',\oy) dy_1 d\oy
\end{align*}
for $j<l$, and
\begin{align*}
        \cG^{k,l} F(x,\overline u, t) &= \gamma(t)\prod_{i=1}^{d} \alpha_i(x_i) \\
        &\times \int_{\R^d} e^{i2^k\Psi^{\overline u}(x,t,u_1,y')} Z_{k,l}^{[2]}(x,t,u_1, \overline u) \prod_{i=2}^{d} \beta_i(y_i) F(u_1, y', \overline u) du_1 dy',
    \end{align*}
    with
    \begin{align*}
        \Psi^{\overline u}(x,t,u_1,y')&= u_1(x_1-tg(\tfrac{x'-y'}{t}) )+ (x,t)^\intercal S^{\overline u} (P^\intercal y'-tg(\tfrac{x'-y'}{t})e_1).
    \end{align*}
    Here, $S^{\overline u} = \sum_{i=1}^{m} \overline u_i S_i$. When $j = l$, we replace $\zeta_{k-j}$ in the definition of $\cF_{k,j}^m$ with $\sum_{j=l}^k \zeta_{k-j}$.
Then, $\cT^{k,l,j}$ is expressed as
\[
\cT^{k,l,j} f(x,t) = 2^{k(m+1)} \cF^m_{k,j} \cG^{k,l}\cF^{m+1}_k f(x,t).
\]
By scaling and Plancherel's theorem, we have
\begin{align*}
    \|\cF^m_{k,j}\|_{L^2(\R^{d+m+1})\to L^2(\R^{d+m+1})}&\lesssim 2^{-\frac{m}{2}k}, \\
    \|\cF^{m+1}_k\|_{L^2(\R^{d+m})\to L^2(\R^{d+m})}&\lesssim 2^{-\frac{m+1}{2}k}.
\end{align*}
We now set $\cG^{k,l}_{\overline u} f(x,t) = \cG^{k,l} f(x,\overline u,t)$. The preceding bounds show that \eqref{e:locall2cT-2} follows from \eqref{e:claim3} and \eqref{e:claim3-b}.

We now turn to the case $r = r^\sharp$. For $r = r^\sharp$, the same argument applies with $\cT^{k,l} = 2^{k(m+1)} \cF^m_k \cG^{k,l}\cF^{m+1}_k$, where
\[
\cF^m_k G(\ux,\ox,t) = \prod_{i=d+1}^{d+m} \alpha_i(x_i) \int_{\R^m} G(\ux,\overline u,t) e^{i2^k\inn{\ox}{\overline u}} d\overline u.
\]
Note that $\|\cF^m_k\|_{L^2(\R^{d+m+1})\to L^2(\R^{d+m+1})} \lesssim 2^{-\frac m2 k}$. Following the above argument, we conclude that \eqref{e:claim4} implies \eqref{e:locall2cT}, provided that $2^{-l}$ is small. If $2^{-l}\sim 1$, a similar conclusion can be derived, except that $|\overline u|\sim 1$ above is replaced by $|\overline u|\lesssim 1$. %See \eqref{e:Gu2jlarge} below.

We now define
\[
\chi_{\overline u}(x,t,y) = t^{-d+1}\gamma(t) \prod_{i=1}^d \alpha_i(x_i) \prod_{i=2}^d \beta_i(y_i) \zeta_1(2(t^2\sigma_{\overline u}(x',t,y_1)^2 + t^4|\overline u|^2)^{1/2}).
\]
By the construction of $\varphi_0, \varphi_1$, the support of $\chi_{\overline u}$ is contained in the set $\{(x,t,y): t\sim 1,\ |x'|,\ |y'|\ll 1,\ |x_1|,\ |y_1|\lesssim 1\}$ for all $|\overline u|\sim 1$, with the implicit constants independent of $\overline u$. Thus, decomposing the cutoff functions $\gamma, \alpha_i, \beta_i$ into finite pieces and using translations if necessary, we may assume that $\chi_{\overline u}$ satisfies the support condition \eqref{r:chisupp}.

\subsection*{Acknowledgement}
This work was supported by the National Research Foundation of Korea (NRF) grant funded by the Korea government (MSIT) (No. RS-2026-25483775). The author thanks Andreas Seeger for their valuable comments and helpful discussions.

\subsection*{Use of AI tools} 
The author used GPT-$6$ Astra to assist in proofreading the manuscript, refining the exposition, and searching for relevant references (\cite{BDDHP17}). The author formulated the new results in this paper and developed the main ideas and structure of their proofs independently of generative AI tools.

\end{document}